\documentclass[10pt,journal,compsoc]{IEEEtran} 
\usepackage{floatrow}
\newfloatcommand{capbtabbox}{table}[][\FBwidth]
\usepackage{times}
\usepackage{epsfig}
\usepackage{graphicx}
\usepackage{amsmath}
\usepackage{amssymb}
\usepackage[english]{babel}

\usepackage[utf8x]{inputenc}
\usepackage[english]{babel}
\usepackage[OT1]{fontenc}
\usepackage{amssymb}
\usepackage{amsmath}
\usepackage{amsthm}
\usepackage{subfig}
\usepackage{epsfig}
\usepackage{epstopdf}
\usepackage{textcomp,booktabs}
\usepackage{colortbl}

\usepackage{graphicx}
\usepackage{algorithm}
\usepackage{algorithmic}
\usepackage{amsfonts}
\usepackage{fancyhdr}
\usepackage{dsfont}
\usepackage{color}

\usepackage{bm}
\usepackage{afterpage}
\usepackage{multirow}
\usepackage{diagbox}
\usepackage[font=small,skip=0pt]{caption}
\definecolor{colortwo}{RGB}{255,113,31}
\definecolor{colorone}{RGB}{218,218,217}

\usepackage[utf8x]{inputenc}
\usepackage[english]{babel}
\usepackage[OT1]{fontenc}
\usepackage{amssymb}
\usepackage{amsmath}
\usepackage{amsthm}
\usepackage{subfig}
\usepackage{epsfig}
\usepackage{epstopdf}
\usepackage{xcolor,colortbl}
\usepackage{graphicx}
\usepackage{algorithm}
\usepackage{algorithmic}
\usepackage{amsfonts}
\usepackage{fancyhdr}
\usepackage{dsfont}
\usepackage{color}
\usepackage{bbm}
\usepackage{expl3}
\usepackage[breaklinks=true,bookmarks=false]{hyperref}

\newcommand{\bbb}[1]{\boldsymbol{\mathbf{#1}}}
\def\cone{\textcolor[rgb]{0.9961,0,0}}
\def\ctwo{\textcolor[rgb]{0,0.5820,0.9}}
\def\cthree{\textcolor[rgb]{0,0.7,0}}
\definecolor{mygray}{gray}{.9}
\definecolor{mypink}{rgb}{.99,.91,.95}

\definecolor{mycyana}{RGB}{153,204,153}
\definecolor{mycyanb}{RGB}{255,153,102}
\definecolor{mycyanc}{RGB}{255,204,51}
\definecolor{mycyand}{RGB}{204,204,102}

\newtheorem{lemma}{Lemma}
\newtheorem{theorem}{Theorem}

\newtheorem{proposition}{Proposition}
\newtheorem{assumption}{Assumption}

\def\beq{\begin{eqnarray}}
\def\eeq{\end{eqnarray}}

\def\sbeq{\begin{equation}\begin{split}}
\def\seeq{\end{split}\end{equation}}

\def\noi{\noindent}
\def\nn{\nonumber}
\def\la{\langle}
\def\ra{\rangle}

\newcommand{\onlyonefont}{\fontsize{8.7}{8.7}\selectfont}

\def\P{\mathcal{T}}
\def\sign{\rm sign}
\newcommand{\algsize}{\fontsize{9.5}{12}\selectfont}

\ifCLASSOPTIONcompsoc
  \usepackage[nocompress]{cite}
\else
  \usepackage{cite}
\fi

\newcolumntype{M}[1]{>{\centering\arraybackslash}m{#1}}
\newcolumntype{N}{@{}m{0pt}@{}}
\begin{document}
%
\title{A Generalized Matrix Splitting Algorithm}

\author{Ganzhao~Yuan,~Wei-Shi Zheng,~Li Shen,~Bernard~Ghanem

\IEEEcompsocitemizethanks{\IEEEcompsocthanksitem Ganzhao Yuan is with the School of Data and Computer Science, Sun Yat-sen
University (SYSU), Guangzhou, Guangdong 510275, China, and also with
Key Laboratory of Machine Intelligence and Advanced Computing, Ministry
of Education, Beijing 221143, China. E-mail: yuanganzhao@gmail.com.

}

\IEEEcompsocitemizethanks{\IEEEcompsocthanksitem

Wei-Shi Zheng is with the School of Data and Computer Science, Sun Yat-sen
University (SYSU), Guangzhou, Guangdong 510275, China, and also with
Key Laboratory of Machine Intelligence and Advanced Computing, Ministry
of Education, Beijing 221143, China. E-mail: wszheng@ieee.org.

}

\IEEEcompsocitemizethanks{\IEEEcompsocthanksitem Li Shen (mathshenli@gmail.com) is with Tencent AI Lab, Shenzhen, China.
}

\IEEEcompsocitemizethanks{\IEEEcompsocthanksitem Bernard Ghanem (bernard.ghanem@kaust.edu.sa) is with Visual Computing Center, King Abdullah University of Science and Technology (KAUST), Saudi Arabia. \protect\\
}

\thanks{Manuscript received 01 JAN. 2018; revised 01 JAN. 2018; accepted 01 JAN.
2018. Date of publication 01 JAN. 2018; date of current version 0.0000.\hfil\break
Recommended for acceptance by XXX XXX. \hfil\break
For information on obtaining reprints of this article, please send e-mail to: reprints@ieee.org, and reference the Digital Object Identifier below.\hfil\break
Digital Object Identifier no. 00.0000/JOURNAL.2018.0000000}}


\markboth{}
{Shell \MakeLowercase{\textit{et al.}}: Bare Demo of IEEEtran.cls for Computer Society Journals}

\IEEEtitleabstractindextext{%
\begin{abstract}
Composite function minimization captures a wide spectrum of applications in both computer vision and machine learning. It includes bound constrained optimization, $\ell_1$ norm regularized optimization, and $\ell_0$ norm regularized optimization as special cases. This paper proposes and analyzes a new Generalized Matrix Splitting Algorithm (GMSA) for minimizing composite functions. It can be viewed as a generalization of the classical Gauss-Seidel method and the Successive Over-Relaxation method for solving linear systems in the literature. Our algorithm is derived from a novel triangle operator mapping, which can be computed exactly using a new generalized Gaussian elimination procedure. We establish the global convergence, convergence rate, and iteration complexity of GMSA for convex problems. In addition, we also discuss several important extensions of GMSA. Finally, we validate the performance of our proposed method on three particular applications: nonnegative matrix factorization, $\ell_0$ norm regularized sparse coding, and $\ell_1$ norm regularized Dantzig selector problem. Extensive experiments show that our method achieves state-of-the-art performance in term of both efficiency and efficacy.


\end{abstract}

}

 \maketitle
\IEEEdisplaynontitleabstractindextext
\IEEEpeerreviewmaketitle

\renewcommand{\arraystretch}{1}
\section{Introduction}
In this paper, we focus on the following composite function minimization problem:
\beq \label{eq:main}
\min_{\bbb{x}}~f(\bbb{x}) \triangleq q(\bbb{x}) + h(\bbb{x});~q(\bbb{x}) = \tfrac{1}{2}\bbb{x}^T\bbb{A}\bbb{x} + \bbb{x}^T\bbb{b}
\eeq
\noi where $\bbb{x} \in\mathbb{R}^n,~\bbb{b} \in \mathbb{R}^n$, $\bbb{A}\in \mathbb{R}^{n\times n}$ is a symmetric positive semidefinite matrix, and $h(\bbb{x}): \mathbb{R}^{n}\mapsto \mathbb{R}$ is a piecewise separable function (\emph{i.e.}~$h(\bbb{x})=\sum_{i=1}^n h_i (\bbb{x}_i)$) but not necessarily convex. Typical examples of $h(\bbb{x})$ include the bound constrained function and the $\ell_0$ and $\ell_1$ norm functions. We assume that $f(\bbb{x})$ is bounded below for any feasible solution $\bbb{x}$.

The optimization in (\ref{eq:main}) is flexible enough to model a variety of applications of interest in both computer vision and machine learning, including compressive sensing \cite{Donoho06}, nonnegative matrix factorization \cite{lee1999learning,lin2007projected,guan2012nenmf}, sparse coding \cite{lee2006efficient,aharon2006img,BaoJQS16,Quan2016CVPR}, support vector machine \cite{hsieh2008dual}, logistic regression \cite{yu2011dual}, subspace clustering \cite{elhamifar2013sparse}, to name a few. Although we only focus on the quadratic function $q(\cdot)$, our method can be extended to handle general non-quadratic composite functions by considering a Newton approximation of the objective \cite{TsengY09,yuan2014newton} and to solve general linear constrained problems by using its associated augmented Lagrangian function of the problem \cite{HeY12,HeY15}.

The most popular method for solving problem (\ref{eq:main}) is perhaps the proximal gradient method \cite{nesterov2013introductory,beck2009fast}. It considers a fixed-point proximal iterative procedure $\bbb{x}^{k+1}= \text{prox}_{\gamma h}\left(\bbb{x}^k - \gamma \nabla q(\bbb{x}^k)\right)$ based on the current gradient $\nabla q(\bbb{x}^k)$. Here the proximal operator $\text{prox}_{\tilde{h}}(\bbb{a}) = \arg \min_{\bbb{x}} ~\frac{1}{2}\|\bbb{x}-\bbb{a}\|_2^2 + \tilde{h}(\bbb{x})$ can often be evaluated analytically, $\gamma={1}/{L}$ is the step size with $L$ being the local (or global) Lipschitz constant. It is guaranteed to decrease the objective at a rate of $\mathcal{O}({L}/{k})$, where $k$ is the iteration number. The accelerated proximal gradient method can further boost the rate to $\mathcal{O}({L}/{k^2})$. Tighter estimates of the local Lipschitz constant leads to better convergence rate, but it scarifies additional computation overhead to compute $L$. Our method is also a fixed-point iterative method, but it does not rely on a sparse eigenvalue solver or line search backtracking to compute such a Lipschitz constant, and it can exploit the specified structure of the quadratic Hessian matrix $\mathbf{A}$.

The proposed method is essentially a generalization of the classical Gauss-Seidel (GS) method and Successive Over-Relaxation (SOR) method \cite{demmel1997applied,saad2003iterative}. In numerical linear algebra, the Gauss-Seidel method, also known as the successive displacement method, is a fast iterative method for solving a linear system of equations. It works by solving a sequence of triangular matrix equations. The method of SOR is a variant of the GS method and it often leads to faster convergence. Similar iterative methods for solving linear systems include the Jacobi method and symmetric SOR. Our proposed method can solve versatile composite function minimization problems, while inheriting the efficiency of modern linear algebra techniques.


%

Our method is closely related to coordinate gradient descent and its variants such as randomized coordinate descent \cite{hsieh2008dual,patrascu2014iteration}, cyclic coordinate descent \cite{sun2015improved}, block coordinate descent \cite{nesterov2012efficiency,beck2013convergence,hong2013iteration}, randomized block coordinate descent \cite{richtarik2014iteration,lu2015complexity}, accelerated randomized coordinate descent \cite{nesterov2012efficiency,lin2015accelerated,lu2013randomized} and others \cite{liu2015asynchronous,johnson2013accelerating,zeng2015gauss}. However, all these work are based on gradient-descent type iterations and a constant Lipschitz step size. They work by solving a first-order majorization/surrogate function via closed form updates. Their algorithm design and convergence result cannot be applied here. In contrast, our method does not rely on computing the Lipschicz constant step size, yet it adopts a triangle matrix factorization strategy, where the triangle subproblem can be solved by an alternating cyclic coordinate strategy.


We are aware that matrix splitting algorithm has been considered to solve symmetric linear complementarity problems \cite{Luo1991,Luo1992,Iusem1993} and second-order cone complementarity problems \cite{zhang2014efficient} in the literature. However, we focus on minimizing a general separable nonsmooth composite function which is different from theirs. In addition, our algorithm is derived from a novel triangle operator mapping, which can be computed exactly using a new Gaussian elimination procedure. It is worthwhile to mention that matrix splitting has been extended to operator splitting recently to solve multi-term nonsmooth convex composite optimization problems \cite{ShenLYM17}.

\textbf{Contributions.} \textbf{(i)} We propose a new Generalized Matrix Splitting Algorithm (GMSA) for composite function minimization (See Section \ref{sect:alg}). Our method is derived from a novel triangle proximal operator (See Subsection \ref{sect:alg:subprob}). We establish the global convergence, convergence rate, and iteration complexity of GMSA for convex problems (See Subsection \ref{sect:alg:convergence}). \textbf{(ii)} We discuss several important extensions of GMSA (see Section \ref{sect:extension}). First, we consider a new correction strategy to achieve pointwise contraction for the proposed method (See Subsection \ref{sect:ext:strong}). Second, we discuss using Richardson extrapolation technique to further accelerate GMSA  (See Subsection \ref{sect:ext:acc}). Third, we extend GMSA to solve nonconvex problems with global convergence guarantee (See Subsection \ref{sect:ext:noncvx}). Fourth, we discuss how to adapt GMSA to minimize non-quadratic functions (See Subsection \ref{sect:ext:nonq}). Fifth, we show how to incorporate GMSA into the general optimization framework of Alternating Direction Method of Multipliers (ADMM) (See Subsection \ref{sect:ext:admm}). \textbf{(iii)} Our extensive experiments on nonnegative matrix factorization, sparse coding and Danzig selectors have shown that GMSA achieves state-of-the-art performance in term of both efficiency and efficacy (See Section \ref{sect:exp}). A preliminary version of this paper appeared in \cite{YuanZG17}.



\bbb{Notation.} We use lowercase and uppercase boldfaced letters to denote real vectors and matrices respectively. The Euclidean inner product between $\bbb{x}$ and $\bbb{y}$ is denoted by $\la \bbb{x},\bbb{y}\ra$ or $\bbb{x}^T\bbb{y}$. We denote $\|\bbb{x}\|=\|\bbb{x}\|_2=\sqrt{\la \bbb{x},\bbb{x} \ra}$, $\|\bbb{x}\|_{\bbb{A}} = \sqrt{\bbb{x}^T\bbb{A}\bbb{x}}$, and $\|\bbb{C}\|$ as the spectral norm (\emph{i.e.} the largest singular value) of $\bbb{C}$. We denote the $i^{\text{th}}$ element of vector $\bbb{x}$ as $\bbb{x}_i$ and the $(i,j)^{\text{th}}$ element of matrix $\mathbf{C}$ as $\mathbf{C}_{i,j}$. $diag(\bbb{D}) \in \mathbb{R}^{n}$ is a column vector formed from the main diagonal of $\bbb{D}\in\mathbb{R}^{n\times n}$. $\bbb{C}\succeq0$ and $\bbb{C}\succ0$ indicate that the matrix $\bbb{C}\in\mathbb{R}^{n\times n}$ is positive semidefinite and positive definite, respectively. Here $\bbb{C}$ is not necessarily symmetric \footnote{$\bbb{C}\succeq0 \Leftrightarrow \forall \bbb{x},~\bbb{x}^T\bbb{Cx}\geq 0 \Leftrightarrow \forall \bbb{x},~\frac{1}{2}\bbb{x}^T(\bbb{C}+\bbb{C}^T)\bbb{x}\geq 0$}. We denote $\bbb{D}$ as a diagonal matrix of $\bbb{A}$ and $\bbb{L}$ as a strictly lower triangle matrix of $\bbb{A}$ \footnote{For example, when $n=3$, $\bbb{D}$ and $\bbb{L}$ take the following form: \\$\textstyle \renewcommand\arraystretch{1}
\setlength{\arraycolsep}{3pt}
\bbb{D}=\begin{bmatrix}
  \bbb{A}_{1,1} & 0 & 0    \\
  0 & \bbb{A}_{2,2} & 0  \\
  0 & 0 & \bbb{A}_{3,3}   \\
\end{bmatrix},~\bbb{L}=\begin{bmatrix}
  0 & 0 & 0    \\
  \bbb{A}_{2,1} & 0 & 0  \\
  \bbb{A}_{3,1} & \bbb{A}_{3,2} &0    \\
\end{bmatrix}$}. Thus, we have $\bbb{A} = \bbb{L}+\bbb{D}+\bbb{L}^T$. Throughout this paper, $\bbb{x}^k$ denotes the value of $\bbb{x}$ at $k$-th iteration if $\bbb{x}\in\mathbb{R}^n$ is a variable, and ${x}^k$ denotes the $k$-th power of ${x}$ if ${x}\in\mathbb{R}$ is a constant scalar. We use $\bbb{x}^*$ to denote any solution of the optimal solution set of (\ref{eq:main}). For notation simplicity, we denote:
\begin{equation}\begin{split}
&~~~~~~~~~\bbb{r}^k\triangleq \bbb{x}^{k} - \bbb{x}^*,~\bbb{d}^k \triangleq \bbb{x}^{k+1}-\bbb{x}^k \nn\\
&u^k \triangleq f(\bbb{x}^{k})-f(\bbb{x}^*),~f^k \triangleq f(\bbb{x}^k),~f^* \triangleq f(\bbb{x}^*)\nn
\end{split}\end{equation}

\renewcommand{\arraystretch}{1}

\section{Proposed Algorithm} \label{sect:alg}

This section presents our proposed Generalized Matrix Splitting Algorithm (GMSA) for solving (\ref{eq:main}). Throughout this section, we assume that $h(\bbb{x})$ is convex and postpone the discussion for nonconvex $h(\bbb{x})$ to Section \ref{sect:ext:noncvx}.

Our solution algorithm is derived from a fixed-point iterative method based on the first-order optimal condition of (\ref{alg:main}). It is not hard to validate that a solution $\bbb{x}$ is the optimal solution of (\ref{alg:main}) if and only if $\bbb{x}$ satisfies the following nonlinear equation (``$\triangleq$'' means define):
\begin{equation}\begin{split} \label{eq:opt:cond}
\textstyle  \bbb{0} &\in \partial f(\bbb{x}) \\
&= \nabla q(\bbb{x}) + \partial h(\bbb{x})= \bbb{Ax} + \bbb{b} + \partial h(\bbb{x})
\end{split}\end{equation}

\noi where $\nabla q(\bbb{x})$ and $\partial h(\bbb{x})$ denote the gradient of $q(\cdot)$ and the sub-gradient of $h(\cdot)$ in $\bbb{x}$, respectively. In numerical analysis, a point $\bbb{x}$ is called a fixed point if it satisfies the equation $\bbb{x} \in \P(\bbb{x})$, for some operator $\P(\cdot)$. Converting the transcendental equation $\bbb{0}\in\partial f(\bbb{x})$ algebraically into the form $\bbb{x} \in \P(\bbb{x})$, we obtain the following iterative scheme with recursive relation:
\beq \label{eq:iterative}
\textstyle \bbb{x}^{k+1}\in \P(\bbb{x}^k),~k = 0, 1, 2,...
\eeq
We now discuss how to adapt our algorithm into the iterative scheme in (\ref{eq:iterative}). First, we split the matrix $\bbb{A}$ in (\ref{eq:opt:cond}) using the following strategy:
\beq\label{eq:matrix:dec}
\textstyle  \bbb{A} = \underbrace{\bbb{L}+\tfrac{1}{\omega}\bbb{D}+ \epsilon \bbb{I}}_{\bbb{B}} + \underbrace{ \bbb{L}^T + \tfrac{\omega-1}{\omega}\bbb{D} -  \epsilon \bbb{I}}_{\bbb{C}}
\eeq
\noi Here, $\omega \in (0,2)$ is a relaxation parameter and $\epsilon\in[0,\infty)$ is a parameter for strong convexity that enforces $diag(\bbb{B})>\bbb{0}$. These parameters are specified by the user beforehand. Using these notations, we obtain the following optimality condition which is equivalent to (\ref{eq:opt:cond}):
\beq \label{eq:kkt:BC}
\textstyle \textstyle-\bbb{Cx} -  \bbb{b} \in  (\bbb{B}+\partial h ) (\bbb{x})   \nn
\eeq
\noi Then, we have the following equivalent fixed-point equation:
\beq\label{eq:P}
\textstyle \bbb{x} \in \P(\bbb{x};\bbb{A},\bbb{b},h) \triangleq  (\bbb{B}+\partial h )^{-1}(-\bbb{Cx} -  \bbb{b})
\eeq
\noi For notation simplicity, we denote $\P(\bbb{x};\bbb{A},\bbb{b},h)$ as $\P(\bbb{x})$ since $\{\bbb{A},\bbb{b},h\}$ can be known from the context.

We name $\P$ the triangle proximal operator, which is novel in this paper\footnote{This is in contrast with Moreau's proximal operator \cite{parikh2014proximal}: $\text{prox}_{h}(\bbb{a}) = \arg\min_{\bbb{x}} ~\frac{1}{2}\|\bbb{x}-\bbb{a}\|_2^2 + h(\bbb{x})=(\bbb{I}+\partial h)^{-1}(\bbb{a})$, where the mapping $(\bbb{I}+\partial h)^{-1}$ is called the resolvent of the subdifferential operator $\partial h$.}.  Due to the triangle property of the matrix $\bbb{B}$ and the element-wise separable structure of $h(\cdot)$, the triangle proximal operator $\P(\bbb{x})$ in (\ref{eq:P}) can be computed exactly and analytically, by a generalized Gaussian elimination procedure (discussed later in Section \ref{sect:alg:subprob}). Our generalized matrix splitting algorithm iteratively applies $\bbb{x}^{k+1} \Leftarrow\P(\bbb{x}^k)$ until convergence. We summarize our algorithm in Algorithm \ref{alg:main}.


In what follows, we show how to compute $\P(\bbb{x})$ in (\ref{eq:P}) in Section \ref{sect:alg:subprob}, and then we study the convergence properties of Algorithm \ref{alg:main} in Section \ref{sect:alg:convergence}.




\subsection{Computing the Triangle Proximal Operator}\label{sect:alg:subprob}

We now present how to compute the triangle proximal operator in (\ref{eq:P}), which is based on a new generalized Gaussian elimination procedure. Notice that (\ref{eq:P}) seeks a solution $\bbb{z}^*\triangleq \P(\bbb{x}^k)$ that satisfies the following nonlinear system:
\beq \label{eq:subproblem}
\textstyle \bbb{0} \in \bbb{B}\bbb{z}^*   +  \bbb{u} +  \partial h(\bbb{z}^*)   ,~\text{where}~\bbb{u}=\bbb{b} + \bbb{C}\bbb{x}^k
\eeq
\noi By taking advantage of the triangular form of $\bbb{B}$ and the element-wise/decomposable structure of $h(\cdot)$, the elements of $\bbb{z}^*$ can be computed sequentially using forward substitution. Specifically, the above equation can be written as a system of nonlinear equations:

{
\onlyonefont
\renewcommand\arraystretch{1}
\setlength{\arraycolsep}{1.5pt}
\beq \label{eq:sub:nonlinear}
\bbb{0} \in \begin{bmatrix}
\bbb{B}_{1,1}& 0 & 0 & 0 &0  \\
\bbb{B}_{2,1}& \bbb{B}_{2,2} & 0 & 0 &0 \\
\vdots & \vdots  & \ddots &    0 &   0 \\
\bbb{B}_{n-1,1} &\bbb{B}_{n-1,2} & \cdots & \bbb{B}_{n-1,n-1}  & 0\\
\bbb{B}_{n,1}& \bbb{B}_{n,2} & \cdots & \bbb{B}_{n,n-1} &\bbb{B}_{n,n}  \\
\end{bmatrix} \begin{bmatrix}
  \bbb{z}^*_{1} \\
  \bbb{z}^*_{2}\\
  \vdots \\
  \bbb{z}^*_{n-1}   \\
  \bbb{z}^*_{n}   \\
\end{bmatrix} + \bbb{u} + \partial h(\bbb{z}^*) \nn
\eeq}{\normalsize}\noi If $\bbb{z}^*$ satisfies the equations above, it must solve the following one-dimensional subproblems:
\beq
0 \in \bbb{B}_{j,j} \bbb{z}^*_j + \bbb{w}_j +  \partial  h_j{(\bbb{z}^*_j)},~\forall j=1,2,~...~,n,\nn\\
\textstyle \bbb{w}_j=\bbb{u}_j +  \sum_{i=1}^{j-1} \bbb{B}_{j,i}\bbb{z}^*_{i}~~~~~~~~~~~~~~~~~\nn
\eeq
\noi This is equivalent to solving the following one-dimensional problem for all $j=1,2,...,n$:
\beq\label{eq:1d:subp}
\textstyle \bbb{z}^*_j= t^* \triangleq \underset{t}{\arg\min}~~\tfrac{1}{2}\bbb{B}_{j,j} t^2 + \bbb{w}_j t +  h_j(t)
\eeq
\noi Note that the computation of $\bbb{z}^{*}$ uses only the elements of $\bbb{z}^{*}$ that have already been computed and a successive displacement strategy is applied to find $\bbb{z}^{*}$.

We remark that the one-dimensional subproblem in (\ref{eq:1d:subp}) often admits a closed form solution for many problems of interest. For example, when $h_j(t)=I_{[lb,ub]}(t),~\forall j =1,2,...,n$ with $I_{[lb,ub]}(t)$ denoting an indicator function on the box constraint $lb\leq t\leq ub$, the optimal solution can be computed as: $t^* = \min(ub,\max(lb,-\bbb{w}_j/\bbb{B}_{j,j}))$; when $h_j(t)=\lambda |t|,~\forall j=1,2,...,n$ (\emph{i.e.} in the case of the $\ell_1$ norm), the optimal solution can be computed as: $t^* =   - \max\left(0,| \bbb{w}_j/\bbb{B}_{j,j}|-\lambda/\bbb{B}_{j,j}\right) \cdot \sign\left(\bbb{w}_j/\bbb{B}_{j,j}\right) $.

Our generalized Gaussian elimination procedure for computing $\P(\bbb{x}^k)$ is summarized in Algorithm \ref{alg:sub}. Note that its computational complexity is $\mathcal{O}(n^2)$, which is the same as computing a matrix-vector product.

\begin{algorithm} [!t]
\algsize
\caption{\label{alg:main} {\bbb{GMSA}: A Generalized Matrix Splitting Algorithm for Solving the Composite Function Minimization Problem in (\ref{eq:main})}}
\begin{algorithmic}[1]
\STATE Choose suitable parameters $\{\omega,\epsilon\}$.~Initialize $\bbb{x}^0$, ~$k=0$.\\
\STATE \text{while not converge}\\
\STATE~~~$\bbb{x}^{k+1} = \P(\bbb{x}^k)$ (Solve (\ref{eq:subproblem}) by Algorithm \ref{alg:sub})
\STATE~~~$k = {k+1}$
\STATE \text{end while}\\
\STATE Output $\bbb{x}^{k+1}$\\
\end{algorithmic}
\end{algorithm}

\begin{algorithm} [!t]
\algsize
\caption{\label{alg:sub} {A Generalized Gaussian Elimination Procedure for Computing the Triangle Proximal Operator $\P(\bbb{x}^k)$.}}
\begin{algorithmic}[1]
\STATE Input $\bbb{x}^k$\\
\STATE Initialization: compute $\bbb{u}= \bbb{b} + \bbb{C}\bbb{x}^k$\\
\STATE $\bbb{x}_1=\arg\min_{t}\frac{1}{2}\bbb{B}_{1,1}t^2 + (\bbb{u}_1) t +  h_1(t)$\\
\STATE $\bbb{x}_2=\arg\min_{t}\frac{1}{2}\bbb{B}_{2,2}t^2 + (\bbb{u}_2   + \bbb{B}_{2,1}\bbb{x}_{1}) t +  h_2(t)$\\
\STATE $\bbb{x}_3=\arg\min_{t}\frac{1}{2}\bbb{B}_{3,3}t^2 + (\bbb{u}_3+ \bbb{B}_{3,1}\bbb{x}_{1} + \bbb{B}_{3,2}\bbb{x}_{2}) t +  h_3(t)$\\
\STATE ...\\
\STATE $\bbb{x}_{n}=\arg\min_{t}\frac{1}{2}\bbb{B}_{n,n}t^2 + (\bbb{u}_n +  \sum_{i=1}^{n-1} \bbb{B}_{n,i}\bbb{x}_{i}) t +  h_n(t)$\\
\STATE Collect $(\bbb{x}_1,\bbb{x}_2,\bbb{x}_3,...,\bbb{x}_n)^T$ as $\bbb{x}^{k+1}$ and Output $\bbb{x}^{k+1}$\\
\end{algorithmic}
\end{algorithm}

%


\subsection{Convergence Analysis}\label{sect:alg:convergence}
In what follows, we present our convergence analysis for Algorithm \ref{alg:main}.

\noindent The following lemma characterizes the optimality of the triangle proximal operator $\P(\bbb{x})$ for any $\bbb{x}$.

\begin{lemma} \label{lemma:opt:ineq}
For all $\bbb{x},\bbb{y}\in\mathbb{R}^n$, it holds that:
\beq
\text{(i)}~ \bbb{0} \in ~\nabla q(\P(\bbb{x})) + \partial h(\P(\bbb{x})) + \bbb{C}(\bbb{x}-\P(\bbb{x}))   ~ \label{eq:opt:bound0}
\eeq
\vspace{-16pt}
\begin{align}
\begin{split}
\text{(ii)}~h(\P(\bbb{x})) -h(\bbb{y}) + \la \nabla q(\P(\bbb{x})),\P(\bbb{x}) - \bbb{y} \ra~~~~~~~~ \\
 \leq \la \bbb{C}(\P(\bbb{x})-\bbb{x}),~\P(\bbb{x})-\bbb{y} \ra  ~~~~~~~~~~~~~  \label{eq:opt:bound00}
\end{split}
\end{align}
\vspace{-5pt}
\begin{equation}\begin{split}\label{eq:opt:bound}
\text{(iii)}~&~f(\P(\bbb{x}))- f(\bbb{y}) ~~~~~~~~~~~~~~~~~~~~~~~~~~~~~~~~~~~~~~~~~~~\\
\leq~& ~  \la \bbb{C}(\P(\bbb{x})-\bbb{x}) ,\P(\bbb{x})-\bbb{y} \ra - \tfrac{1}{2}\|\P(\bbb{x})-\bbb{y}\|_{\bbb{A}}^2
\end{split}\end{equation}

\begin{proof}
(i) Using the optimality of $\P(\bbb{x})$ in (\ref{eq:subproblem}), we derive the following results: $\bbb{0} \in \bbb{B}\P(\bbb{x}) +  \partial h(\P(\bbb{x})) + \bbb{b} + \bbb{C} \bbb{x}  \overset{(a)}{\Rightarrow} \bbb{0} \in \bbb{A}\P(\bbb{x}) +  \partial h(\P(\bbb{x})) + \bbb{b} + \bbb{C} (\bbb{x}-\P(\bbb{x}))\overset{(b)}{\Rightarrow} \bbb{0} \in   \nabla q(\P(\bbb{x})) +  \partial h(\P(\bbb{x}))  + \bbb{C} (\bbb{x}-\P(\bbb{x}))$, where step $(a)$ uses $\bbb{B}=\bbb{A}-\bbb{C}$ and step $(b)$ uses the definition of $\nabla q(\cdot)$ in (\ref{eq:opt:cond}).

(ii) Since $h(\cdot)$ is convex, we have:
\beq\label{eq:fail}
\forall \bbb{s},~\bbb{z},~   h(\bbb{s}) -  h(\bbb{z}) \leq   \la h' ,\bbb{s} - \bbb{z} \ra,~\forall h' \in \partial h(\bbb{s}).
\eeq
\noi Letting $s=\P(\bbb{x}),~\bbb{z}=\bbb{y}$, we derive the following inequalities: $\forall h' \in \partial h(\P(\bbb{x})),~h(\P(\bbb{x})) -  h(\bbb{y})\overset{}{\leq}~ \la h' ,\P(\bbb{x}) - \bbb{y} \ra
\overset{(a)}{\leq}~ \la - \nabla q(\P(\bbb{x})) - \bbb{C}(\bbb{x}-\P(\bbb{x}))  ,\P(\bbb{x}) - \bbb{y} \ra$, where step $(a)$ uses (\ref{eq:opt:bound0}).

(iii) Since $q(\cdot)$ is a quadratic function, we have:
\beq \label{eq:quadratic}
\forall \bbb{s},~\bbb{z},~q(\bbb{s}) - q(\bbb{z})= \la \nabla q(\bbb{s}),\bbb{s}-\bbb{z} \ra - \tfrac{1}{2}\|\bbb{s}-\bbb{z}\|_{\bbb{A}}^2
\eeq
\noi We naturally derive the following results: $f(\P(\bbb{x}))- f(\bbb{y})$$\overset{(a)}{=}$$h(\P(\bbb{x})) - h(\bbb{y})+q(\P(\bbb{x})) - q(\bbb{y})$$\overset{(b)}{=}$$h(\P(\bbb{x})) - h(\bbb{y})+\la \nabla q(\P(\bbb{x})),\P(\bbb{x})-\bbb{y} \ra$$ - \tfrac{1}{2}\|\P(\bbb{x})-\bbb{y}\|_{\bbb{A}}^2$$\overset{(c)}{\leq}$$\la \bbb{C}(\P(\bbb{x})-\bbb{x}),\P(\bbb{x})-\bbb{y} \ra - $$\tfrac{1}{2}\|\P(\bbb{x})-\bbb{y}\|_{\bbb{A}}^2$, where step $(a)$ uses the definition of $f(\cdot)$; step $(b)$ uses (\ref{eq:quadratic}) with $\bbb{s}=\P(\bbb{x})$ and $\bbb{z}=\bbb{y}$; step $(c)$ uses (\ref{eq:opt:bound00}).

\end{proof}
\end{lemma}

\noi \bbb{Remarks.}  Both (\ref{eq:opt:bound0}) and (\ref{eq:opt:bound00}) can be used to characterize the optimality of (\ref{eq:main}). Recall that we have the following sufficient and necessary conditions for the optimal solution: $\bbb{x}^* \text{~is the optimal solution} \Leftrightarrow \bbb{0} \in \nabla q(\P(\bbb{x}^*)) + \partial h(\P(\bbb{x}^*)) $ $ \Leftrightarrow \la \nabla q(\P(\bbb{x}^*)),\P(\bbb{x}^*) - \bbb{y} \ra  + h(\P(\bbb{x}^*)) -h(\bbb{y}) \leq 0,~\forall \bbb{y}$. When $\bbb{x}=\P(\bbb{x})$ occurs, (\ref{eq:opt:bound0}) and (\ref{eq:opt:bound00}) coincide with the optimal condition and one can conclude that $\bbb{x}$ is the optimal solution.

\begin{theorem} \label{theorem:1}
(Proof of Global Convergence) We define $\delta \triangleq 2 \epsilon + \tfrac{2-\omega}{\omega}  \min(diag(\bbb{D}))$ and let $\{\omega,\epsilon\}$ be chosen such that $\delta\in(0,\infty)$. Algorithm \ref{alg:main} is globally convergent.

\begin{proof}
(i) First, the following results hold for all $\bbb{z}\in\mathbb{R}^n$:

\begin{align}  \label{eq:upperbound:0}
\bbb{z}^T(\bbb{A}-2\bbb{C})\bbb{z} = &~\textstyle \bbb{z}^T( \bbb{B}-\bbb{C}) \bbb{z} \nn\\
 =&~\textstyle \textstyle \bbb{z}^T(\bbb{L} - \bbb{L}^T  + \tfrac{2-\omega}{\omega}\bbb{D} + 2 \epsilon  \bbb{I})\bbb{z} \nn\\
=&~\textstyle \bbb{z}^T( {2 \epsilon} \bbb{I} + \tfrac{2-\omega}{\omega}\bbb{D} )\bbb{z} \geq  \delta \|\bbb{z}\|_2^2
\end{align}
\noi where we have used the definition of $\bbb{A}$ and $\bbb{C}$, and the fact that $\bbb{z}^T \bbb{L}\bbb{z} = (\bbb{z}^T \bbb{L}\bbb{z} )^T = \bbb{z}^T\bbb{L}^T\bbb{z},~\forall \bbb{z}$.

We invoke (\ref{eq:opt:bound}) in Lemma \ref{lemma:opt:ineq} with $\bbb{x}=\bbb{x}^k,~\bbb{y}=\bbb{x}^k$ and combine the inequality in (\ref{eq:upperbound:0}) to obtain:
\beq \label{eq:descent}
 \textstyle f^{k+1} - f^k  \textstyle \leq  -\frac{1}{2} \la \bbb{d}^{k},(\bbb{A}-2\bbb{C})\bbb{d}^{k}\ra \leq \textstyle -\frac{\delta}{2} \|\bbb{d}^k\|_2^2
\eeq
\noi (ii) Second, summing (\ref{eq:descent}) over $i=0,...,k-1$, we have:
\begin{align}
& ~\textstyle \tfrac{\delta}{2}\sum_{i=0}^{k-1} \|\bbb{d}^{i}\|_2^2 \leq  f^0  - f^{k}\overset{(a)}{\leq} \textstyle  f^0  - f^*\nn\\
\Rightarrow &~\textstyle \tfrac{\delta}{2} \min_{i=0,...,k-1}~\|\bbb{d}^{i}\|_2^2 \leq (f^0  - f^*)/k\nn
\end{align}
\noi where step $(a)$ uses the fact that $f^*\leq f^{k}$. Note that $f^*$ is bounded below. As $k\rightarrow \infty$, we have $ \bbb{d}^{k}  \triangleq \bbb{x}^{k+1}-\bbb{x}^k \rightarrow \bbb{0}$, which implies the convergence of the algorithm. Invoking (\ref{eq:opt:bound0}) in Lemma \ref{lemma:opt:ineq} with $\bbb{x}=\bbb{x}^k$, we obtain: $\nabla q(\bbb{x}^{k+1}) + \partial h(\bbb{x}^{k+1}) \ni - \bbb{C}(\bbb{x}^k- \bbb{x}^{k+1} ) \rightarrow \bbb{0}$. The fact that $\nabla q(\bbb{x}^{k+1}) + \partial h(\bbb{x}^{k+1}) \ni \bbb{0}$ implies that $\bbb{x}^{k+1}$ is the global optimal solution of the convex problem.

Note that guaranteeing $\delta\in(0,\infty)$ can be achieved by simply choosing $\omega\in(0,2)$ and setting $\epsilon$ to a small number.
\end{proof}
\end{theorem}
\noi \bbb{Remarks.} \textbf{(i)} When $h(\cdot)$ is empty and $\epsilon=0$, Algorithm \ref{alg:main} reduces to the classical Gauss-Seidel method ($\omega=1$) and Successive Over-Relaxation method ($\omega\neq1$). \textbf{(ii)} When $\bbb{A}$ contains zeros in its diagonal entries, one needs to set $ \epsilon$ to a strictly positive number. This is to guarantee the strong convexity of the one dimensional subproblem and a bounded solution for any $h(\cdot)$ in (\ref{eq:1d:subp}). The introduction of the parameter $\epsilon$ is novel in this paper and it removes the assumption that $\bbb{A}$ is strictly positive-definite or strictly diagonally dominant, which is used in the classical result of GS and SOR method \cite{saad2003iterative,demmel1997applied}.

We now prove the convergence rate of Algorithm \ref{alg:main}. We make the following assumption, which characterizes the relations between $\P(\bbb{x})$ and $\bbb{x}^*$ for any $\bbb{x}$.

\begin{assumption}\label{lemma:local:bound}
If $\bbb{x}$ is not the optimum of (\ref{eq:main}), there exists a constant $\eta\in(0,\infty)$ such that $\|\bbb{x} - \bbb{x}^* \| \leq \eta \|\bbb{x} - \P(\bbb{x})\|$.
\end{assumption}

\noi \bbb{Remarks.} Assumption \ref{lemma:local:bound} is similar to the classical local proximal error bound assumption in the literature \cite{luo1993error,TsengY09,Tseng10,yun2011block}, and it is mild. Firstly, if $\bbb{x}$ is not the optimum, we have $\bbb{x} \neq \P(\bbb{x})$. This is because when $\bbb{x} = \P(\bbb{x})$, we have $\bbb{0} = -\bbb{C}(\bbb{x}-\P(\bbb{x})) \in  \nabla q(\P(\bbb{x})) + \partial h(\P(\bbb{x}))$ (refer to the optimal condition of $\P(\bbb{x})$ in (\ref{eq:opt:bound0})), which contradicts with the condition that $\bbb{x}$ is not the optimal solution. Secondly, by the boundedness of $\bbb{x}$ and $\bbb{x}^*$, there exists a sufficiently large constant $\eta\in(0,\infty)$ such that $\|\bbb{x} - \bbb{x}^* \| \leq \eta \|\bbb{x} - \P(\bbb{x})\|$.


%
%
%

We now prove the convergence rate of Algorithm \ref{alg:main}.

\begin{theorem} \label{theorem:general:rate}
(Proof of Convergence Rate) We define $\delta \triangleq {2 \epsilon} + \tfrac{2-\omega}{\omega}  \min(diag(\bbb{D}))$ and let $\{\omega,~\epsilon\}$ be chosen such that $\delta\in(0,\infty)$. Assuming that $\bbb{x}^k$ is bound for all $k$, we have:
\beq
 f(\bbb{x}^{k}) - f(\bbb{x}^*) \leq \left(\frac{C_1}{1+C_1}\right)^k [f(\bbb{x}^{0}) - f(\bbb{x}^*)], \label{eq:QQ2}\\
 \|\bbb{x}^k - \bbb{x}^{k+1}\|_2^2  \leq  \frac{2}{\delta} \left(\frac{C_1}{1+C_1}\right)^k [f(\bbb{x}^{0}) - f(\bbb{x}^*)].  \label{eq:QQ3} 
\eeq
\noi where $C_1 \triangleq 2 \|\bbb{B}\|\eta/\delta - 1$.
\begin{proof}


Invoking Assumption \ref{lemma:local:bound} with $\bbb{x}=\bbb{x}^k$, we obtain:
\beq \label{eq:opt:bound:ineq}
\|\bbb{x}^k - \bbb{x}^*\| \leq \eta \|\bbb{x}^k - \P(\bbb{x}^k)\| ~\Rightarrow~\|\bbb{r}^k\| \leq \eta \|\bbb{d}^k\|
\eeq
\noi We derive the following inequalities:
\begin{align}
& ~~~~~~~f^{k+1} - f^* \nn \\
&\overset{(a)}{\leq}~ \textstyle  \la \bbb{r}^{k+1} ,\bbb{C} \bbb{d}^k  \ra   - \tfrac{1}{2}\la \bbb{r}^{k+1},\bbb{A}\bbb{r}^{k+1}\ra  \label{eq:linear:conv0}\\
&\overset{(b)}{=}~ \textstyle  \la \bbb{r}^k,~(\bbb{C}-\bbb{A}) \bbb{d}^k \ra  -\tfrac{1}{2} \|\bbb{r}^{k}\|_{\bbb{A}}^2 + \tfrac{1}{2} \|\bbb{d}^{k}\|_{2\bbb{C}-\bbb{A}}^2  \nn \\
&\overset{(c)}{\leq}~ \textstyle   - \la \bbb{r}^k,~\bbb{B} \bbb{d}^k \ra + 0 - \tfrac{\delta}{2} \|\bbb{d}^{k}\|_{2}^2    \nn\\
&\overset{(d)}{\leq}~ \textstyle  \|\bbb{r}^k\| \|\bbb{B}\|  \|\bbb{d}^k\|  - \tfrac{\delta}{2} \|\bbb{d}^{k}\|_{2}^2  \overset{(e)}{\leq}~ \textstyle    (\eta  \|\bbb{B}\|- \tfrac{\delta}{2}      )  \|\bbb{d}^k\|_2^2 \nn \\
&\overset{(f)}{\leq}~ \textstyle    (\eta  \|\bbb{B}\|- \tfrac{\delta}{2}      )  \frac{2}{\delta}  (f^k-f^{k+1}) \overset{(g)}{=}~ \textstyle      C_1 (f^k-f^{k+1})\label{eq:linear:conv}
\end{align}
\noi where step $(a)$ uses (\ref{eq:opt:bound}) in Lemma \ref{lemma:opt:ineq} with $\bbb{x}=\bbb{x}^*,~\bbb{y}=\bbb{x}^k$; step $(b)$ uses the fact that $\bbb{r}^{k+1}=\bbb{r}^{k}+\bbb{d}^{k}$ and $\bbb{A}=\bbb{B}+\bbb{C}$; step $(c)$ uses $\bbb{A} \succeq \bbb{0}$ and the inequality that $\bbb{A}-2\bbb{C}\succeq \delta \bbb{I}$ which is due to (\ref{eq:upperbound:0}); step ($d$) uses the Cauchy-Schwarz inequality $\la \bbb{x},\bbb{y} \ra\leq \|\bbb{x}\|\|\bbb{y}\|,~\forall \bbb{x},\bbb{y}\in\mathbb{R}^n$ and the norm inequality $\|\bbb{Bx}\|\leq \|\bbb{B}\|\|\bbb{x}\|,~\forall \bbb{x} \in\mathbb{R}^n$; step ($e$) uses (\ref{eq:opt:bound:ineq}); step $(f)$ uses the descent condition in (\ref{eq:descent}); step $(g)$ uses the definition of $C_1$.

Rearranging the last inequality in (\ref{eq:linear:conv}), we have $ f^{k+1} - f^* \leq C_1 (f^k-f^{k+1}) = \textstyle C_1  (f^k-f^*) - C_1(f^{k+1}-f^*) \Rightarrow (1+C_1)[f(\bbb{x}^{k+1}) - f(\bbb{x}^*)] \leq C_1 [f(\bbb{x}^k)-f(\bbb{x}^*)]$, leading to: $\tfrac{f(\bbb{x}^{k+1}) - f(\bbb{x}^*)}{f(\bbb{x}^k)-f(\bbb{x}^*)} \leq \tfrac{C_1}{1+C_1} < 1$. Solving this recursive formulation, we obtain (\ref{eq:QQ2}). In other words, the sequence $\{f(\bbb{x}^k)\}_{k=0}^{\infty}$ converges to $f(\bbb{x}^*)$ linearly in the quotient sense. Using (\ref{eq:descent}), we derive the following inequalities: $\|\bbb{x}^k - \bbb{x}^{k+1}\|_2^2 \leq \frac{2 (f^k - f^{k+1})}{\delta} \leq \frac{2 (f^k - f^* ) }{\delta}$. Therefore, we obtain (\ref{eq:QQ3}).

\end{proof}
\end{theorem}

The following lemma is useful in our proof of iteration complexity. 

\begin{lemma} \label{lemma:quadratic:recursive}
Suppose a nonnegative sequence $\{u^k\}_{k=0}^{\infty}$ satisfies $u^{k+1} \leq -2 C + 2C \sqrt{1+ \frac{u^k}{C} }$ for some constant $C> 0$. It holds that: $u^{k+1} \leq  \frac{\max(8C,\sqrt{4 Cu^0})}{k+1}$.
\begin{proof}
 The proof of this lemma can be obtained by mathematical induction. We denote $\chi \triangleq \max(8C,\sqrt{4Cu^0})$. (i) When $k=0$, we have $u^{1}  \leq -2 C + 2C \sqrt{1+ \frac{1}{C} u^0} \leq -2C + 2C (1+\sqrt{\frac{u^0}{C}} )=2\sqrt{Cu^0}\leq \frac{\chi}{k+1}$. (ii) When $k\geq 1$, we assume that $u^{k} \leq  \frac{\chi}{k}$ holds. We derive the following results: $k\geq 1 \Rightarrow \frac{k+1}{k} \leq 2$ $\overset{(a)}{\Rightarrow} 4C \frac{k+1}{k} \leq 8C \leq \chi$ $\overset{(b)}{\Rightarrow}  \frac{4C}{k(k+1)} =  4C (\frac{1}{k} - \frac{1}{k+1} ) \leq \frac{\chi}{(k+1)^2}$ $\Rightarrow   \frac{4C}{k} \leq    \frac{4C}{k+1} + \frac{\chi}{(k+1)^2}$ $\Rightarrow   4C^2 ( 1+ \frac{\chi}{k C})  \leq  4C^2 +   \tfrac{4\chi C}{k+1} + \tfrac{\chi^2}{(k+1)^2}$ $\Rightarrow   2C \sqrt{1+\frac{\chi}{k C}} \leq 2C + \frac{\chi}{k+1}$ $\Rightarrow   -2C + 2C \sqrt{1+\frac{\chi}{kC}} \leq \frac{\chi}{k+1}$ $\overset{(c)}{\Rightarrow}   -2C + 2C \sqrt{1+\frac{u^k}{C}} \leq \frac{\chi}{k+1}$ $\Rightarrow   u^{k+1} \leq \frac{\chi}{k+1}$. Here, step $(a)$ uses $8C\leq \chi$; step $(b)$ uses the inequality that $\frac{1}{k(k+1)} = \frac{1}{k}-\frac{1}{k+1}$; step $(c)$ uses $u^k \leq \frac{\chi}{k}$.

\end{proof}
\end{lemma}

We now prove the iteration complexity of Algorithm \ref{alg:main}.

\begin{theorem} \label{theorem:general:rate2}
(Proof of Iteration Complexity) We define $\delta \triangleq {2 \epsilon} + \tfrac{2-\omega}{\omega}  \min(diag(\bbb{D}))$ and let $\{\omega,~\epsilon\}$ be chosen such that $\delta \in(0,\infty)$. Assuming that $\|\bbb{x}^k\|\leq R$ for all $k$, we have:
\beq
\textstyle
f^k - f^* \leq
   \begin{cases}
   u^0(\frac{2C_3}{2C_3+1})^k,  &\mbox{if~$\sqrt{f^{k}-f^{k+1}} \geq \frac{C_2}{C_3}$},~\forall k\leq \breve{k}\nn \\
   \frac{C_4}{k}, &\mbox{if~$\sqrt{f^{k}-f^{k+1}} < \frac{C_2}{C_3}$}, ~\text{else}
   \end{cases}
\eeq
\noi where $C_2 \triangleq 2 R\|\bbb{C}\| \sqrt{{2}/{\delta}}$, $C_3 \triangleq \frac{2}{\delta}\|\bbb{C}\|$, $C_4 \triangleq \max(8C_2^2,\sqrt{4C_2^2u^0})$, and $\breve{k}$ is some unknown iteration index.

\begin{proof}
We have the following inequalities:
\begin{align}
&~~~~~~ u^{k+1} \nn\\
& \overset{(a)}{\leq} \textstyle \la \bbb{r}^{k+1} ,\bbb{C} \bbb{d}^k  \ra   - \tfrac{1}{2}\la \bbb{r}^{k+1},\bbb{A}\bbb{r}^{k+1}\ra \nn\\
& \overset{(b)}{\leq} \textstyle    \la \bbb{r}^{k}+\bbb{d}^{k} ,\bbb{C} \bbb{d}^k  \ra   +0 \nn \\
& \overset{(d)}{\leq} \textstyle   \|\bbb{r}^{k}\| \cdot \|\bbb{C}\| \cdot\|\bbb{d}^k\|   + \|\bbb{C}\| \cdot \|\bbb{d}^k\|_2^2 \nn \\
& \overset{(d)}{\leq} \textstyle   2 R\|\bbb{C}\| \cdot\|\bbb{d}^k\|   + \|\bbb{C}\| \cdot \|\bbb{d}^k\|_2^2 \nn \\
& \overset{(e)}{\leq} \textstyle   2 R\|\bbb{C}\| \cdot \sqrt{\tfrac{2}{\delta} (u^k-u^{k+1})}   + \|\bbb{C}\| \cdot \tfrac{2}{\delta} \cdot (u^k-u^{k+1})  \nn \\
&\overset{(f)}{=} C_2 \sqrt{u^k-u^{k+1}} + C_3 (u^k-u^{k+1})\label{eq:recursion:u}
\end{align}
\noi where step $(a)$ uses (\ref{eq:linear:conv0}); step $(b)$ uses the fact that $\bbb{r}^{k+1}=\bbb{r}^{k}+\bbb{d}^{k},~\bbb{A}\succeq \bbb{0}$; step $(c)$ uses the Cauchy-Schwarz inequality and the norm inequality; step $(d)$ uses the fact that $\|\bbb{r}^{k}\|_2 \leq \|\bbb{x}^{k}\|_2 + \|\bbb{x}^{*}\|_2\leq 2R$; step $(e)$ uses (\ref{eq:descent}); step $(f)$ uses the definition of $C_2$ and $C_3$.

Now we consider the two cases for the recursion formula in (\ref{eq:recursion:u}): (i) $\sqrt{u^k-u^{k+1}} \geq \frac{C_2}{C_3}$ for some $k\leq \breve{k}$ (ii) $\sqrt{u^k-u^{k+1}} \leq \frac{C_2}{C_3}$ for some $k>\breve{k}$. In case (i), (\ref{eq:recursion:u}) implies that we have $u^{k+1}\leq  2 C_3 (u^{k}-u^{k+1})$ and rearranging terms gives: $u^{k+1}\leq \frac{2C_3}{2C_3+1} u^k$. Thus, we have: $u^{k+1}\leq (\frac{2C_3}{2C_3+1})^{k+1} u^0$. We now focus on case (ii). When $\sqrt{u^k-u^{k+1}} \leq \frac{C_2}{C_3}$, (\ref{eq:recursion:u}) implies that we have $u^{k+1}\leq  2 C_2 \sqrt{u^{k}-u^{k+1}}$ and rearranging terms yields:$\frac{(u^{k+1})^2}{4 C_2^2} + u^{k+1}  \leq   u^{k}$. Solving this quadratic inequality, we have: $u^{k+1} \leq  -2 C_2^2 + 2 C_2^2 \sqrt{1+\frac{1}{C_2^2} u^k}$; solving this recursive formulation by Lemma \ref{lemma:quadratic:recursive}, we obtain $ u^{k+1} \leq \frac{C_4}{k+1}$.

\end{proof}
\end{theorem}

\noi \bbb{Remarks.} The convergence result in Theorem \ref{theorem:general:rate2} is weaker than that in Theorem \ref{theorem:general:rate}, however, it does not rely on Assumption \ref{lemma:local:bound} and the unknown constant $\eta$.

We now derive the convergence rate when $q (\cdot)$ is strongly convex.

\begin{theorem}\label{theorem:general:rate3}
(Proof of Convergence Rate when $q(\cdot)$ is Strongly Convex) We define $\delta \triangleq {2 \epsilon} + \tfrac{2-\omega}{\omega}  \min(diag(\bbb{D}))$ and let $\{\omega,~\epsilon\}$ be chosen such that $\delta\in(0,\infty)$. Assuming that $q(\bbb{x})$ is strongly convex with respect to $\bbb{x}$ such that $\bbb{A} \succeq \sigma \bbb{I}$ with $\sigma>0$ and $\|\bbb{x}^k\|\leq R$ for all $k$, we have:
\beq \label{eq:aaa}
f(\bbb{x}^{k}) - f(\bbb{x}^*) \leq \left(\tfrac{C_5}{1+C_5}\right)^k [f(\bbb{x}^{0}) - f(\bbb{x}^*)], \label{eq:QQ22} \label{eq:conv:rate:strong:f} \\
\|\bbb{x}^k - \bbb{x}^*\|^2_2 \leq     \frac{8\| \bbb{C}\|^2}{\sigma^2\delta} \left(\tfrac{C_5}{1+C_5}\right)^k [f(\bbb{x}^{0}) - f(\bbb{x}^*)].~~~ \label{eq:conv:rate:strong:x}
\eeq
\noi where $C_5 \triangleq {\|\bbb{C}\|^2}/{(\delta\sigma)}$.

\begin{proof}
Invoking (\ref{eq:opt:bound}) in Lemma \ref{lemma:opt:ineq} with $\bbb{x}=\bbb{x}^k,~\bbb{y}=\bbb{x}^*$, we derive the following inequalities:
\begin{align} \label{eq:maximization}
&~f(\bbb{x}^{k+1}) - f(\bbb{x}^*)~~~~~~~~~~~~~~~~~~~~~~~~~~~~~~~\nn\\
\leq& ~ \la \bbb{C}(\bbb{x}^{k+1}-\bbb{x}^{k}),\bbb{x}^{k+1}-\bbb{x}^* \ra  - \tfrac{1}{2}\|\bbb{x}^{k+1}-\bbb{x}^{*}\|_{\bbb{A}}^2 \nn\\
\leq&~  \la \bbb{C}(\bbb{x}^{k+1}-\bbb{x}^{k}),\bbb{x}^{k+1}-\bbb{x}^* \ra - \tfrac{\sigma}{2}\|\bbb{x}^{k+1}-\bbb{x}^{*}\|_2^2
\end{align}

\noi We notice that the right-hand side in (\ref{eq:maximization}) is concave. Maximizing over $\bbb{x}^*$, we obtain:
\begin{align} \label{eq:optimal:solution}
& ~~\sigma (\bbb{x}^{*}-\bbb{x}^{k+1}) + \bbb{C}(\bbb{x}^{k+1}-\bbb{x}^{k}) = 0 \nn\\
 \Rightarrow & ~~\bbb{x}^{*} = \bbb{x}^{k+1} -  { \bbb{C}(\bbb{x}^{k+1}-\bbb{x}^{k})}/{\sigma}
\end{align}
\noi Putting (\ref{eq:optimal:solution}) into (\ref{eq:maximization}), we further derive the following inequalities: 
\begin{align}
f(\bbb{x}^{k+1}) - f(\bbb{x}^*)\overset{}{\leq} \textstyle \tfrac{\|\bbb{C}(\bbb{x}^{k+1}-\bbb{x}^{k})\|_2^2}{2 \sigma} \overset{(a)}{\leq} \textstyle \frac{\|\bbb{C}\|^2\cdot \|\bbb{x}^k-\bbb{x}^{k+1}\|_2^2}{2\sigma} \nn\\
\overset{(b)}{\leq} \textstyle \frac{\|\bbb{C}\|^2 \cdot [f(\bbb{x}^k)-f(\bbb{x}^{k+1})]}{\delta\sigma} \overset{(c)}{=}  C_5 [f(\bbb{x}^k)-f(\bbb{x}^{k+1})]\nn
\end{align}
\noi where step $(a)$ uses the norm inequality $\|\bbb{Cx}\|\leq \|\bbb{C}\|\cdot\|\bbb{x}\|$; step $(b)$ uses (\ref{eq:descent}); step $(c)$ uses the definition of $C_5$.  Finally, we obtain: $\tfrac{f(\bbb{x}^{k+1}) - f(\bbb{x}^{*})}{f(\bbb{x}^{k}) - f(\bbb{x}^{*}) } \leq \tfrac{C_5}{1+C_5}$. Solving the recursive formulation, we obtain the result in (\ref{eq:conv:rate:strong:f}).

Using the similar strategy for deriving (\ref{eq:QQ3}), we have:
\beq \label{eq:similar}
\|\bbb{x}^k - \bbb{x}^{k+1}\|_2^2 \leq \tfrac{2}{\delta} (\tfrac{C_5}{1+C_5})^k [f(\bbb{x}^{0}) - f(\bbb{x}^*)]
\eeq

\noi Finally, we derive the following inequalities:
\begin{align}
&~\tfrac{\sigma}{2}\|\bbb{x}^{k+1} - \bbb{x}^*\|^2_2 \nn\\
\overset{(a)}{\leq} &~\la \bbb{C}(\bbb{x}^{k+1}-\bbb{x}^{k}),\bbb{x}^{k+1}-\bbb{x}^* \ra + f(\bbb{x}^*) - f(\bbb{x}^{k+1}) \nn\\
\overset{(b)}{\leq} & ~\la \bbb{C}(\bbb{x}^{k+1}-\bbb{x}^{k}),\bbb{x}^{k+1}-\bbb{x}^* \ra \nn\\
\overset{(c)}{\leq} & ~\| \bbb{C}\|\cdot\|\bbb{x}^{k+1}-\bbb{x}^*\| \cdot \|\bbb{x}^{k+1}-\bbb{x}^{k}\| \nn
\nn
\end{align}
\noi where step $(a)$ uses (\ref{eq:maximization}); step $(b)$ uses the fact that $f(\bbb{x}^*)\leq f(\bbb{x}^{k+1})$; step $(c)$ uses the norm inequality. Therefore, we obtain:
\beq 
\tfrac{\sigma}{2} \| \bbb{x}^{k+1} - \bbb{x}^* \|_2 \leq \|\bbb{C}\|\|\bbb{x}^{k+1} - \bbb{x}^k\| \nn
\eeq
\noi Combining with (\ref{eq:similar}), we obtain (\ref{eq:conv:rate:strong:x}).

\end{proof}

\end{theorem}

\noi \bbb{Remarks.} Thanks to the strongly convexity of $q(\cdot)$, we can characterize the convergence rate for both $\|\bbb{x}^k-\bbb{x}^*\|$ and $\|\bbb{x}^k-\bbb{x}^{k+1}\|$ in Theorem \ref{theorem:general:rate3} without using Assumption \ref{lemma:local:bound} and the unknown constant $\eta$. Therefore, the convergence result in Theorem \ref{theorem:general:rate3} is stronger than that in Theorem \ref{theorem:general:rate}.


\section{Extensions}\label{sect:extension}
This section discusses several extensions of our proposed generalized matrix splitting algorithm.

\subsection{Pointwise Contraction via a Correction Strategy}\label{sect:ext:strong}

This section considers a new correction strategy to achieve pointwise contraction for the proposed method to solve (\ref{eq:main}). One remarkable feature of this strategy is that the resulting iterated solutions of $\bbb{x}^k$ always satisfy the monotone/contractive property that $\|\bbb{x}^{k+1}-\bbb{x}^*\|_2^2<\|\bbb{x}^{k}-\bbb{x}^*\|_2^2$ for all $k$ if $\bbb{x}^{k}$ is not the optimal solution. We summarize our new algorithm in Algorithm \ref{alg:gmsa:C}. 

We provide detailed theoretical analysis for Algorithm \ref{alg:gmsa:C}. The following lemmas are useful in our proof.
\begin{algorithm} [!t]
\algsize
\caption{\label{alg:gmsa:C} {GMSA-C: Generalized Matrix Splitting Algorithm with Correction Strategy for Solving (\ref{eq:main}).}}
\begin{algorithmic}[1]
\STATE Choose suitable parameters $\{\omega,~\epsilon\}$.~Initialize $\bbb{x}^0$.\\
\STATE \text{for $k=0,1,2,...$}\\
\STATE~~~$\bbb{y}^{k} = \P(\bbb{x}^k)$
\STATE~~~Choose a suitable parameter $\alpha^k$ (\text{e.g.}~$\alpha^k = \frac{2 \epsilon + \tfrac{2-\omega}{\omega}  \min(diag(\bbb{D}))}{\|\bbb{B}^T\bbb{B}\|}$\text{~or~}$\alpha^k = \frac{\|\bbb{y}^k - \bbb{x}^k\|_{2\bbb{B}-0.5\bbb{A}}^2}{\|\bbb{y}^k - \bbb{x}^k\|_{ 2\bbb{B}^T\bbb{B}}^2}$)
\STATE~~~$\bbb{x}^{k+1} = \bbb{x}^k + \alpha^k \bbb{B} (\bbb{y}^{k} - \bbb{x}^k )$
\STATE \text{end for}\\
\STATE Output $\bbb{x}^{k+1}$\\
\end{algorithmic}
\end{algorithm}
\begin{lemma}\label{eq:dis:key}
Assume that $h(\cdot)$ is convex. For all $\bbb{x} \in \mathbb{R}^n$, it holds that: 
$$\la \bbb{Ax}^* - \bbb{A}\P(\bbb{x}) + \bbb{C} (\P(\bbb{x})-\bbb{x}), ~\bbb{x}^* - \P(\bbb{x}) \ra \leq 0.$$

\begin{proof}
By the optimality of $\bbb{x}^*$, we have: $\bbb{0} \in ~\nabla q(\bbb{x}^*) + \partial h(\bbb{x}^*)$. Combining with (\ref{eq:opt:bound0}) in Lemma (\ref{lemma:opt:ineq}), we obtain:
\begin{align} \label{eq:dis:1}
\bbb{0}~  \in& ~\nabla q(\P(\bbb{x})) + \partial h(\P(\bbb{x})) + \bbb{C}(\bbb{x}-\P(\bbb{x})) \nn\\
&~ - \nabla q(\bbb{x}^*) - \partial h(\bbb{x}^*)
\end{align}
\noi Using the monotonicity of $\partial h(\cdot)$, we obtain: $\la h'-h'', \bbb{x}^*-\P(\bbb{x}) \ra \geq 0,~\forall h' \in \partial h(\bbb{x}^{*}),~h'' \in \partial h(\P(\bbb{x}))$. Combining with (\ref{eq:dis:1}), we conclude this lemma.

\end{proof}

\end{lemma}

\begin{lemma} \label{eq:hhh}
Assume $\bbb{A}\succeq \bbb{0}$. For all $\bbb{x},~\bbb{y},~\bbb{z} \in \mathbb{R}^n$, it holds that: 
$$\la \bbb{x}- \bbb{z}  , \bbb{Az} -\bbb{Ay}\ra \leq  \tfrac{1}{4} \|\bbb{x}-\bbb{y}\|_{\bbb{A}}^2.$$

\begin{proof}
Using the variable substitution  that $\bbb{z}-\bbb{x}=\bbb{p},~\bbb{z}-\bbb{y} = \bbb{u},~\bbb{x}-\bbb{y} = \bbb{u}-\bbb{p}$, we have the following equivalent inequalities: $\la -\bbb{p},\bbb{Au} \ra \leq \tfrac{1}{4}\|\bbb{u}-\bbb{p}\|_{\bbb{A}}^2 \Leftrightarrow  \|\bbb{u}-\bbb{p}\|_{\bbb{A}}^2 + 4 \la \bbb{Ap},\bbb{u} \ra \geq 0 \Leftrightarrow \|\bbb{u}+\bbb{p}\|_{\bbb{A}}^2  \geq 0$. Clearly, these inequalities hold since $\bbb{A}\succeq \bbb{0}$.

\end{proof}
\end{lemma}

The following theorem provides important theoretical insights on choosing suitable parameters $\{\omega,~\epsilon,~\alpha^k\}$ to guarantee convergence of Algorithm \ref{alg:gmsa:C}.

\begin{theorem} \label{theorem:strong}
We define $\delta \triangleq 2 \epsilon + \tfrac{2-\omega}{\omega}  \min(diag(\bbb{D}))$ and let $\{\omega,~\epsilon\}$ be chosen such that $\delta\in(0,\infty)$. Assuming that $h(\cdot)$ is convex and $\bbb{x}^k$ generated by Algorithm \ref{alg:gmsa:C} is not the optimal solution, we have the following results: (i)
\beq \label{eq:dis:ok}
&~~~~\|\bbb{x}^{k+1} - \bbb{x}^*\|_2^2 - \|\bbb{x}^{k} - \bbb{x}^*\|_2^2 \leq \|\bbb{y}^k - \bbb{x}^k\|_{\bbb{G}^k}^2 \\
&\bbb{G}^k\triangleq \frac{1}{2} (\alpha^k)^2 \bbb{P} - \alpha^k\bbb{Q},~\bbb{P} \triangleq 2\bbb{B}^T\bbb{B},~\bbb{Q} \triangleq 2 \bbb{B}-\tfrac{1}{2}\bbb{A}  \nn
\eeq
\noi (ii) If we choose a global constant $0<\alpha^k<\frac{\delta}{\|\bbb{B}^T\bbb{B}\|}$, we have $\bbb{G}^k\prec 0$ and $\|\bbb{x}^{k+1} - \bbb{x}^*\|_2^2 - \|\bbb{x}^{k} - \bbb{x}^*\|_2^2 <0$.

\noi (iii) If we choose a local constant $\alpha^k = \frac{\|\bbb{y}^k - \bbb{x}^k\|_{\bbb{Q}}^2}{\|\bbb{y}^k - \bbb{x}^k\|_{\bbb{P}}^2}$, we have $
\|\bbb{x}^{k+1} - \bbb{x}^*\|_2^2 - \|\bbb{x}^{k} - \bbb{x}^*\|_2^2 \leq - \tfrac{\delta^2 \|\bbb{y}^k - \bbb{x}^k\|_2^2 }{4 \|\bbb{BB}^T\|} <0$.

\begin{proof}
(i) First of all, we derive the following inequalities:
\begin{align} \label{eq:dis:prop:1}
&~\la \bbb{y}^{k}-\bbb{x}^{*}, \bbb{B} (\bbb{y}^{k} - \bbb{x}^k )\ra \nn\\
\overset{(a)}{=}&~ \la \bbb{y}^{k}-\bbb{x}^{*}, \bbb{A} (\bbb{y}^{k} - \bbb{x}^k )\ra  + \la \bbb{y}^{k}-\bbb{x}^{*}, \bbb{C} (\bbb{x}^k - \bbb{y}^{k} )\ra \nn\\
\overset{(b)}{\leq}&~\la \bbb{y}^{k}-\bbb{x}^{*}, \bbb{A} (\bbb{y}^{k} - \bbb{x}^k ) \ra + \la \bbb{y}^{k}-\bbb{x}^{*}, \bbb{A} (\bbb{x}^{*} - \bbb{y}^k )\ra   \nn\\
\overset{(c)}{=}&~\la \bbb{y}^{k}-\bbb{x}^{*}, \bbb{A} (\bbb{x}^{*} - \bbb{x}^k )\ra  \overset{(d)}{\leq}  \tfrac{1}{4} \|\bbb{y}^{k}-\bbb{x}^{k}\|_{\bbb{A}}^2
\end{align}
\noi where step ($a$) uses the fact that $\bbb{B}=\bbb{-C}+\bbb{A}$; step ($b$) Lemma \ref{eq:dis:key} with $\bbb{x}=\bbb{x}^k$; step $(c)$ uses the fact that $(\bbb{x}^{*} - \bbb{y}^k )+(\bbb{y}^{k} - \bbb{x}^k ) = (\bbb{x}^{*} - \bbb{x}^k )$; step $(d)$ uses Lemma \ref{eq:hhh}. We then have the following results:
\begin{align}
&~\|\bbb{x}^{k+1}-\bbb{x}^*\|_{2}^2 - \|\bbb{x}^{k}-\bbb{x}^*\|_{2}^2 \nn\\
\overset{(a)}{=}&~\|\bbb{x}^{k+1}-\bbb{x}^k\|_{2}^2 +  2\la \bbb{x}^k-\bbb{x}^{*},\bbb{x}^{k+1} - \bbb{x}^{k} \ra \nn  \\
\overset{(b)}{=}&~\|\alpha^k \bbb{B} (\bbb{y}^{k} - \bbb{x}^k )\|_{2}^2 + 2\alpha^k  \la  \bbb{x}^k-\bbb{x}^{*}, \bbb{B} (\bbb{y}^{k} - \bbb{x}^k )\ra\nn\\
\overset{(c)}{=}&~\|\alpha^k \bbb{B} (\bbb{y}^{k} - \bbb{x}^k )\|_{2}^2 + 2\alpha^k  \la  \bbb{x}^k - \bbb{y}^{k} , \bbb{B} (\bbb{y}^{k} - \bbb{x}^k )\ra  \nn\\
&~ + 2\alpha^k  \la \bbb{y}^{k}-\bbb{x}^{*}, \bbb{B} (\bbb{y}^{k} - \bbb{x}^k )\ra \nn\\
\overset{(d)}{\leq}&~\|\alpha^k \bbb{B} (\bbb{y}^{k} - \bbb{x}^k )\|_{2}^2 + 2\alpha^k  \la  \bbb{x}^k - \bbb{y}^{k} , \bbb{B} (\bbb{y}^{k} - \bbb{x}^k )\ra \nn\\
&~  + \tfrac{1}{2}\alpha^k \|\bbb{y}^{k} - \bbb{x}^k\|_{\bbb{A}}^2 \overset{(e)}{=}\|\bbb{y}^{k} - \bbb{x}^k\|_{\bbb{G}^k}^2\nn
\end{align}
\noi where step $(a)$ uses Pythagoras relation that $\|\bbb{x}-\bbb{z}\|_{2}^2 - \|\bbb{y}-\bbb{z}\|_{2}^2 =\|\bbb{x}-\bbb{y}\|_{2}^2 +  2\la \bbb{y}-\bbb{z},\bbb{x} - \bbb{y} \ra,~\forall \bbb{x},~\bbb{y},~\bbb{z}$; step $(b)$ uses the update rule for $\bbb{x}^{k+1}$; step $(c)$ uses the fact that $\bbb{x}^k - \bbb{x}^*=(\bbb{x}^k-\bbb{y}^k)+(\bbb{y}^k-\bbb{x}^*)$; step $(d)$ uses (\ref{eq:dis:prop:1}); step $(e)$ uses the definition of $\bbb{G}^k$.

(ii) We have the following inequalities:
\beq \label{eq:Q:delta}
\bbb{Q} \overset{(a)}{=} 2\bbb{B} - \tfrac{1}{2}\bbb{A} \overset{(b)}{\succeq} 2\bbb{B} - \bbb{A} \overset{(c)}{=} \bbb{B} - \bbb{C} \overset{(d)}{\succeq} \delta \bbb{I}
\eeq
\noi where $(a)$ uses the definition of $\bbb{Q}$ in (\ref{eq:dis:ok}); step $(b)$ uses $\tfrac{1}{2} \bbb{A} \preceq \bbb{A}$; step $(c)$ uses the fact that $\bbb{A}=\bbb{B}+\bbb{C}$; step $(d)$ uses the fact that $\forall \bbb{z},~\bbb{z}^T(\bbb{B}-\bbb{C})\bbb{z}\geq \delta\|\bbb{z}\|_2^2$, which is due to (\ref{eq:upperbound:0}).
\noi Then we derive the following inequalities:
\begin{equation}\begin{split} \label{eq:strong:negative}
\bbb{G}^k =  \alpha^k (\tfrac{1}{2}\alpha^k \bbb{P} -\bbb{Q}) &\overset{(a)}{\preceq } \alpha^k (\alpha^k \bbb{B}^T\bbb{B} - \delta \bbb{I} ) \overset{(b)}{\prec}\bbb{0}\nn\\
\end{split}\end{equation}
\noi where step $(a)$ uses (\ref{eq:Q:delta}); step $(b)$ uses the choice that $0<\alpha^k<\frac{\delta}{\|\bbb{B}^T\bbb{B}\|}$.


(iii) We define $\bbb{v}\triangleq\bbb{y}^k-\bbb{x}^k$. Minimizing the right-hand side of (\ref{eq:dis:ok}) over the variable $\alpha$, we obtain (\ref{eq:dis:final2}).
\begin{equation}\begin{split} \label{eq:dis:final2}
\alpha^* = \arg \min_{\alpha}~\psi(\alpha) \triangleq \tfrac{1}{2}(\|\bbb{v}\|^2_{\bbb{P}}) \alpha^2 - (\|\bbb{v}\|_{\bbb{Q}}^2) \alpha.\\
\end{split}\end{equation}
\noi Setting the gradient of quadratic function $\psi(\alpha)$ to zero, we obtain the optimal solution for $\alpha^* = { \|\bbb{v}\|_{\bbb{Q}}^2 }/\|\bbb{v}\|_{\bbb{P}}^2$. We obtain $\psi(\alpha^*) =
\tfrac{1}{2} \|\bbb{v}\|_{\bbb{P}}^2 \cdot \frac{ \|\bbb{v}\|_{\bbb{Q}}^2 }{\|\bbb{v}\|_{\bbb{P}}^2} \cdot \frac{ \|\bbb{v}\|_{\bbb{Q}}^2 }{\|\bbb{v}\|_{\bbb{P}}^2} - \|\bbb{v}\|_{\bbb{Q}}^2\cdot \frac{ \|\bbb{v}\|_{\bbb{Q}}^2 }{\|\bbb{v}\|_{\bbb{P}}^2} =   -\tfrac{1}{2}\tfrac{\|\bbb{v}\|_{\bbb{Q}}^4}{ \|\bbb{v}\|_{\bbb{P}}^2}$. Therefore, we have $\psi(\alpha^*) \leq - \tfrac{\delta^2 \|\bbb{v}\|_2^2 }{4 \|\bbb{BB}^T\|}$.

\end{proof}

\end{theorem}

\noi \textbf{Remarks.} \textbf{(i)} The iterated solutions generated by Algorithm \ref{alg:gmsa:C} satisfy the monotone/contractive property. Therefore, the convergence properties in Theorem \ref{theorem:strong} are stronger than the results in Theorem \ref{theorem:general:rate} and Theorem \ref{theorem:general:rate2}. \textbf{(ii)} There are two methods to decide the value of $\alpha^k$ in Algorithm \ref{alg:gmsa:C}. One method is to use the global constant as indicated in part (ii) in Theorem \ref{theorem:strong}, and another method is to use a local constant parameter as shown in part (iii) in Theorem \ref{theorem:strong}. We remark that a local constant is more desirable in practice since it provides a better estimation of the local structure of the problem and does not require any sparse eigenvalue solver.


\subsection{Acceleration via Richardson Extrapolation} \label{sect:ext:acc}
This subsection discusses using Richardson extrapolation to further accelerate GMSA.

We introduce a parameter $\theta^k \in (0,\infty)$ and consider the following iterative procedure:
\begin{align} \begin{split} \label{eq:acceleration}
\bbb{y}^{k} = ~& \P(\bbb{x}^k)  \\
\bbb{x}^{k+1} =~&  \bbb{x}^{k} + \theta^k (\bbb{y}^{k}-\bbb{x}^{k}).
\end{split}\end{align}
\noi Note that the SOC update rule is not a special case of this formulation. Values of $0<\theta^k <1$ are often used to help establish convergence of the iterative procedure, while values of $\theta^k >1$ are used to speed up convergence of a slow-converging procedure which is also known as Richardson extrapolation \cite{Varga1962Matrix}. Such strategy is closely related to Nesterov's extrapolation acceleration strategy \cite{beck2009fast,nesterov2013introductory}.


The following proposition provides important insights on how to choose the parameter $\theta^k$.

\begin{proposition}
We define:
\beq\label{eq:optimal:theta}
\min_{\theta}~ \varphi (\theta) &\triangleq& \|\bbb{x}^{k+1}-\bbb{x}^*\|_2^2 - \|\bbb{x}^k - \bbb{x}^*\|_2^2~\\
&=&\|\bbb{x}^{k} + \theta^k (\bbb{y}^k - \bbb{x}^{k})-\bbb{x}^*\|_2^2 - \|\bbb{x}^k - \bbb{x}^*\|_2^2 \nn
\eeq

\noi The optimal solution of (\ref{eq:optimal:theta}) can be computed as $\theta^*=\frac{\la \bbb{x}^{k}-\bbb{x}^{*},~\bbb{x}^{k}-\bbb{y}^{k} \ra}{\|\bbb{x}^{k}-\bbb{y}^{k}\|_2^2}$. In addition, if $\la \bbb{x}^{k}-\bbb{x}^{*},~\bbb{x}^{k}-\bbb{y}^{k} \ra \neq 0$, there exists a constant $0<\nu<1$ such that $\|\bbb{x}^{k+1}-\bbb{x}^*\|_2^2 \leq \nu^{k+1}\|\bbb{x}^{0}-\bbb{x}^*\|_2^2$.

\begin{proof}

From (\ref{eq:optimal:theta}), we have:
\beq \label{eq:jjjj}
\varphi (\theta)&=&\|\theta (\bbb{y}^k-\bbb{x}^k)+\bbb{x}^k-\bbb{x}^*\|_2^2 - \|\bbb{x}^k - \bbb{x}^*\|_2^2 \nn\\
&=&\theta^2\| \bbb{y}^k-\bbb{x}^k\|_2^2 + 2\theta\la \bbb{x}^k-\bbb{x}^*,  \bbb{y}^k-\bbb{x}^k\ra~~~~~~~~
\eeq

\noi Setting the gradient of the quadratic function $\varphi (\theta)$ to zero, we obtain: $\theta^*=\frac{\la \bbb{x}^{k}-\bbb{y}^{k+1},~\bbb{x}^{k}-\bbb{y}^{k} \ra}{\|\bbb{x}^{k}-\bbb{y}^{k}\|_2^2}$. Putting the optimal solution $\theta^*$ to (\ref{eq:jjjj}), we obtain: $\varphi (\theta^*)  ~=~ \frac{\la \bbb{x}^{k}-\bbb{x}^{*},~\bbb{x}^{k}-\bbb{y}^{k} \ra^2}{\|\bbb{x}^{k}-\bbb{y}^{k}\|_2^2} ~+ 2 \frac{\la \bbb{x}^{k}-\bbb{x}^{*},~\bbb{x}^{k}-\bbb{y}^{k} \ra}{\|\bbb{x}^{k}-\bbb{y}^{k}\|_2^2} \la \bbb{x}^k-\bbb{x}^*,  \bbb{y}^k-\bbb{x}^k\ra  =~ -\frac{\la \bbb{x}^{k}-\bbb{x}^{*},~\bbb{x}^{k}-\bbb{y}^{k} \ra^2}{\|\bbb{x}^{k}-\bbb{y}^{k}\|_2^2}$. Under the assumption that $\la \bbb{x}^{k}-\bbb{x}^{*},~\bbb{x}^{k}-\bbb{y}^{k} \ra \neq 0$, we have $\|\bbb{x}^{k+1}-\bbb{x}^*\|_2^2 - \|\bbb{x}^k - \bbb{x}^*\|_2^2<0$ for all $k$. There exists a constant such that $\|\bbb{x}^{k+1}-\bbb{x}^*\|_2^2\leq \nu \|\bbb{x}^k - \bbb{x}^*\|_2^2$ for all $k$. Solving this recursive inequality, we obtain: $\|\bbb{x}^{k+1}-\bbb{x}^*\|_2^2 \leq  \nu \|\bbb{x}^k - \bbb{x}^*\|_2^2 \leq \nu^2 \|\bbb{x}^{k-1} - \bbb{x}^*\|_2^2 \leq ...\leq \nu^{k+1} \|\bbb{x}^0 - \bbb{x}^*\|_2^2$.

\end{proof}
\end{proposition}

\noi \textbf{Remarks.} \textbf{(i)} The assumption $\la \bbb{x}^{k}-\bbb{x}^{*},~\bbb{x}^{k}-\bbb{y}^{k} \ra \neq 0$ also implies that $\bbb{x}^{k}$ is not the optimal solution since it holds that $\bbb{x}^{k} \neq \bbb{x}^{*}$ and $\bbb{x}^{k} \neq \bbb{y}^{k}$ when $\bbb{x}^{k}$ is not the optimal solution. \textbf{(ii)} The step size selection strategy is attractive since it guarantees contractive property for $\|\bbb{x}^k - \bbb{x}\|^2_2$. However, it is not practical since the optimal solution $\bbb{x}^*$ is unknown.

In what follows, we consider the following solution. Since $\bbb{y}^k$ is the current available solution which is the closest to $\bbb{x}^*$, we replace $\bbb{x}^*$ with $\bbb{y}^k$, and $k$ with $k-1$. We obtain the follow update rule for $\theta^k$:
\beq
\theta^k = \frac{\la \bbb{x}^{k-1}-\bbb{y}^{k},~\bbb{x}^{k-1}-\bbb{y}^{k-1} \ra}{\|\bbb{x}^{k-1}-\bbb{y}^{k-1}\|_2^2}. \nn
\eeq

\noi We summarize our accelerated generalized matrix splitting algorithm in Algorithm \ref{alg:acc:gmsa}. Note that we set $\theta^k=1$ in the first iteration ($k=0$) and introduce two additional parameters $L$ and $U$ to avoid $\theta^k$ to become arbitrarily small or large. Since the strategy in (\ref{eq:acceleration}) is expected to achieve acceleration when $\theta^k>1$, we set $L=1$ and $U=10$ as the default parameters for Algorithm \ref{alg:acc:gmsa}.

\begin{algorithm} [!h]
\algsize
\caption{\label{alg:acc:gmsa} {\bbb{GMSA-A}: \bbb{G}eneralized \bbb{M}atrix \bbb{S}plitting \bbb{A}lgorithm with Richardson Extrapolation \bbb{A}cceleration for Solving (\ref{eq:main}).}}
\begin{algorithmic}[1]
\STATE Choose suitable parameters $\{\omega,~\epsilon\}$.~Initialize $\bbb{x}^0$.\\
\STATE \text{for $k=0,1,2,...$}\\
\STATE~~~$\bbb{y}^{k} = \P(\bbb{x}^k)$
\STATE~~~\text{if $k=0$}
\STATE~~~\text{~~~~$\theta^k=1$}
\STATE~~~\text{else}
\STATE~~~\text{~~~~$\theta^k=\frac{\la \bbb{x}^{k-1}-\bbb{y}^{k},~\bbb{x}^{k-1}-\bbb{y}^{k-1} \ra}{\|\bbb{x}^{k-1}-\bbb{y}^{k-1}\|_2^2}$}
\STATE~~~\text{~~~~$\theta^k = \min[U,\max(L,\theta^k)]$}
\STATE~~~\text{end~if}
\STATE~~~$\bbb{x}^{k+1} =  \bbb{x}^{k} + \theta^k (\bbb{y}^{k}-\bbb{x}^{k}) $
\STATE \text{end for}\\
\STATE Output $\bbb{x}^{k+1}$\\
\end{algorithmic}
\end{algorithm}

\subsection{When h is Nonconvex} \label{sect:ext:noncvx}
When $h(\bbb{x})$ is nonconvex, our theoretical analysis breaks down in (\ref{eq:fail}) and the exact solution to the triangle proximal operator $\P(\bbb{x}^k)$ in (\ref{eq:subproblem}) cannot be guaranteed. However, our Gaussian elimination procedure in Algorithm \ref{alg:sub} can still be applied. What one needs is to solve a one-dimensional nonconvex subproblem in (\ref{eq:1d:subp}). For example, when $h_j(t)=\lambda |t|_0,~\forall j=1,2,...,n$ (\emph{e.g.} in the case of the $\ell_0$ norm), it has an analytical solution: $t^* = {\tiny \left\{
        \begin{array}{cc}
          -\bbb{w}_j/\bbb{B}_{j,j}, & {\bbb{w}_j^2 > 2\lambda  \bbb{B}_{j,j}} \\
          0, & {\bbb{w}_j^2 \leq 2\lambda  \bbb{B}_{j,j}}
        \end{array}
      \right.}$; when $h_j(t)=\lambda |{t}|^p,~\forall j=1,2,...,n$ and $p<1$, it admits a closed form solution for some special values \cite{xu2012regularization,Cao0X13}, such as $p=\frac{1}{2}$ or $\frac{2}{3}$.

Our generalized matrix splitting algorithm is guaranteed to converge even when ${h}(\cdot)$ is nonconvex. Specifically, we present the following theorem.

\begin{theorem} 
(Proof of Global Convergence when ${h}(\cdot)$ is Nonconvex) We define $\delta\triangleq \epsilon +\frac{1-\omega}{\omega} \min(diag(\bbb{D}))$ and let $\{\omega,~\epsilon\}$ be chosen such that $\delta\in(0,\infty)$. Assuming the nonconvex one-dimensional subproblem in (\ref{eq:1d:subp}) can be solved globally and analytically, we have: (i)
\beq \label{eq:nonconvex:suf:dec}
\textstyle f(\bbb{x}^{k+1}) - f(\bbb{x}^k) \leq - \frac{\delta}{2} \|\bbb{x}^{k+1}-\bbb{x}^k\|_2^2 \leq 0
\eeq
\noi (ii) Algorithm \ref{alg:main} is globally convergent.
\begin{proof}
(i) Due to the optimality of the one-dimensional subproblem in (\ref{eq:1d:subp}), for all $j=1,2,...,n$, we have:
\beq
\textstyle  \tfrac{1}{2}\bbb{B}_{j,j}(\bbb{x}^{k+1}_j)^2 + (\bbb{u}_j +  \sum_{i=1}^{j-1} \bbb{B}_{j,i}\bbb{x}^{k+1}_{i}) \bbb{x}^{k+1}_j +  h(\bbb{x}^{k+1}_j) \nn \\
\leq  \textstyle\tfrac{1}{2}\bbb{B}_{j,j}{t}_j^2 + (\bbb{u}_j +  \sum_{i=1}^{j-1} \bbb{B}_{j,i}\bbb{x}^{k+1}_{i}) {t}_j +  h({t}_j),~\forall {t}_j~~~~\nn
\eeq
\noi Letting ${t}_1=\bbb{x}^k_1,~{t}_2=\bbb{x}^k_2,~...~,{t}_n=\bbb{x}^k_n$, we obtain:
\beq
&\textstyle \tfrac{1}{2} \sum_{i}^n\bbb{B}_{i,i} (\bbb{x}_i^{k+1})^2 + \la \bbb{u} + \bbb{Lx}^{k+1},\bbb{x}^{k+1} \ra  +  h(\bbb{x}^{k+1})\nn\\
&\textstyle \leq\tfrac{1}{2}\sum_{i}^n\bbb{B}_{i,i} (\bbb{x}_i^{k})^2 + \la \bbb{u}+\bbb{Lx}^{k+1},\bbb{x}^{k} \ra +  h(\bbb{x}^{k})\nn
\eeq
\noi Since $\bbb{u} = \bbb{b}+\bbb{C}\bbb{x}^k$, we obtain the following inequality:
\beq
\textstyle f^{k+1} + \tfrac{1}{2}\la \bbb{x}^{k+1},(\tfrac{1}{\omega}\bbb{D}+ \epsilon \bbb{I}+2\bbb{L}-\bbb{A})\bbb{x}^{k+1} + 2 \bbb{C}\bbb{x}^k \ra\nn\\
\leq   f^k + \tfrac{1}{2} \la \bbb{x}^{k},(\tfrac{1}{\omega}\bbb{D}+\epsilon \bbb{I}+2\bbb{C}-\bbb{A}) \bbb{x}^{k} + 2\bbb{L}\bbb{x}^{k+1}\ra~~~~\nn
\eeq
\noi By denoting $\bbb{S} \triangleq \bbb{L}-\bbb{L}^T$ and $\bbb{T}\triangleq(\omega-1)/\omega\bbb{D}-\epsilon \bbb{I}$, we have: $\tfrac{1}{\omega}\bbb{D}+\epsilon \bbb{I}+2\bbb{L}-\bbb{A}=\bbb{T}-\bbb{S}$, $\tfrac{1}{\omega}\bbb{D}+\epsilon \bbb{I}+2\bbb{C}-\bbb{A}=\bbb{S}-\bbb{T}$, and $\bbb{L}-\bbb{C}^T=-\bbb{T}$. Therefore, we have the following inequalities:
\begin{align}
&~\textstyle f^{k+1} - f^k~~~~ \nn\\
\leq&~~  \tfrac{1}{2}\la \bbb{x}^{k+1},\bbb{T}\bbb{x}^{k+1}\ra  + \tfrac{1}{2}\la \bbb{x}^k,\bbb{T}\bbb{x}^k\ra  - \la \bbb{x}^k, \bbb{T}\bbb{x}^{k+1}\ra \nn\\
&~-\tfrac{1}{2}\la \bbb{x}^{k+1}, \bbb{S}\bbb{x}^{k+1}\ra  - \tfrac{1}{2}\la \bbb{x}^k, \bbb{S}\bbb{x}^k\ra \nn\\
\overset{(a)}{=}&~\tfrac{1}{2}\la \bbb{x}^k-\bbb{x}^{k+1}, \bbb{T}(\bbb{x}^k-\bbb{x}^{k+1})\ra\overset{(b)}{\leq}-\tfrac{\delta }{2}\|\bbb{x}^{k+1}-\bbb{x}^{k}\|_2^2\nn
\end{align}
\noi where step ($a$) uses $\la \bbb{x},\bbb{S}\bbb{x}\ra =0~\forall \bbb{x}$, since $\bbb{S}$ is a Skew-Hermitian matrix; step ($b$) uses $\bbb{T} + \delta \bbb{I} \preceq \bbb{0}$, since $\bbb{x}+ \min(-\bbb{x})\leq \bbb{0}~\forall \bbb{x}$. Thus, we obtain the sufficient decrease inequality in (\ref{eq:nonconvex:suf:dec}). 

(ii) Based on the sufficient decrease inequality in (\ref{eq:nonconvex:suf:dec}), we have: $f(\bbb{x}^k)$ is a non-increasing sequence, $\|\bbb{x}^k-\bbb{x}^{k+1}\|\rightarrow 0$, and $f(\bbb{x}^{k+1})<f(\bbb{x}^k)$ if $\bbb{x}^k\neq\bbb{x}^{k+1}$. We note that (\ref{eq:opt:bound0}) can be still applied even ${h}(\cdot)$ is nonconvex. Using the same methodology as in the second part of Theorem \ref{theorem:1}, we obtain that $\nabla q(\bbb{x}^{k+1}) + \partial h(\bbb{x}^{k+1})  \rightarrow \bbb{0}$, which implies the convergence of the algorithm.

Note that guaranteeing $\delta\in(0,\infty)$ can be achieved by simply choosing $\omega\in(0,1)$ and setting $\epsilon$ to a small number.

\end{proof}
\end{theorem}






\subsection{When q is not Quadratic} \label{sect:ext:nonq}

This subsection discusses how to adapt GMSA to solve (\ref{eq:main}) even when $q(\cdot)$ is not quadratic but it is convex and twice differentiable. Following previous work \cite{TsengY09,yuan2014newton}, we keep the nonsmooth function $h(\cdot)$ and approximate the smooth function $q(\cdot)$ around the current solution $\bbb{x}^k$ by its second-order Taylor expansion:
\begin{align}
\begin{split}
\mathcal{Q}(\bbb{y},\bbb{x}^k)~\triangleq& ~h(\bbb{y})+ q(\bbb{\bbb{x}}^k) + \la \nabla q(\bbb{\bbb{x}}^k) , \bbb{y} - \bbb{\bbb{x}}^k\ra + \nn \\
&~\tfrac{1}{2} (\bbb{y} - \bbb{\bbb{x}}^k)^T \nabla^2 q(\bbb{\bbb{x}}^k) (\bbb{y} - \bbb{\bbb{x}}^k)\nn
\end{split}
\end{align}
\noi where $\nabla q(\bbb{\bbb{x}^k})$ and $\nabla^2 q(\bbb{\bbb{x}^k})$ denote the gradient and Hessian of $q(\bbb{x})$ at $\bbb{x}^k$, respectively. In order to generate the next solution that decreases the objective, one can minimize the quadratic model above by solving
\beq \label{eq:newton:exact}
\textstyle \bbb{y}^{k} = \arg\min_{\bbb{y}}~\mathcal{Q}(\bbb{y},\bbb{x}^k).
\eeq
\noi And then one performs line search by the update: $\bbb{x}^{k+1} \Leftarrow \bbb{x}^k + \alpha^k (\bbb{x}^{k+1}-\bbb{x}^k)$ for the greatest descent in objective (as in the damped Newton method). Here, $\alpha^k\in(0,1]$ is the step-size selected by backtracking line search. 

\begin{algorithm} [!t]
\algsize
\caption{\label{alg:non:quad} {Generalized Matrix Splitting Algorithm for Minimizing non-Quadratic Composite Functions in (\ref{eq:main}).}}
\begin{algorithmic}[1]
\STATE Choose suitable parameters $\{\omega,~\epsilon\}$.~Initialize $\bbb{x}^0$.\\
\STATE \text{for $k=0,1,2,...$}\\
\STATE~~~Define $\bbb{A}=\nabla^2 q(\bbb{x}^k)$,~$\bbb{b}=\nabla q(\bbb{x}^k) - \nabla^2 q(\bbb{x}^k)\bbb{x}^k$
\STATE~~~$\bbb{y}^k = \P(\bbb{x}^k;\bbb{A},\bbb{b},h)$
\STATE~~~Define $\bbb{d}^k \triangleq \bbb{y}^k-\bbb{x}^k$  \\
\STATE~~~\label{defdef}Define ${\Delta}^k \triangleq \la \nabla q(\bbb{\bbb{x}}^k) ,\bbb{d}^k \ra + h(\bbb{x}^k+\bbb{d}^k) - h(\bbb{x}^k)$
\STATE~~~Find the largest $\alpha^k \in \{\eta^0,\eta^1,...\}$ such that
\beq \label{eq:armijo:cond}
f(\bbb{x}^k+\alpha^k \bbb{d}^k) \leq f(\bbb{x}^k) + \alpha^k \tau \Delta^k
\eeq
\vspace{-10pt}
\STATE~~~$\bbb{x}^{k+1}=\bbb{x}^k + \alpha^k \bbb{d}^k$
\STATE \text{end for}\\
\STATE Output $\bbb{x}^{k+1}$\\
\end{algorithmic}
\end{algorithm}

In practice, one does not need to solve the Newton approximation subproblem in (\ref{eq:newton:exact}) exactly and one iteration suffices for global convergence. We use $\bbb{x}^{k} = \P(\bbb{x}^k;\bbb{A}^k,\bbb{b}^k,h)$ to denote one iteration of GMSA with the parameter $\{\bbb{A}^k,\bbb{b}^k,h\}$. Note that both $\bbb{A}^k$ and $\bbb{b}^k$ are changing with $k$. We use $\{\bbb{B}^k,\bbb{C}^k,\bbb{D}^k,\bbb{L}^k\}$ to denote the associated matrices of $\bbb{A}^k$ using the same splitting strategy as in (\ref{eq:matrix:dec}). In some situation, we drop the superscript $k$ for simplicity since it can be known from the context. We summarized our algorithm to solve the general convex composite problem in Algorithm \ref{alg:non:quad}.

Algorithm \ref{alg:non:quad} iteratively calls Algorithm \ref{alg:main} for one time to compute the next point $\bbb{y}^k$. Based on the search direction $\bbb{d}^k = \bbb{y}^k-\bbb{x}^k$, we employ Armijo's rule and try step size $\alpha \in \{\eta^0,\eta^1,...\}$ with a constant decrease rate $0<\eta<1$ until we find the smallest $t\in \mathbb{N}$ with $\alpha=\eta^t$ such that $\bbb{x}^k+\alpha^k \bbb{d}^k$ satisfies the sufficient decrease condition. A typical choice for the parameters $\{\eta,~\tau\}$ is $\{0.1,~0.25\}$.

In what follows, we present our convergence analysis for Algorithm \ref{alg:non:quad}. The following lemma is useful in our proof.

\begin{lemma}\label{lemma:tail:bound}
Let $\Delta^k$ be defined in Line \ref{defdef} of Algorithm \ref{alg:non:quad}. It holds that:
$$\Delta^k \leq  -\la \bbb{d}^k, \bbb{B}^k \bbb{d}^k\ra.$$

\begin{proof}
It is not hard to notice that GMSA reduces to the following inclusion problem: 
\beq
\bbb{0} \in \bbb{A}^k \bbb{y}^k + \nabla q(\bbb{x}^k) - \bbb{A}^k\bbb{x}^k + \partial h(\bbb{y}^k) + \bbb{C}^k(\bbb{x}^k-\bbb{y}^k) \nn
\eeq
\noi Using the definition that $\bbb{d}^k \triangleq \bbb{y}^k-\bbb{x}^k$ and $\bbb{A}^k = \bbb{B}^k + \bbb{C}^k$, we have:
\beq \label{eq:hbbb}
\bbb{0} \in  \bbb{B}^k \bbb{d}^k + \nabla q(\bbb{x}^k)  + \partial h(\bbb{y}^k)
\eeq
\noi We further derive the following inequalities: 
\begin{align}
&h(\bbb{x}^k+\bbb{d}^k) - h(\bbb{x}^k) \overset{(a)}{=} h(\bbb{y}^k) - h(\bbb{x}^k) \nn\\
\overset{(b)}{\leq}& \la \bbb{y}^k-\bbb{x}^k, h' \ra,~\forall h' \in \partial h(\bbb{y}^k) \overset{(c)}{=} \la \bbb{d}^k, -\bbb{B}^k \bbb{d}^k - \nabla q(\bbb{x}^k) \ra \nn
\end{align} 
\noi where step $(a)$ uses the fact that $\bbb{x}^k+\bbb{d}^k=\bbb{y}^k$; step $(b)$ uses the convexity of $h(\cdot)$; step $(c)$ uses (\ref{eq:hbbb}) and $\bbb{y}^k-\bbb{x}^k=\bbb{d}^k$. Rearranging terms we finish the proof of this lemma.

\end{proof}

\end{lemma}

\begin{theorem}
\label{theorem:conv}
We define $\delta^k \triangleq ({1}/{\omega}-{1}/{2})\min(diag(\bbb{D}^k))+\epsilon$ and let $\{\omega,~\epsilon\}$ be chosen such that $\delta^k\in(0,\infty)$. Assuming that the gradient of $q(\cdot)$ is $L$-Lipschitz continuous, we have the following results: (i) There always exists a strictly positive constant $\alpha^k$ such that the descent condition in (\ref{eq:armijo:cond}) is satisfied. (ii) The sequence $\{f(\bbb{x}^k)\}_{k=0}^{\infty}$ is nonincreasing and Algorithm \ref{alg:non:quad} is globally convergent.

\begin{proof}

For simplicity, we drop the iteration counter $k$ as it can be inferred from the context. First of all, for all $\bbb{v}\in\mathbb{R}^n$, we have: 
\begin{align} \label{eq:bound:normB}
\bbb{v}^T\bbb{B} \bbb{v} & \overset{(a)}{=}  \bbb{v}^T(\tfrac{1}{2}\bbb{L}+ \tfrac{1}{2}\bbb{L}^T+ \tfrac{1}{2}\bbb{D} + (\tfrac{1}{\omega}-\tfrac{1}{2}) \bbb{D}+\epsilon \bbb{I} )\bbb{v} \nn\\
& \overset{(b)}{=}  \tfrac{1}{2} \bbb{v}^T\bbb{A}\bbb{v} + \bbb{v}^T[(\tfrac{1}{\omega}-\tfrac{1}{2})\bbb{D}+\epsilon \bbb{I} ]\bbb{v} \nn\\
&\overset{(c)}{\geq} 0+ \delta^k\|\bbb{v}\|_2^2
\end{align} 
\noi where step $(a)$ uses the definition of $\bbb{B}$ that $\bbb{B}=\bbb{L} + \frac{1}{\omega} \bbb{D}+\epsilon \bbb{I}$ and the fact that $\bbb{v}^T\bbb{Lv}=\bbb{v}^T\bbb{L}^T\bbb{v}$; step $(b)$ uses $\bbb{A}=\bbb{L}+\bbb{L}^T+\bbb{D}$; step $(c)$ uses the fact that $\bbb{A}$ is positive semidefinite for convex problems.

For any $\alpha\in(0,1]$, we have the following results:
\begin{align}
&~f(\bbb{x}+\alpha \bbb{d}) - f(\bbb{x}) \nn\\
=&~ q(\bbb{x}+\alpha \bbb{d}) - q(\bbb{x}) + h(\bbb{x}+\alpha \bbb{d}) - h(\bbb{x}) \nn\\
\overset{(a)}{\leq}&~ q(\bbb{x}+\alpha \bbb{d}) - q(\bbb{x}) + \alpha [h(\bbb{x}+ \bbb{d}) - h(\bbb{x})] \nn\\
\overset{(b)}{\leq}& ~\tfrac{\alpha^2L}{2} \|\bbb{d}\|_2^2  + \alpha [\la \bbb{d} , \nabla q(\bbb{x}) \ra + h(\bbb{x}+ \bbb{d}) - h(\bbb{x})] \nn\\
\overset{(c)}{=} &~ \tfrac{\alpha^2L}{2} \|\bbb{d}\|_2^2 +  \alpha \Delta^k \overset{(d)}{\leq}  - \tfrac{\alpha^2L}{2} \tfrac{\Delta^k}{\delta^k} + \alpha \Delta^k \nn\\
\overset{(e)}{=} & ~(1-\tfrac{\alpha L}{2\delta^k})\alpha \Delta^k  \leq \tau^k \alpha \Delta^k  \nn
\end{align}
\noi where step $(a)$ uses the fact that $h(\bbb{x} +\alpha \bbb{d} ) = h(\alpha (\bbb{x}+\bbb{d}) + (1-\alpha)\bbb{x}) \leq \alpha h(\bbb{x}+\bbb{d}) + (1-\alpha) h(\bbb{x})$ which is due to the convexity of $h(\cdot)$; step $(b)$ uses the inequality $q(\bbb{y})\leq f(\bbb{x}) +\la \nabla q(\bbb{x}), \bbb{y}-\bbb{x} \ra + \tfrac{L}{2}\|\bbb{y}-\bbb{x}\|_2^2$ with $\bbb{y}=\bbb{x}+\alpha \bbb{d}$; step $(c)$ uses the definition of $\Delta^k$ in Algorithm \ref{alg:non:quad}; step $(d)$ uses Lemma (\ref{lemma:tail:bound}) that $\Delta^k\leq -\la \bbb{d}^k, \bbb{B}^k\bbb{d}^k \ra \leq -\delta^k \|\bbb{d}^k\|_2^2$ which is due to (\ref{eq:bound:normB}); step $(e)$ uses the choice that $0<\alpha<\min[1,{2\delta^k (1-\tau^k)}/{L}]$.

%
%
%
%

(ii) We obtain the following inequality:
\begin{align} \label{eq:conv}
\forall k,~f(\bbb{x}^{k+1})  - f(\bbb{x}^{k})~\leq ~- \tau^k \alpha^k \|\bbb{d}^k\|_{2}^2 ~~\text{with} ~~\tau^k \alpha^k  >0,\nn
\end{align}
\noi and the sequence $f(\bbb{x}^k)$ is non-increasing. Using the same methodology as in the second part of Theorem \ref{theorem:1}, we have $\bbb{d}^k \rightarrow 0$ as $k\rightarrow \infty$. Therefore, any cluster point of the sequence $\bbb{x}^k$ is a stationary point. Finally, we have $\bbb{d}^k=\bbb{0}$. From (\ref{eq:opt:bound0}) in Lemma (\ref{lemma:opt:ineq}), we obtain: $\bbb{0} = -\bbb{C}(\bbb{x}^k-\P(\bbb{x}^k)) \in \nabla q(\P(\bbb{x}^k)) + \partial h(\P(\bbb{x}^k))$. Therefore, we conclude that $\P(\bbb{x}^k) = \bbb{y}^k = \bbb{x}^k$ is the global optimal solution.

\end{proof}
\end{theorem}



\subsection{Adapting into ADMM Optimization Framework} \label{sect:ext:admm}

This subsection shows how to adapt GMSA into the general optimization framework of ADMM \cite{HeY12,HeY15} to solve the following structured convex problem:
\beq  \label{eq:structured:general:convex2}
\min_{\bbb{x},\bbb{y}}~\tfrac{1}{2}\bbb{x}^T\bbb{Q}\bbb{x} + \bbb{x}^T\bbb{p} + h(\bbb{x}) + r(\bbb{y}),~\bbb{Ex}+\bbb{z}=\bbb{y}
\eeq
\noi where $\bbb{x}\in\mathbb{R}^{n}$ and $\bbb{y}\in\mathbb{R}^{m}$ are decision variables, and $\bbb{Q}\in\mathbb{R}^{n\times n},~\bbb{p}\in\mathbb{R}^n,~\bbb{E}\in\mathbb{R}^{m\times n},~\bbb{z}\in\mathbb{R}^m$ are given. We assume that $\bbb{Q}$ is positive semidefinite and $r(\cdot)$ is simple and convex (may not necessarily be separable) such that its proximal
operator $\text{prox}_{r}(\bbb{a}) = \arg \min_{\bbb{x}} ~\frac{1}{2}\|\bbb{x}-\bbb{a}\|_2^2 + r(\bbb{x})$ can be evaluated analytically. We let $\mathcal{L}: \mathbb{R}^n,~\mathbb{R}^m,~\mathbb{R}^m \mapsto \mathbb{R}$ be the augmented Lagrangian function of (\ref{eq:structured:general:convex2}):
\begin{align}\begin{split}
&\mathcal{L}(\bbb{x},~\bbb{y},~\bbb{\pi}) = \tfrac{1}{2}\bbb{x}^T\bbb{Q}\bbb{x} + \bbb{x}^T\bbb{p} + h(\bbb{x}) + r(\bbb{y})\\
&+ \la \bbb{Ex}+\bbb{z}-\bbb{y},~\bbb{\pi} \ra + \tfrac{\beta}{2}\| \bbb{Ex}-\bbb{z}+\bbb{y}\|_2^2 \nn
\end{split}\end{align}
\noi where $\bbb{\pi}\in\mathbb{R}^m$ is the multiplier associated with the equality constraint $\bbb{Ex}+\bbb{z}=\bbb{y}$, and $\beta>0$ is the penalty parameter.

We summarize our algorithm for solving (\ref{eq:structured:general:convex2}) in Algorithm \ref{alg:admm}. The algorithm optimizes for a set of primal variables at a time and keeps the other primal and dual variables fixed, with the dual variables updating via gradient ascent. For the $\bbb{y}$-subproblem, it admits a closed-form solution. For the $\bbb{x}$-subproblem, since the smooth part is quadratic and the nonsmooth part is separable, it can be solved by Algorithm \ref{alg:main}. Algorithm \ref{alg:admm} is convergent if the $\bbb{x}$-subproblem is solved exactly since it reduces to classical ADMM \cite{HeY12}. We remark that similar to linearized ADMM \cite{HeY12}, Algorithm \ref{alg:admm} also suffices for convergence empirically even if we only solve the $\bbb{x}$-subproblem approximately.

\begin{algorithm} [!t]
\algsize
\caption{\label{alg:admm} {\textbf{GMSA-ADMM}: \textbf{G}eneralized \textbf{M}atrix \textbf{S}plitting \textbf{A}lgorithm-based \textbf{A}lternating \textbf{D}irection \textbf{M}ethod of \textbf{M}ultipliers for Solving (\ref{eq:structured:general:convex2}).}}
\begin{algorithmic}[1]
\STATE Choose $\omega\in(0,2),~ \epsilon\in[0,\infty)$.~Initialize $\bbb{x}^0$.\\
\STATE \text{for $k=0,1,2,...$}\\
\STATE~~~Use Algorithm \ref{alg:main} to solve the following problem:\\
~~~~~~~~~~~~~~~$\bbb{x}^{k+1} =  \min_{\bbb{x}}~\mathcal{L}(\bbb{x},~\bbb{y}^{k},~\bbb{\pi}^k)$
\STATE~~~$\bbb{y}^{k+1} =  \min_{\bbb{y}}~\mathcal{L}(\bbb{x}^{k+1},~\bbb{y},~\bbb{\pi}^k) $
\STATE~~~$\bbb{\pi}^{k+1} =  \bbb{\pi}^{k} + \beta (\bbb{Ex}^{k+1}+\bbb{z}-\bbb{y}^{k+1})$
\STATE \text{end for}\\
\STATE Output $\bbb{x}^{k+1}$\\
\end{algorithmic}
\end{algorithm}

\begin{figure*} [!t]
\centering

\subfloat[ \footnotesize $m=200,~\bbb{x}^0 = \bbb{0}$]{\includegraphics[width=0.25\textwidth,height=0.14\textheight]{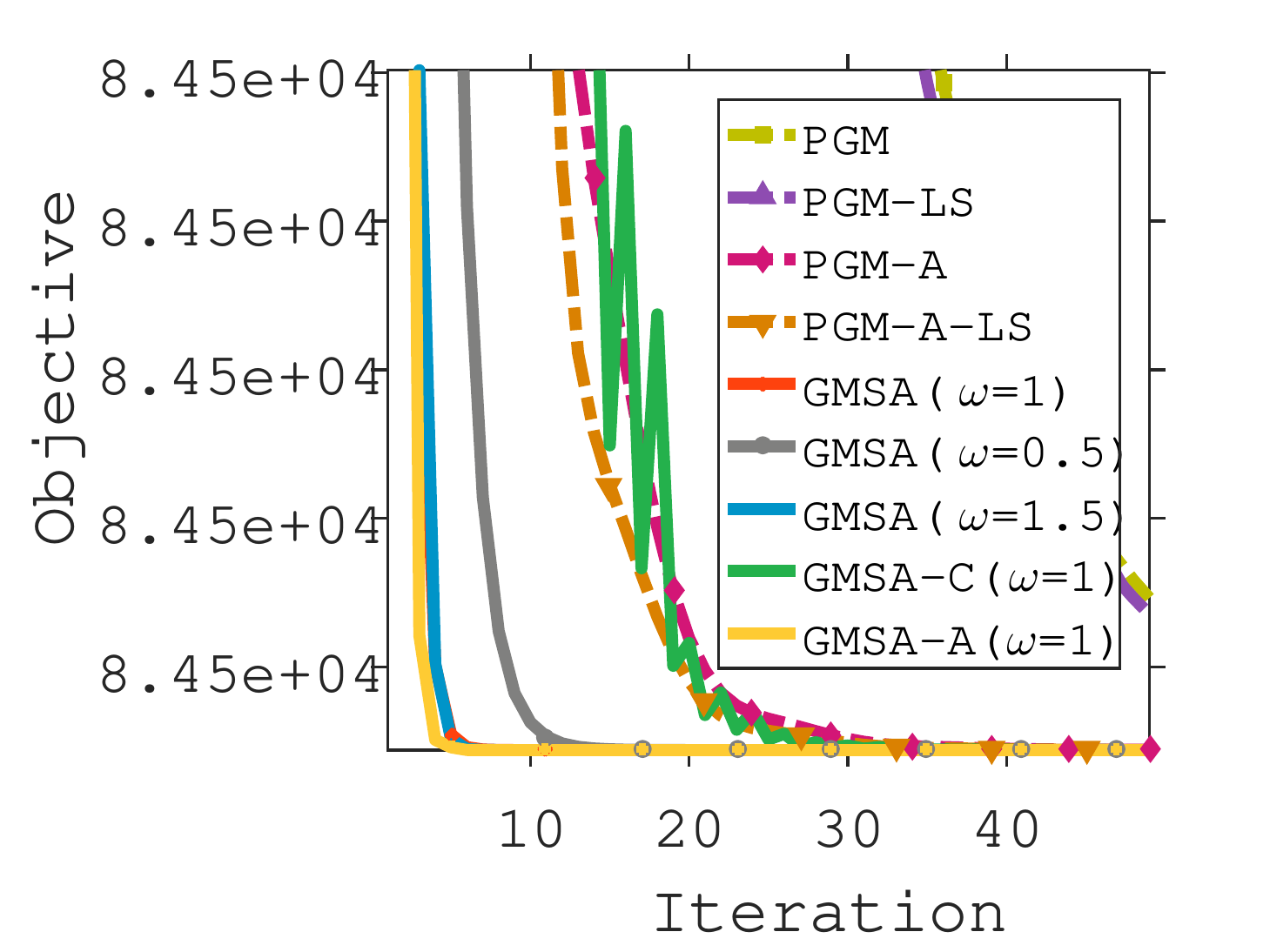}}
\subfloat[ \footnotesize $m=200,~\bbb{x}^0 = rand(n,1)$]{\includegraphics[width=0.25\textwidth,height=0.14\textheight]{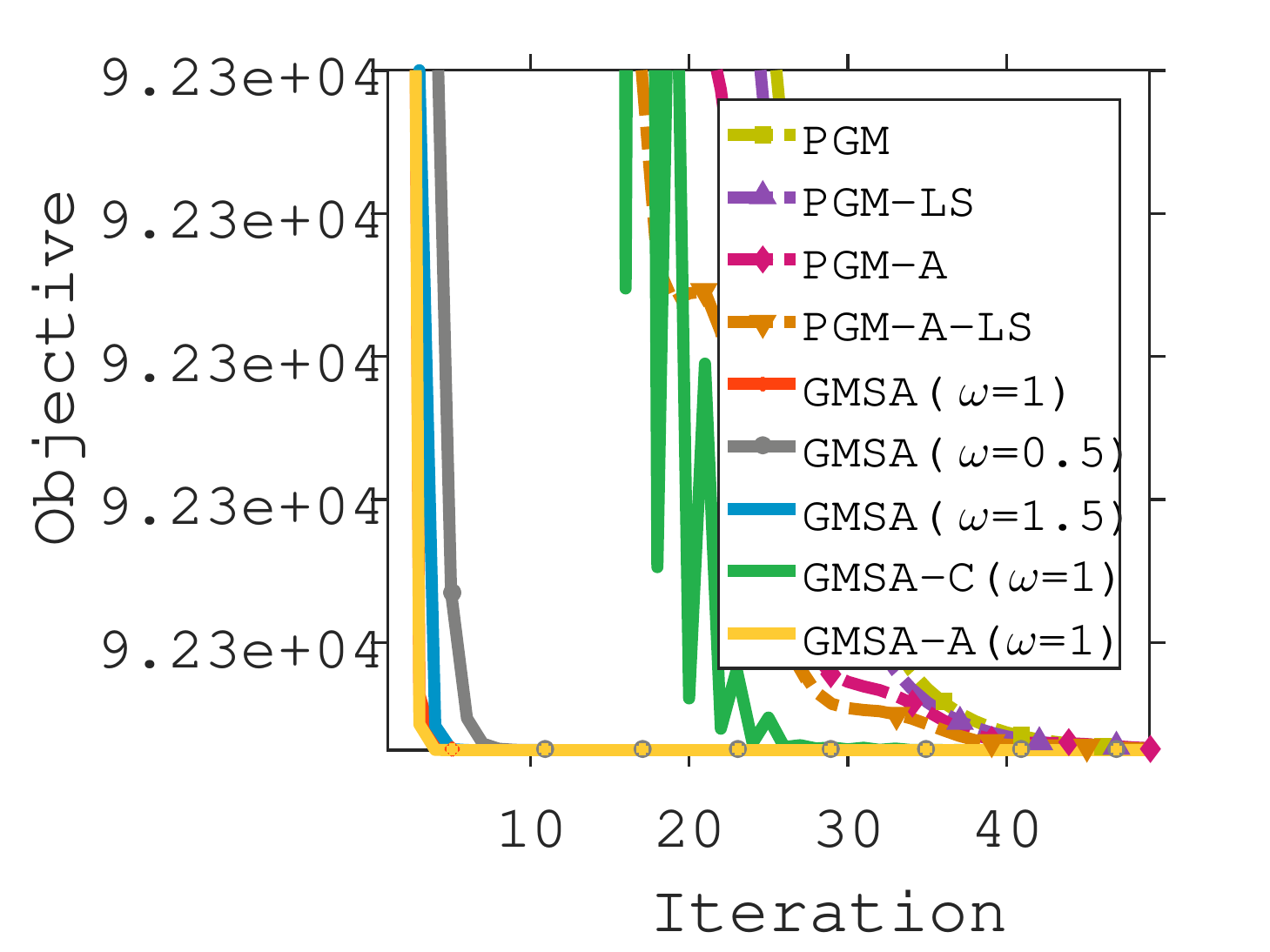}}
\subfloat[ \footnotesize $m=500,~\bbb{x}^0 = \bbb{0}$]{\includegraphics[width=0.25\textwidth,height=0.14\textheight]{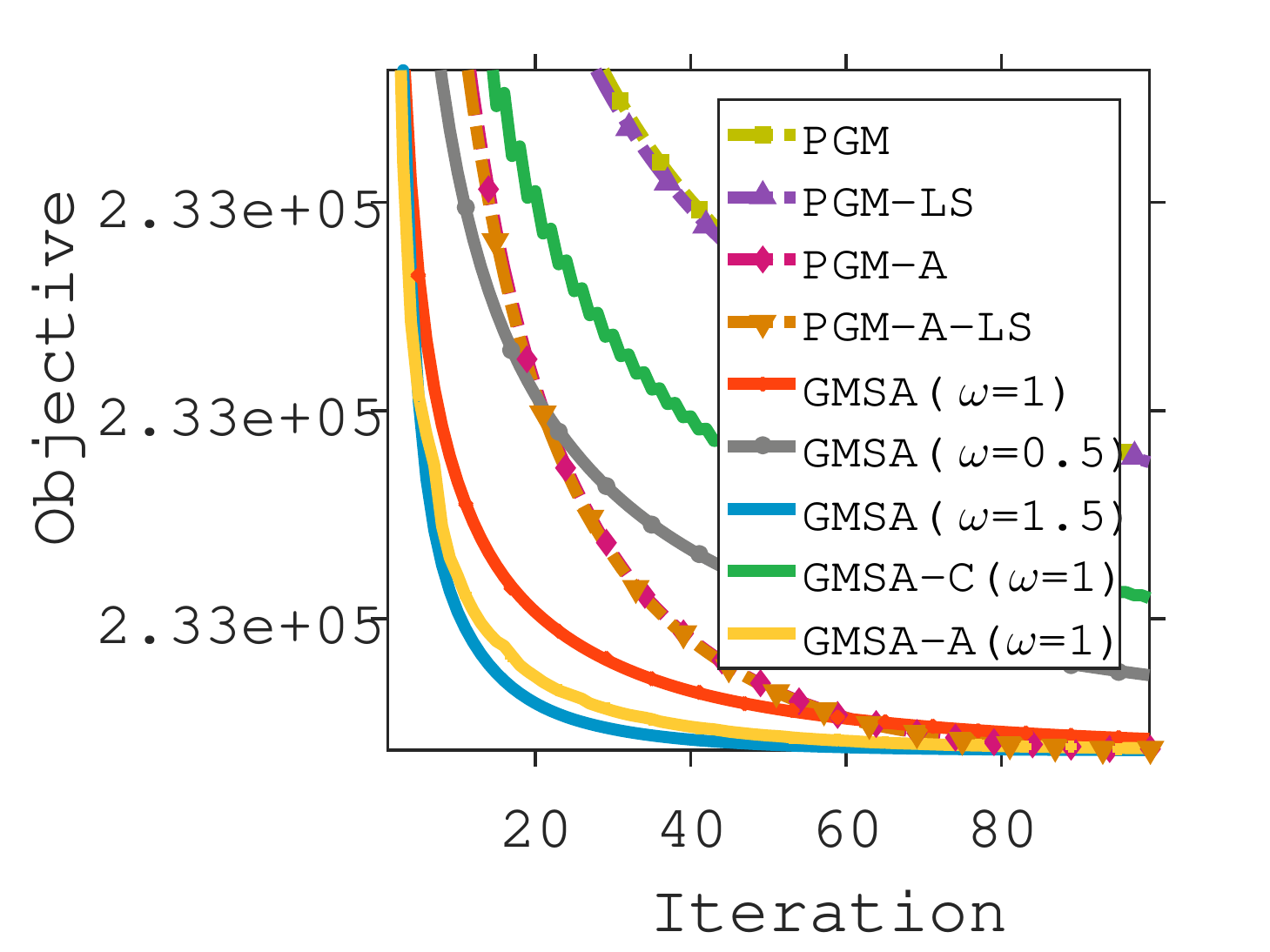}}
\subfloat[ \footnotesize $m=500,~\bbb{x}^0 = rand(n,1)$]{\includegraphics[width=0.25\textwidth,height=0.14\textheight]{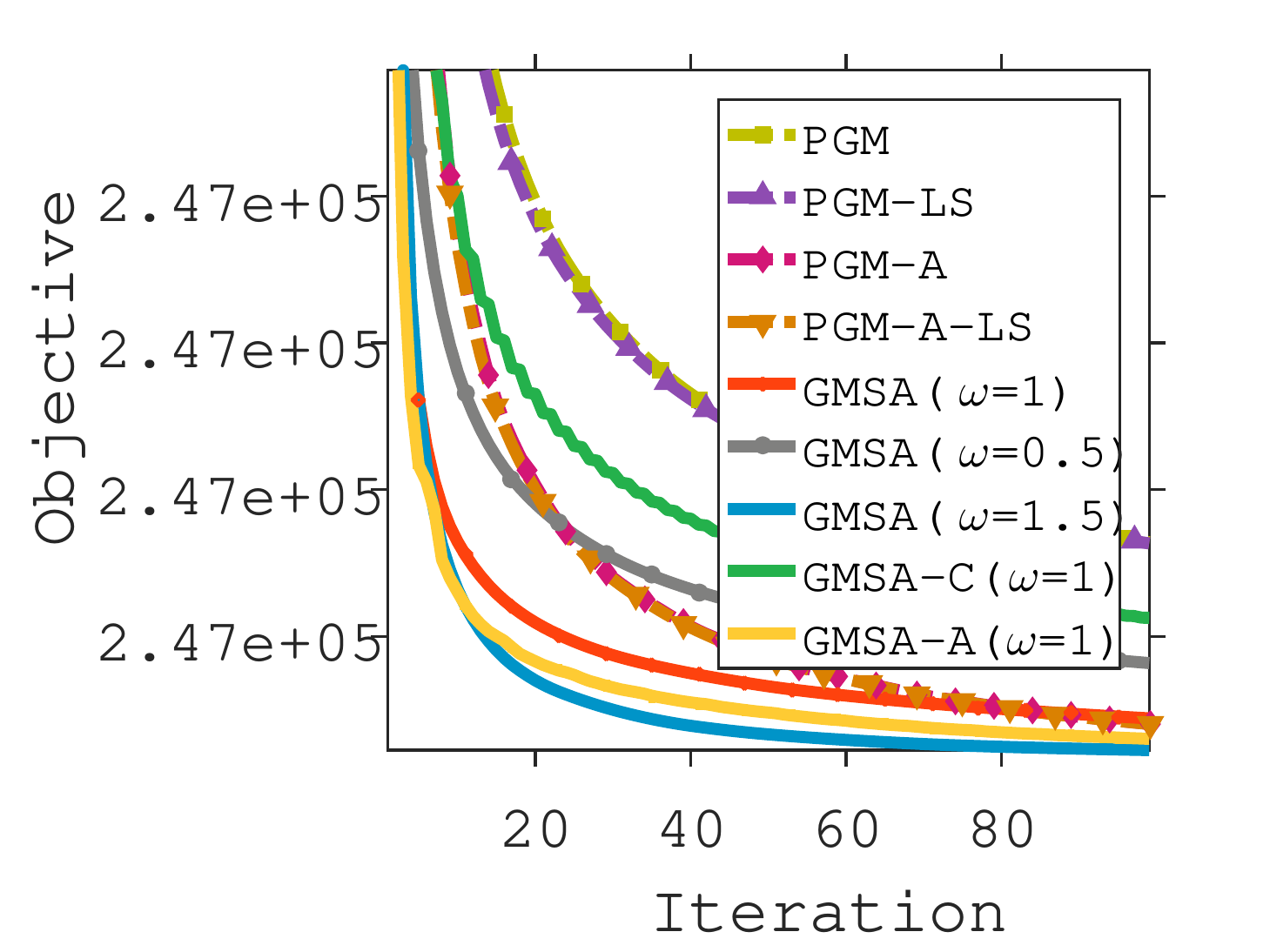}}

\subfloat[\footnotesize $m=200,~\bbb{x}^0 = \bbb{0}$]{\includegraphics[width=0.25\textwidth,height=0.14\textheight]{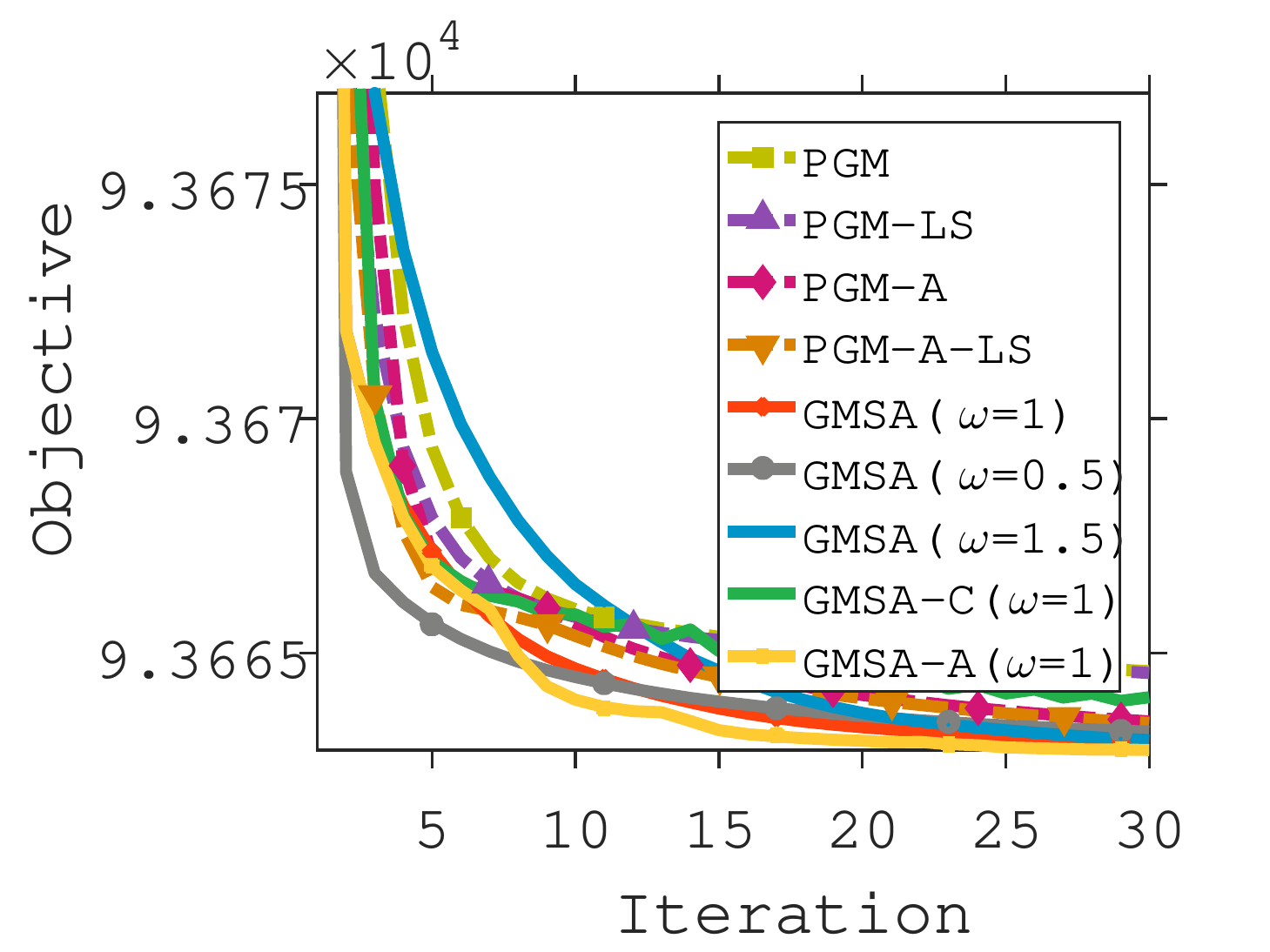}}
\subfloat[ \footnotesize $m=200,~\bbb{x}^0 = rand(n,1)$]{\includegraphics[width=0.25\textwidth,height=0.14\textheight]{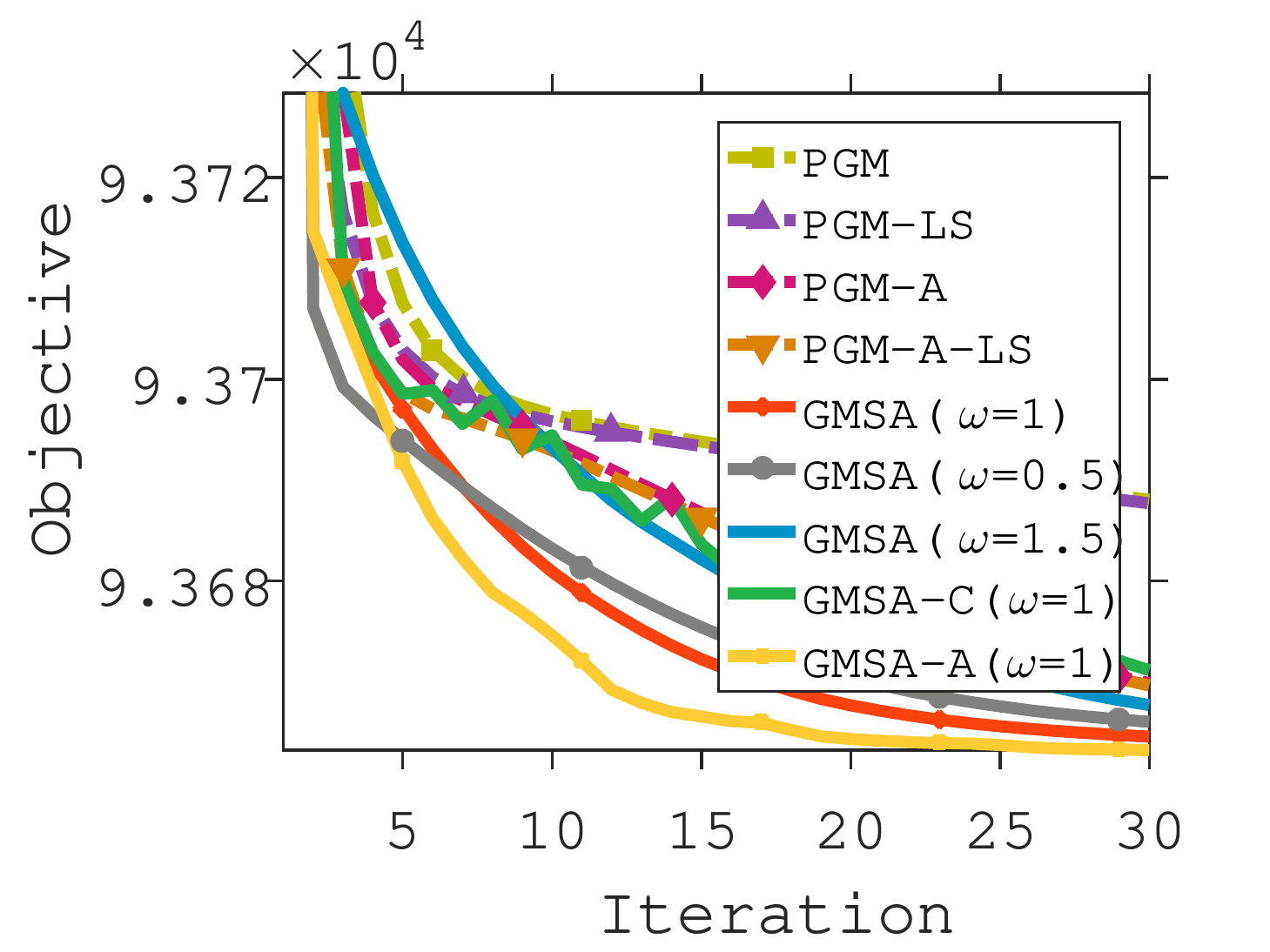}}
\subfloat[ \footnotesize $m=500,~\bbb{x}^0 = \bbb{0}$]{\includegraphics[width=0.25\textwidth,height=0.14\textheight]{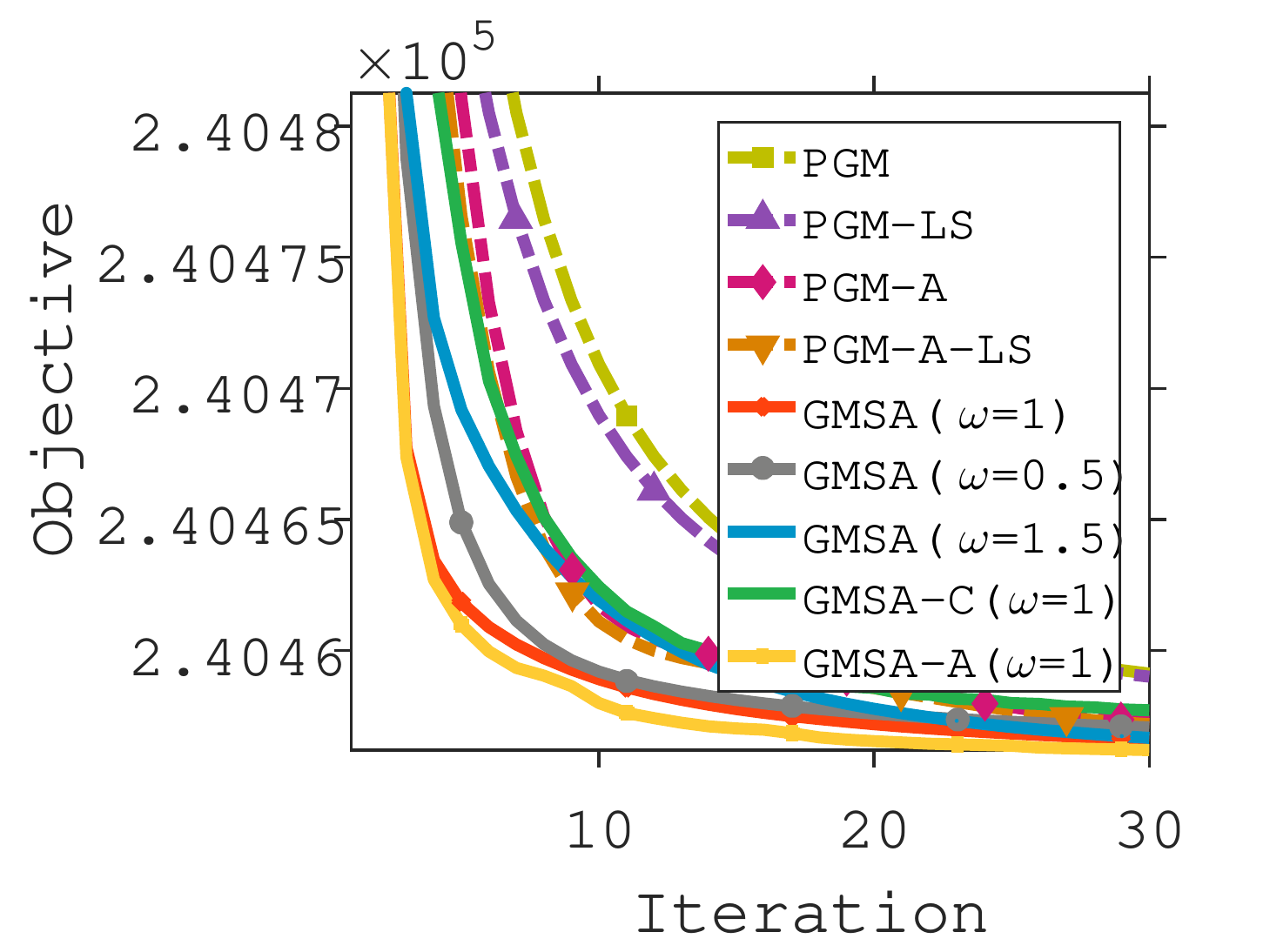}}
\subfloat[ \footnotesize $m=500,~\bbb{x}^0 = rand(n,1)$]{\includegraphics[width=0.25\textwidth,height=0.14\textheight]{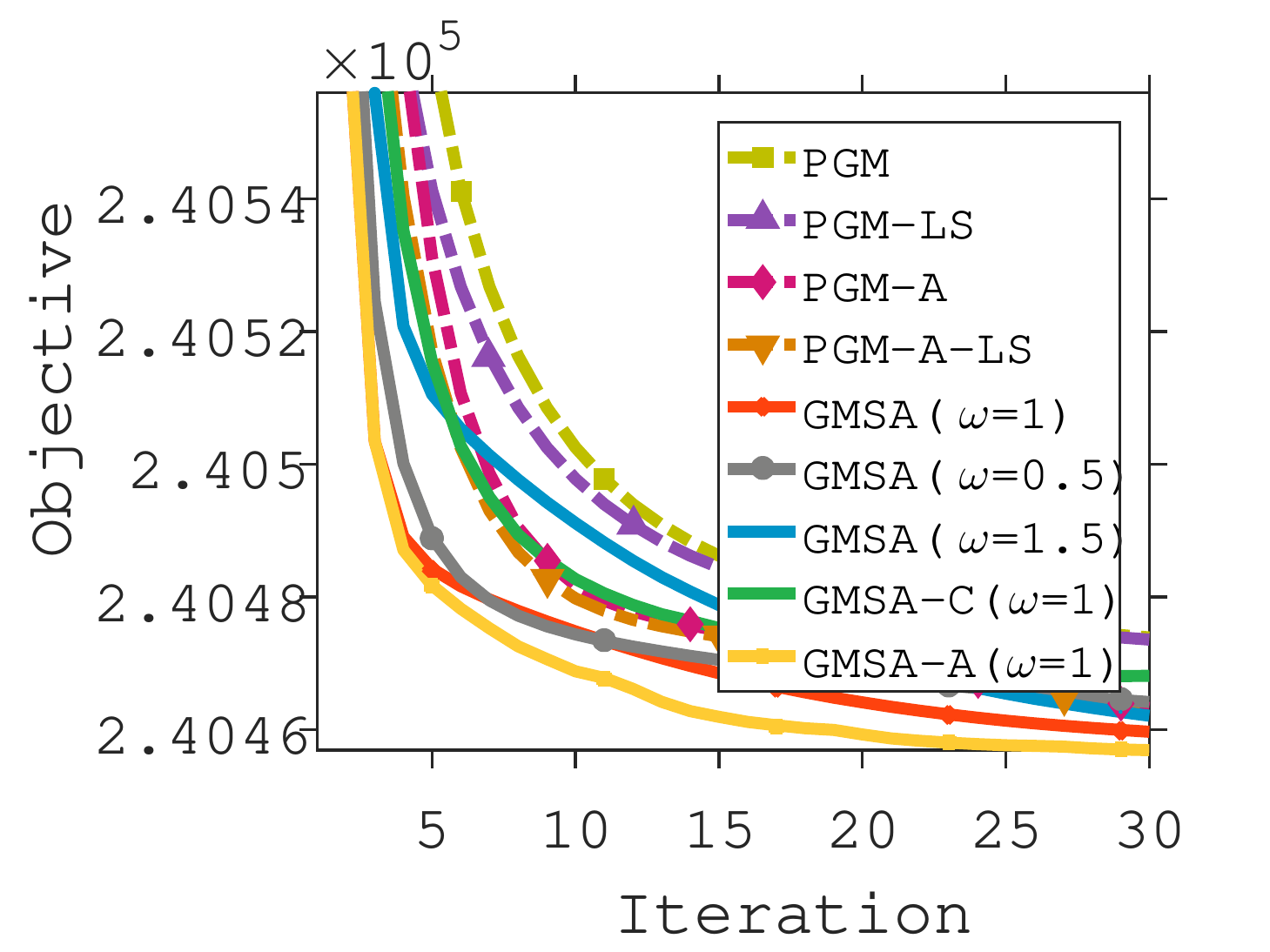}}

\subfloat[ \footnotesize $m=200,~\bbb{x}^0 = \bbb{0}$]{\includegraphics[width=0.25\textwidth,height=0.14\textheight]{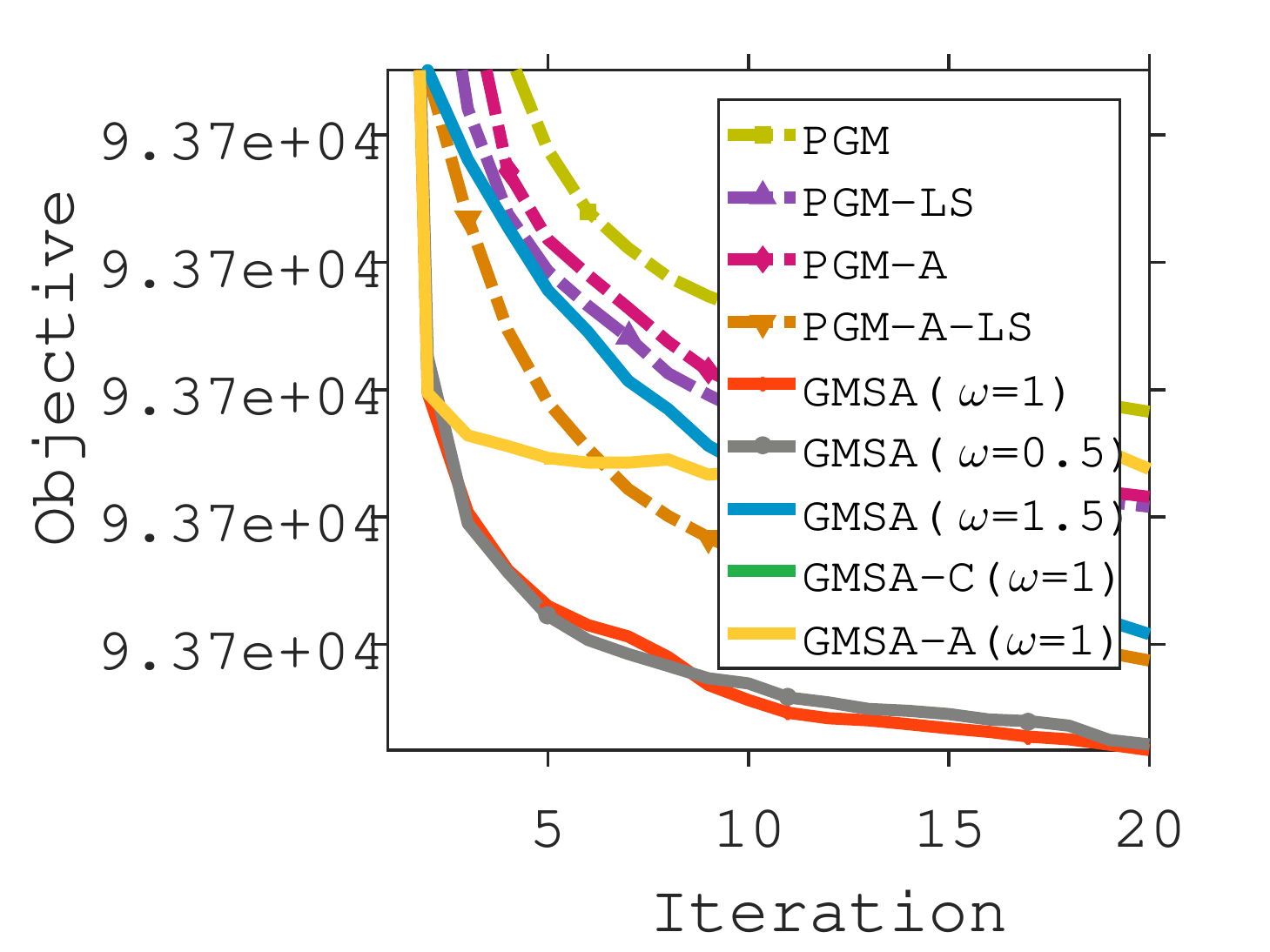}}
\subfloat[ \footnotesize $m=200,~\bbb{x}^0 = rand(n,1)$]{\includegraphics[width=0.25\textwidth,height=0.14\textheight]{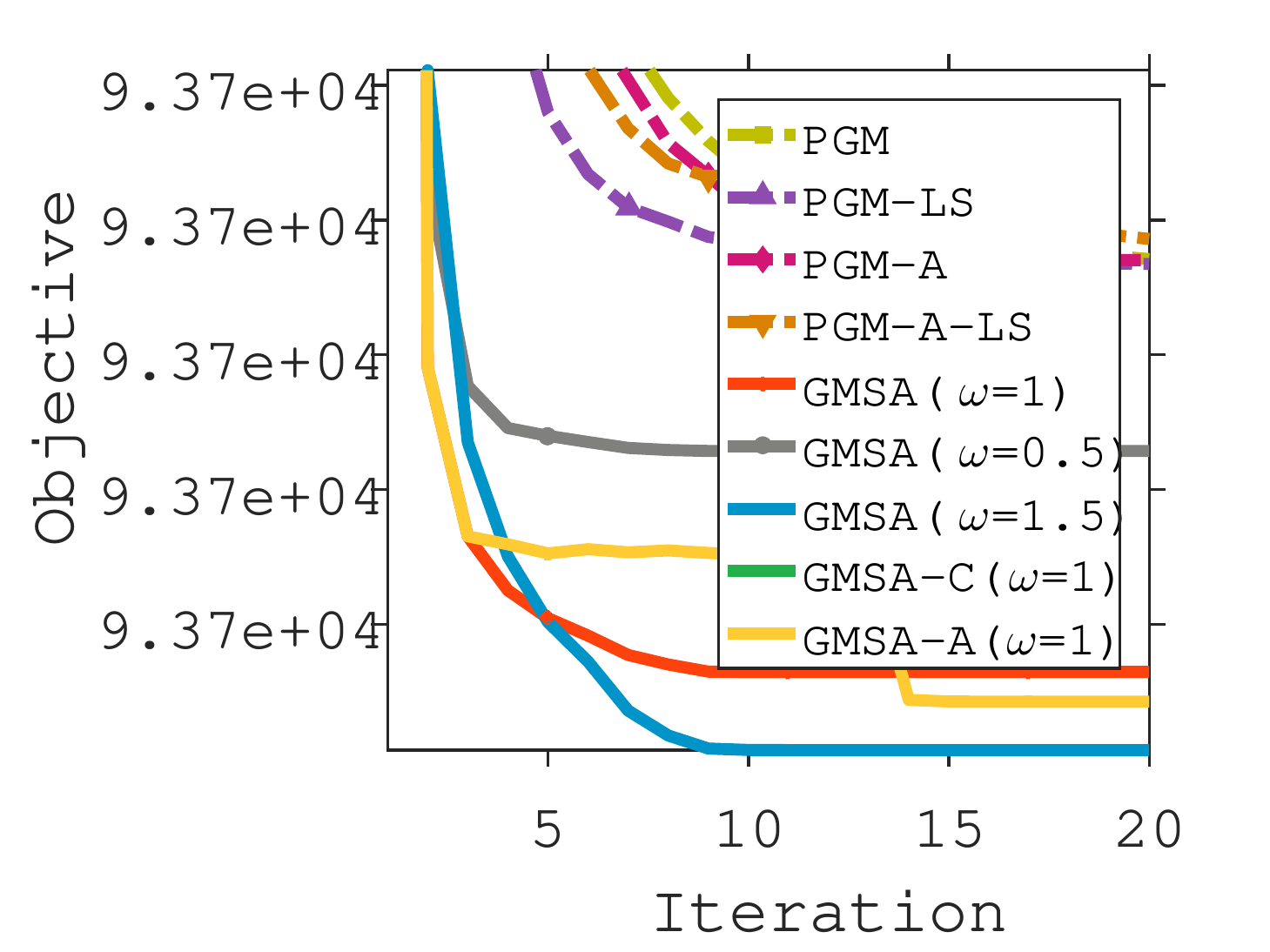}}
\subfloat[ \footnotesize $m=500,~\bbb{x}^0 = \bbb{0}$]{\includegraphics[width=0.25\textwidth,height=0.14\textheight]{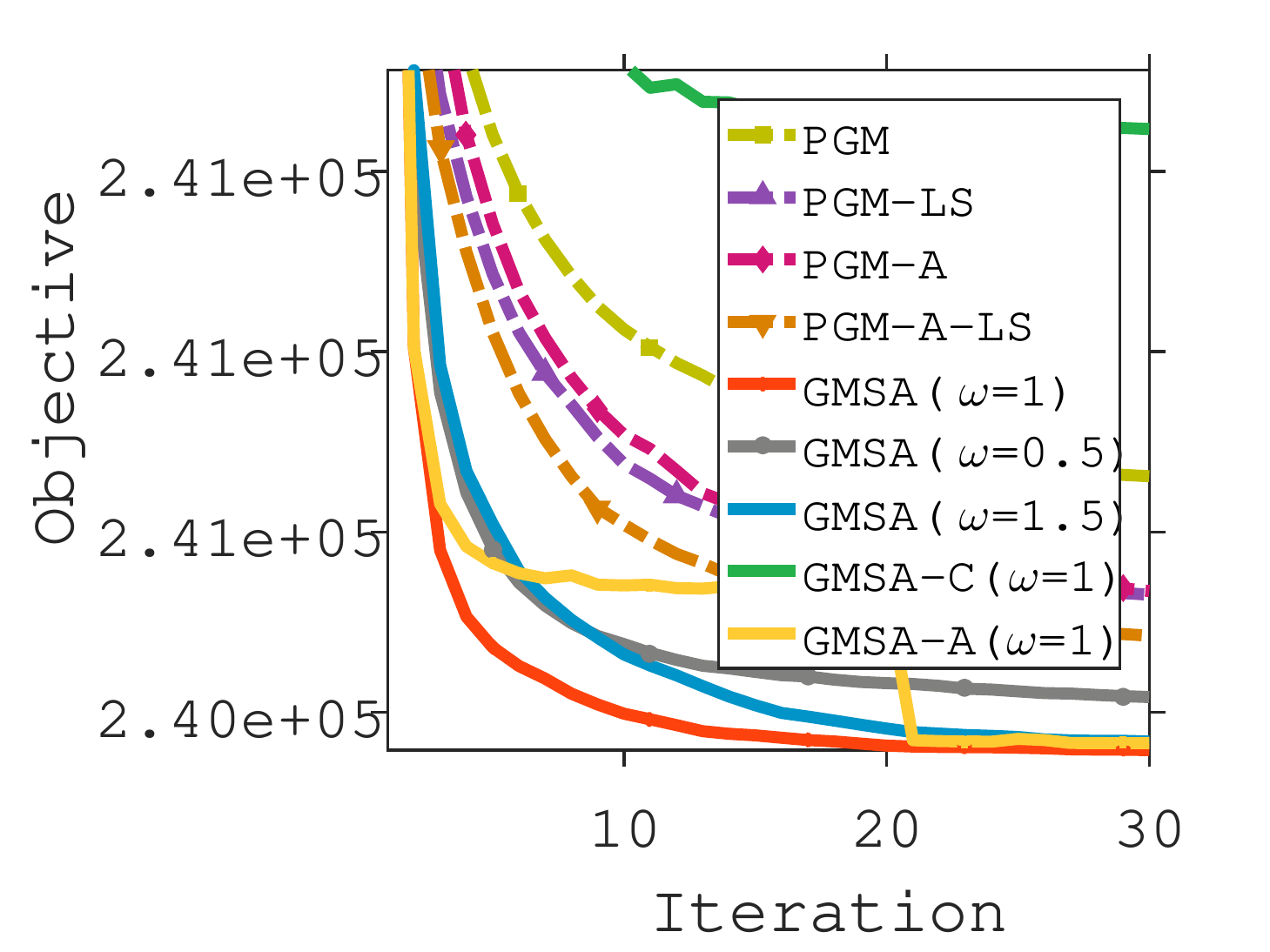}}
\subfloat[ \footnotesize $m=500,~\bbb{x}^0 = rand(n,1)$]{\includegraphics[width=0.25\textwidth,height=0.12\textheight]{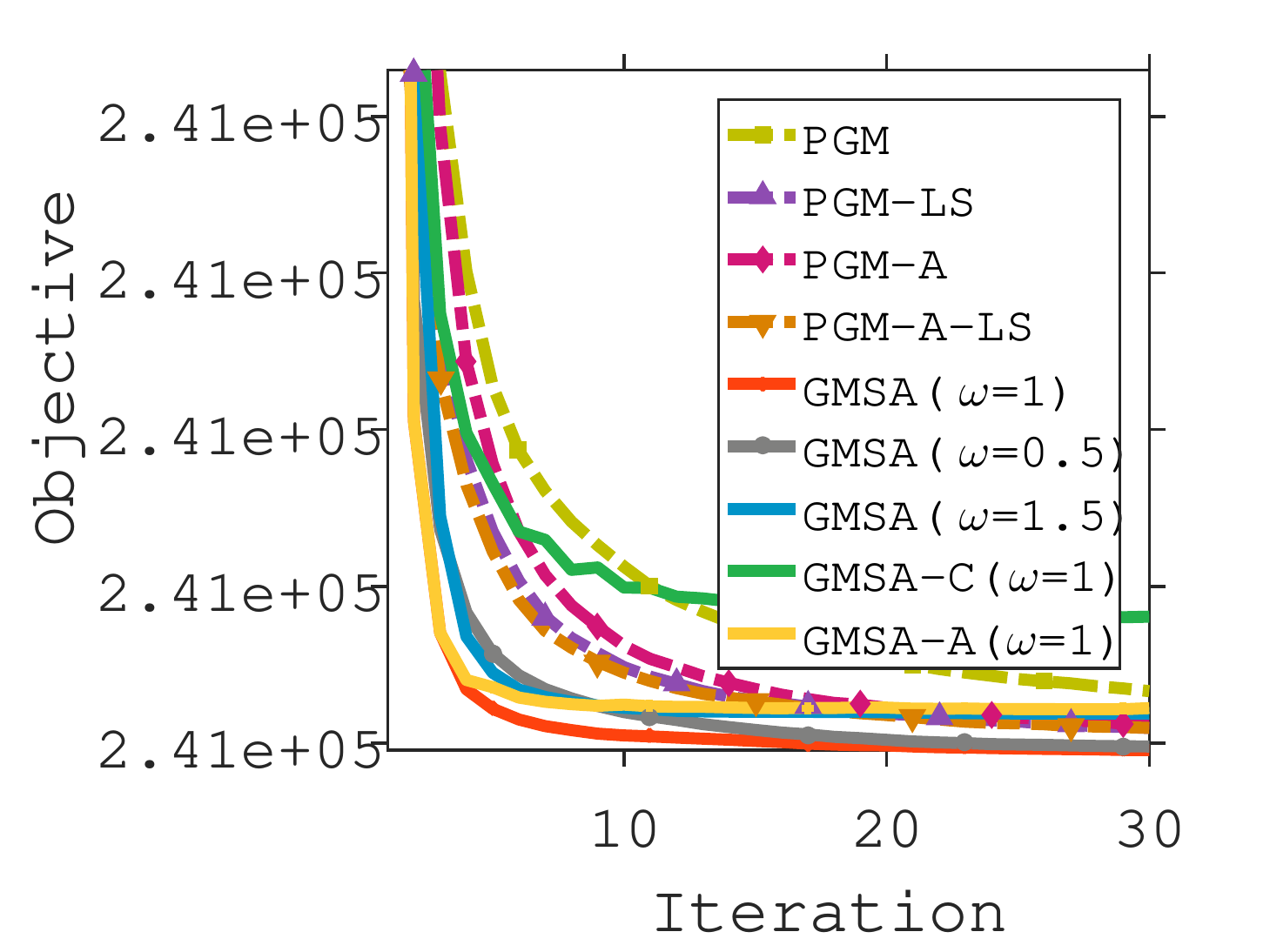}}


\caption{Convergence behavior for solving convex non-negative least squares problem (first row): $\min_{\bbb{x}\geq \bbb{0}}~\tfrac{1}{2}\|\bbb{Cx}-\bbb{d}\|_2^2$, convex $\ell_1$ norm regularized least squares problem (second row): $\min_{\bbb{x}}~\tfrac{1}{2}\|\bbb{Cx}-\bbb{d}\|_2^2+\|\bbb{x}\|_1$, and non-convex $\ell_0$ norm regularized least squares problem: $\min_{\bbb{x}}~\tfrac{1}{2}\|\bbb{Cx}-\bbb{d}\|_2^2+0.1\|\bbb{x}\|_0$ with $\bbb{C}\in\mathbb{R}^{m\times n}$ and $\bbb{d}\in\mathbb{R}^{m}$ being generated from a standard Gaussian distribution. Here $rand(n,1)$ is a function that returns a random vector sampled from a (0-1) uniform distribution. }
\label{fig:sample:l01inf}

\end{figure*}

\subsection{When x is a Matrix}
In many applications (e.g. nonegative matrix factorization and sparse coding), the solutions exist in the matrix form as follows: $\min_{\bbb{X}\in \mathbb{R}^{n\times r}}~~\tfrac{1}{2}tr(\bbb{X}^T\bbb{A}\bbb{X}) + tr(\bbb{X}^T\bbb{R})+ h(\bbb{X})$, where $\bbb{R}\in \mathbb{R}^{n\times r}$. Our matrix splitting algorithm can still be applied in this case. Using the same technique to decompose $\bbb{A}$ as in (\ref{eq:matrix:dec}): $\bbb{A}=\bbb{B}+\bbb{C}$, one needs to replace (\ref{eq:subproblem}) to solve the following nonlinear equation: $\bbb{BZ}^* + \bbb{U} +  \partial h(\bbb{Z}^*) \in 0$, where $\bbb{U}=\bbb{R}+\bbb{C}\bbb{X}^k$. It can be decomposed into $r$ independent components. By updating every column of $\bbb{X}$, the proposed algorithm can be used to solve the matrix problem above. Thus, our algorithm can also make good use of existing parallel architectures to solve the matrix optimization problem.

\section{Experiments}\label{sect:exp}

This section demonstrates the efficiency and efficacy of the proposed Generalized Matrix Splitting Algorithm (GMSA) by considering three important applications: nonnegative matrix factorization (NMF) \cite{lee1999learning,lin2007projected}, $\ell_0$ norm regularized sparse coding \cite{olshausen1996emergence,Quan2016CVPR,lee2006efficient}, and $\ell_1$ norm regularized Danzig selectors. We implement our method in MATLAB on an Intel 2.6 GHz CPU with 8 GB RAM. Only our generalized Gaussian elimination procedure is developed in C and wrapped into the MATLAB code, since it requires an elementwise loop that is quite inefficient in native MATLAB. We consider $\epsilon=0.01$ and $\omega=1$ as our default parameters for GMSA in all our experiments. Some Matlab code can be found in the authors' research webpages.


\subsection{Convergence Behavior of Different Methods}

We demonstrate the convergence behavior of different methods for solving random least squares problems. We compare the following methods. \bbb{(i)} PGM: classical proximal gradient method with constant step size \cite{nesterov2013introductory}; \bbb{(ii)} PGM-LS: classical PGM with line search \cite{beck2009fast}; \bbb{(iii)} PGM-A: accelerated PGM with constant step size \cite{nesterov2013introductory}; \bbb{(iv)} PGM-A-LS: accelerated PGM with line search \cite{beck2009fast,nesterov2013introductory}; \bbb{(v)} GMSA $(\omega=1/0.5/1.5)$: generalized matrix splitting algorithm with varying the parameter $\omega$ described in Algorithm \ref{alg:main}; \bbb{(vi)} GMSA-C: generalized matrix splitting algorithm with correction step described in (\ref{alg:gmsa:C}), where a local constant for computing the step size $\alpha^k$ is used; \bbb{(vii)} GMSA-A: generalized matrix splitting algorithm with Richardson extrapolation acceleration described in (\ref{alg:acc:gmsa}). We report the objective values of the comparing methods for each iteration when they are applied to solve non-negative/$\ell_1$ norm regularized/$\ell_0$ norm regularized least squares problems in Figure \ref{fig:sample:l01inf}. Note that all the methods have the same computational complexity for one iteration.

We have the following observations. \bbb{(i)} GMSA with the default parameters $\omega=1$ and $\theta=0.01$ significantly outperforms proximal gradient method and its variants. \bbb{(ii)} GMSA$(\omega=1.5)$ gives better performance than GMSA $(\omega=0.5)$ for solving non-negative least squares problem but it gives worse performance than GMSA $(\omega=0.5)$ for solving $\ell_1$ norm regularized least squares problem. The choice of the parameter $\omega$ seems to be sensitive to the specific data. \bbb{(iii)} GMSA-C converges slower than GMSA but faster than $\{$PGM,~PGM-LS$\}$. \bbb{(iv)} GMSA-A generally outperforms the other methods in the convex problems. \bbb{(v)} GMSA generally presents the best performance in the nonconvex problems.

Since (i) GMSA with the choice $\epsilon = 0.01$ and $\omega=1$ gives comparable performance to its variants, and (ii) GMSA-A is not necessarily a monotonic algorithm (although it achieves acceleration upon GMSA), we only report the results for GMSA with $\epsilon = 0.01$ and $\omega=1$ in our following experiments.

\begin{table*}[!th]
\scalebox{0.66}{\begin{tabular}{cc}
\begin{tabular}{|c|c|c|c|c|c|c|c|}
\hline
\multicolumn{8}{|>{\columncolor{mycyana}}c|}{\centering time limit=20} \\
\hline
data  & n & \cite{lin2007projected} & \cite{kim2011fast} &  \cite{kim2011fast}  & \cite{guan2012nenmf} & \cite{hsieh2011fast} & [ours] \\
  & & PG & AS &  BPP  & APG  & CGD & GMSA \\
\hline
20news &  20 &5.001e+06 & 2.762e+07 & 8.415e+06 & \cthree{4.528e+06} & \ctwo{4.515e+06} & \cone{4.506e+06} \\
20news &  50 &5.059e+06 & 2.762e+07 & 4.230e+07 & \cthree{3.775e+06} & \ctwo{3.544e+06} & \cone{3.467e+06} \\
20news &  100 &6.955e+06 & 5.779e+06 & 4.453e+07 & \ctwo{3.658e+06} & \cthree{3.971e+06} & \cone{2.902e+06} \\
20news &  200 &7.675e+06 & \ctwo{3.036e+06} & 1.023e+08 & \cthree{4.431e+06} & 3.573e+07 & \cone{2.819e+06} \\
20news &  300 &\cthree{1.997e+07} & 2.762e+07 & 1.956e+08 & \ctwo{4.519e+06} & 4.621e+07 & \cone{3.202e+06} \\
COIL &  20 &2.004e+09 & 5.480e+09 & 2.031e+09 & \cone{1.974e+09} & \cthree{1.976e+09} & \ctwo{1.975e+09} \\
COIL &  50 &1.412e+09 & 1.516e+10 & 6.962e+09 & \cthree{1.291e+09} & \ctwo{1.256e+09} & \cone{1.252e+09} \\
COIL &  100 &2.960e+09 & 2.834e+10 & 3.222e+10 & \cthree{9.919e+08} & \ctwo{8.745e+08} & \cone{8.510e+08} \\
COIL &  200 &3.371e+09 & 2.834e+10 & 5.229e+10 & \cthree{8.495e+08} & \ctwo{5.959e+08} & \cone{5.600e+08} \\
COIL &  300 &3.996e+09 & 2.834e+10 & 1.017e+11 & \cthree{8.493e+08} & \ctwo{5.002e+08} & \cone{4.956e+08} \\
TDT2 &  20 &1.597e+06 & 2.211e+06 & 1.688e+06 & \cone{1.591e+06} & \cthree{1.595e+06} & \ctwo{1.592e+06} \\
TDT2 &  50 &1.408e+06 & 2.211e+06 & 2.895e+06 & \cthree{1.393e+06} & \ctwo{1.390e+06} & \cone{1.385e+06} \\
TDT2 &  100 &1.300e+06 & 2.211e+06 & 6.187e+06 & \ctwo{1.222e+06} & \cthree{1.224e+06} & \cone{1.214e+06} \\
TDT2 &  200 &1.628e+06 & 2.211e+06 & 1.791e+07 & \ctwo{1.119e+06} & \cthree{1.227e+06} & \cone{1.079e+06} \\
TDT2 &  300 &1.915e+06 & \cthree{1.854e+06} & 3.412e+07 & \ctwo{1.172e+06} & 7.902e+06 & \cone{1.066e+06} \\
\hline  \end{tabular}  \hspace{-5.1pt} \begin{tabular}{|c|c|c|c|c|c|c|c|}
  \hline
  \multicolumn{8}{|>{\columncolor{mycyanb}}c|}{\centering time limit=30} \\
    \hline
  data  & n & \cite{lin2007projected} & \cite{kim2011fast} &  \cite{kim2011fast}  & \cite{guan2012nenmf} & \cite{hsieh2011fast} & [ours] \\
    & & PG & AS &  BPP  & APG  & CGD & GMSA \\
  \hline
20news &  20 &4.716e+06 & 2.762e+07 & 7.471e+06 & \cthree{4.510e+06} & \ctwo{4.503e+06} & \cone{4.500e+06} \\
20news &  50 &4.569e+06 & 2.762e+07 & 5.034e+07 & \cthree{3.628e+06} & \ctwo{3.495e+06} & \cone{3.446e+06} \\
20news &  100 &6.639e+06 & 2.762e+07 & 4.316e+07 & \cthree{3.293e+06} & \ctwo{3.223e+06} & \cone{2.817e+06} \\
20news &  200 &\cthree{6.991e+06} & 2.762e+07 & 1.015e+08 & \ctwo{3.609e+06} & 7.676e+06 & \cone{2.507e+06} \\
20news &  300 &\cthree{1.354e+07} & 2.762e+07 & 1.942e+08 & \ctwo{4.519e+06} & 4.621e+07 & \cone{3.097e+06} \\
COIL &  20 &1.992e+09 & 4.405e+09 & 2.014e+09 & \cone{1.974e+09} & \cthree{1.975e+09} & \ctwo{1.975e+09} \\
COIL &  50 &1.335e+09 & 2.420e+10 & 5.772e+09 & \cthree{1.272e+09} & \ctwo{1.252e+09} & \cone{1.250e+09} \\
COIL &  100 &2.936e+09 & 2.834e+10 & 1.814e+10 & \cthree{9.422e+08} & \ctwo{8.623e+08} & \cone{8.458e+08} \\
COIL &  200 &3.362e+09 & 2.834e+10 & 4.627e+10 & \cthree{7.614e+08} & \ctwo{5.720e+08} & \cone{5.392e+08} \\
COIL &  300 &3.946e+09 & 2.834e+10 & 7.417e+10 & \cthree{6.734e+08} & \ctwo{4.609e+08} & \cone{4.544e+08} \\
TDT2 &  20 &1.595e+06 & 2.211e+06 & 1.667e+06 & \cone{1.591e+06} & \cthree{1.594e+06} & \ctwo{1.592e+06} \\
TDT2 &  50 &1.397e+06 & 2.211e+06 & 2.285e+06 & \cthree{1.393e+06} & \ctwo{1.389e+06} & \cone{1.385e+06} \\
TDT2 &  100 &1.241e+06 & 2.211e+06 & 5.702e+06 & \ctwo{1.216e+06} & \cthree{1.219e+06} & \cone{1.212e+06} \\
TDT2 &  200 &1.484e+06 & 1.878e+06 & 1.753e+07 & \ctwo{1.063e+06} & \cthree{1.104e+06} & \cone{1.049e+06} \\
TDT2 &  300 &1.879e+06 & 2.211e+06 & 3.398e+07 & \ctwo{1.060e+06} & \cthree{1.669e+06} & \cone{1.007e+06} \\
\hline  \end{tabular}  \\
\begin{tabular}{|c|c|c|c|c|c|c|c|}
  \hline
  \multicolumn{8}{|>{\columncolor{mycyanc}}c|}{\centering time limit=40} \\
    \hline
  data  & n & \cite{lin2007projected} & \cite{kim2011fast} &  \cite{kim2011fast}  & \cite{guan2012nenmf} & \cite{hsieh2011fast} & [ours] \\
    & & PG & AS &  BPP  & APG  & CGD & GMSA \\
  \hline
20news &  20 &4.622e+06 & 2.762e+07 & 7.547e+06 & \cone{4.495e+06} & \cthree{4.500e+06} & \ctwo{4.496e+06} \\
20news &  50 &4.386e+06 & 2.762e+07 & 1.562e+07 & \cthree{3.564e+06} & \ctwo{3.478e+06} & \cone{3.438e+06} \\
20news &  100 &6.486e+06 & 2.762e+07 & 4.223e+07 & \cthree{3.128e+06} & \ctwo{2.988e+06} & \cone{2.783e+06} \\
20news &  200 &6.731e+06 & 1.934e+07 & 1.003e+08 & \ctwo{3.304e+06} & \cthree{5.744e+06} & \cone{2.407e+06} \\
20news &  300 &\cthree{1.041e+07} & 2.762e+07 & 1.932e+08 & \ctwo{3.621e+06} & 4.621e+07 & \cone{2.543e+06} \\
COIL &  20 &1.987e+09 & 5.141e+09 & 2.010e+09 & \cone{1.974e+09} & \cthree{1.975e+09} & \ctwo{1.975e+09} \\
COIL &  50 &1.308e+09 & 2.403e+10 & 5.032e+09 & \cthree{1.262e+09} & \ctwo{1.250e+09} & \cone{1.248e+09} \\
COIL &  100 &2.922e+09 & 2.834e+10 & 2.086e+10 & \cthree{9.161e+08} & \ctwo{8.555e+08} & \cone{8.430e+08} \\
COIL &  200 &3.361e+09 & 2.834e+10 & 4.116e+10 & \cthree{7.075e+08} & \ctwo{5.584e+08} & \cone{5.289e+08} \\
COIL &  300 &3.920e+09 & 2.834e+10 & 7.040e+10 & \cthree{6.221e+08} & \ctwo{4.384e+08} & \cone{4.294e+08} \\
TDT2 &  20 &1.595e+06 & 2.211e+06 & 1.643e+06 & \cone{1.591e+06} & \cthree{1.594e+06} & \ctwo{1.592e+06} \\
TDT2 &  50 &1.394e+06 & 2.211e+06 & 1.933e+06 & \cthree{1.392e+06} & \ctwo{1.388e+06} & \cone{1.384e+06} \\
TDT2 &  100 &1.229e+06 & 2.211e+06 & 5.259e+06 & \ctwo{1.213e+06} & \cthree{1.216e+06} & \cone{1.211e+06} \\
TDT2 &  200 &1.389e+06 & 1.547e+06 & 1.716e+07 & \ctwo{1.046e+06} & \cthree{1.070e+06} & \cone{1.041e+06} \\
TDT2 &  300 &1.949e+06 & 1.836e+06 & 3.369e+07 & \ctwo{1.008e+06} & \cthree{1.155e+06} & \cone{9.776e+05} \\
\hline  \end{tabular} \hspace{-5.1pt} \begin{tabular}{|c|c|c|c|c|c|c|c|}
  \hline
  \multicolumn{8}{|>{\columncolor{mycyand}}c|}{\centering time limit=50} \\
    \hline
  data  & n & \cite{lin2007projected} & \cite{kim2011fast} &  \cite{kim2011fast}  & \cite{guan2012nenmf} & \cite{hsieh2011fast} & [ours] \\
    & & PG & AS &  BPP  & APG  & CGD & GMSA \\
  \hline
20news &  20 &4.565e+06 & 2.762e+07 & 6.939e+06 & \cone{4.488e+06} & \cthree{4.498e+06} & \ctwo{4.494e+06} \\
20news &  50 &4.343e+06 & 2.762e+07 & 1.813e+07 & \cthree{3.525e+06} & \ctwo{3.469e+06} & \cone{3.432e+06} \\
20news &  100 &6.404e+06 & 2.762e+07 & 3.955e+07 & \cthree{3.046e+06} & \ctwo{2.878e+06} & \cone{2.765e+06} \\
20news &  200 &5.939e+06 & 2.762e+07 & 9.925e+07 & \ctwo{3.121e+06} & \cthree{4.538e+06} & \cone{2.359e+06} \\
20news &  300 &\cthree{9.258e+06} & 2.762e+07 & 1.912e+08 & \ctwo{3.621e+06} & 2.323e+07 & \cone{2.331e+06} \\
COIL &  20 &1.982e+09 & 7.136e+09 & 2.033e+09 & \cone{1.974e+09} & \cthree{1.975e+09} & \ctwo{1.975e+09} \\
COIL &  50 &1.298e+09 & 2.834e+10 & 4.365e+09 & \cthree{1.258e+09} & \ctwo{1.248e+09} & \cone{1.248e+09} \\
COIL &  100 &1.945e+09 & 2.834e+10 & 1.428e+10 & \cthree{9.014e+08} & \ctwo{8.516e+08} & \cone{8.414e+08} \\
COIL &  200 &3.362e+09 & 2.834e+10 & 3.760e+10 & \cthree{6.771e+08} & \ctwo{5.491e+08} & \cone{5.231e+08} \\
COIL &  300 &3.905e+09 & 2.834e+10 & 6.741e+10 & \cthree{5.805e+08} & \ctwo{4.226e+08} & \cone{4.127e+08} \\
TDT2 &  20 &1.595e+06 & 2.211e+06 & 1.622e+06 & \cone{1.591e+06} & \cthree{1.594e+06} & \ctwo{1.592e+06} \\
TDT2 &  50 &1.393e+06 & 2.211e+06 & 1.875e+06 & \cthree{1.392e+06} & \ctwo{1.386e+06} & \cone{1.384e+06} \\
TDT2 &  100 &1.223e+06 & 2.211e+06 & 4.831e+06 & \ctwo{1.212e+06} & \cthree{1.214e+06} & \cone{1.210e+06} \\
TDT2 &  200 &1.267e+06 & 2.211e+06 & 1.671e+07 & \ctwo{1.040e+06} & \cthree{1.054e+06} & \cone{1.036e+06} \\
TDT2 &  300 &1.903e+06 & 2.211e+06 & 3.328e+07 & \ctwo{9.775e+05} & \cthree{1.045e+06} & \cone{9.606e+05} \\
\hline  \end{tabular}  \\
\end{tabular}}
\caption{Comparisons of objective values for non-negative matrix factorization for all the compared methods. The $1^{st}$, $2^{nd}$, and $3^{rd}$ best results are colored with \cone{red}, \ctwo{blue} and \cthree{green}, respectively.}
\label{tab:nmf}
\end{table*}

\subsection{Nonnegative Matrix Factorization }

Nonnegative matrix factorization \cite{lee1999learning} is a very useful tool for feature extraction and identification in the fields of text mining and image understanding. It is formulated as the following optimization problem:
\beq
\textstyle \underset{\bbb{W},\bbb{H}}\min~~\frac{1}{2}\|\bbb{Y}-\bbb{WH}\|_F^2 ,~~s.t.~~\bbb{W}\geq 0,~\bbb{H}\geq 0 \nn
\eeq
\noi where $\bbb{W}\in\mathbb{R}^{m\times n}$ and $\bbb{H}\in\mathbb{R}^{n\times d}$. Following previous work \cite{kim2011fast,guan2012nenmf,lin2007projected,hsieh2011fast}, we alternatively minimize the objective while keeping one of the two variables fixed. In each alternating subproblem, we solve a convex nonnegative least squares problem, where our GMSA is used. We conduct experiments on three datasets \footnote{\url{http://www.cad.zju.edu.cn/home/dengcai/Data/TextData.html}} 20news, COIL, and TDT2. The size of the datasets are $18774\times 61188,~7200\times 1024,~9394\times 36771$, respectively. We compare GMSA against the following state-of-the-art methods: (1) Projective Gradient (PG)  \cite{lin2007projected,bertsekas1999nonlinear} that updates the current solution via steep gradient descent and then maps a point back to the bounded feasible region \footnote{\url{https://www.csie.ntu.edu.tw/~cjlin/libmf/}}; (2) Active Set (AS) method \cite{kim2011fast} and (3) Block Principal Pivoting (BPP) method \cite{kim2011fast} \footnote{\url{http://www.cc.gatech.edu/~hpark/nmfsoftware.php}} that iteratively identify an active and passive set by a principal pivoting procedure and solve a reduced linear system; (4) Accelerated Proximal Gradient (APG) \cite{guan2012nenmf} \footnote{\url{https://sites.google.com/site/nmfsolvers/}} that applies Nesterov's momentum strategy with a constant step size to solve the convex sub-problems; (5) Coordinate Gradient Descent (CGD) \cite{hsieh2011fast} \footnote{\url{http://www.cs.utexas.edu/~cjhsieh/nmf/}} that greedily selects one coordinate by measuring the objective reduction and optimizes for a single variable via closed-form update. Similar to our method, the core procedure of CGD is developed in C and wrapped into the MATLAB code, while all other methods are implemented using builtin MATLAB functions.

We use the same settings as in \cite{lin2007projected}. We compare objective values after running $t$ seconds with $t$ varying from 20 to 50. Table \ref{tab:nmf} presents average results of using 10 random initial points, which are generated from a standard normal distribution. While the other methods may quickly lower objective values when $n$ is small ($n=20$), GMSA catches up very quickly and achieves a faster convergence speed when $n$ is large. It generally achieves the best performance in terms of objective value among all the methods.

\subsection{$\ell_0$ Norm Regularized Sparse Coding}
Sparse coding is a popular unsupervised feature learning technique for data representation that is widely used in computer vision and medical imaging. Motivated by recent success in $\ell_0$ norm modeling \cite{yuan2015l0tv,BaoJQS16,Yang2016}, we consider the following $\ell_0$ norm regularized (\emph{i.e.} cardinality) sparse coding problem:
\beq \label{eq:card:sparse:coding}
 \underset{\bbb{W},\bbb{H}}{\min}~~\tfrac{1}{2}\|\bbb{Y}-\bbb{WH}\|_F^2 + \lambda \|\bbb{H}\|_0,~s.t.~\|\bbb{W}(:,i)\|=1,~\forall i,
\eeq with $\bbb{W}\in \mathbb{R}^{m \times n}$ and $\bbb{H}\in \mathbb{R}^{n \times d}$. Existing solutions for this problem are mostly based on the family of proximal point methods \cite{nesterov2013introductory,BaoJQS16}. We compare GMSA with the following methods: (1) Proximal Gradient Method (PGM) with constant step size, (2) PGM with line search, (3) accelerated PGM with constant step size, and (4) accelerated PGM with line search.

We evaluate all the methods for the application of image denoising. Following \cite{aharon2006img,BaoJQS16}, we set the dimension of the dictionary to $n = 256$. The dictionary is learned from $m=1000$ image patches randomly chosen from the noisy input image. The patch size is $8\times 8$, leading to $d=64$. The experiments are conducted on 16 conventional test images with different noise standard deviations $\sigma$. We tune the regularization parameter $\lambda$ and compare the resulting objective values and the Signalto-Noise Ratio (SNR) values for all methods. We do not include the comparison of SNR values here. Interested readers can refer to Section 4.2 of the conference version of this paper \cite{YuanZG17}.

\begin{figure*} [!t]
\centering

\subfloat[ \footnotesize $\lambda=50$]{\includegraphics[width=0.25\textwidth,height=0.14\textheight]{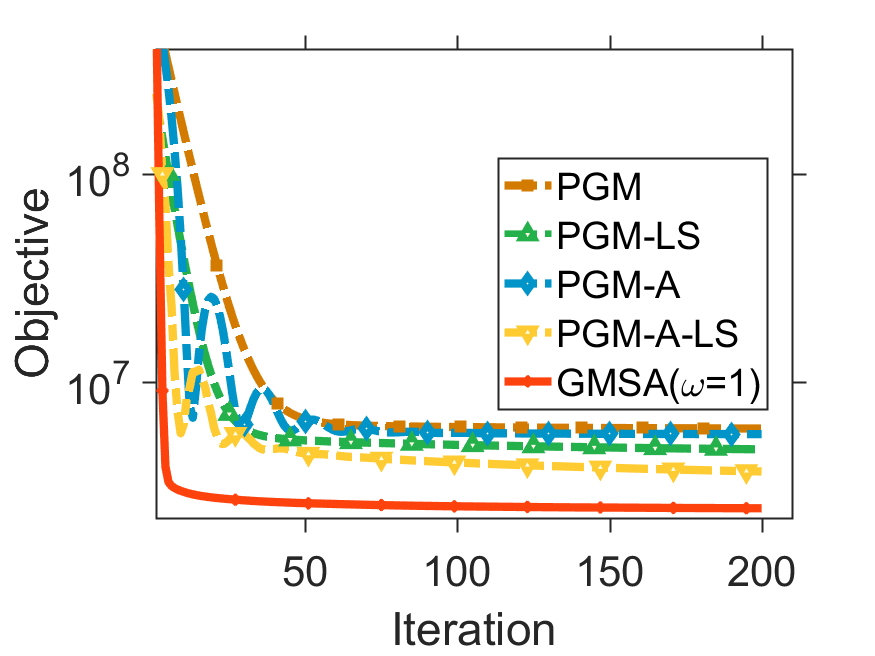}}
\subfloat[ \footnotesize $\lambda=500$]{\includegraphics[width=0.25\textwidth,height=0.14\textheight]{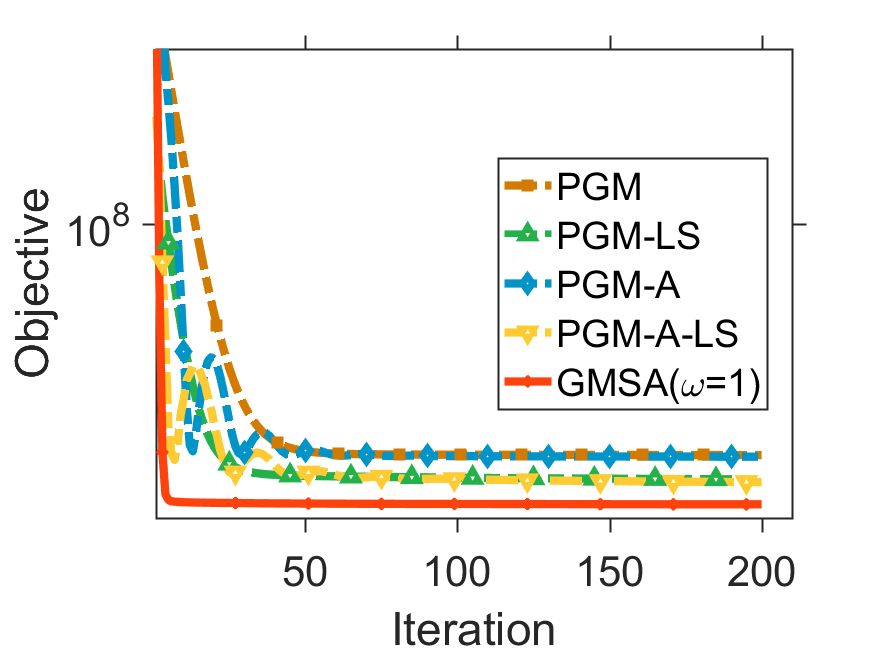}}
\subfloat[ \footnotesize $\lambda=5000$]{\includegraphics[width=0.25\textwidth,height=0.14\textheight]{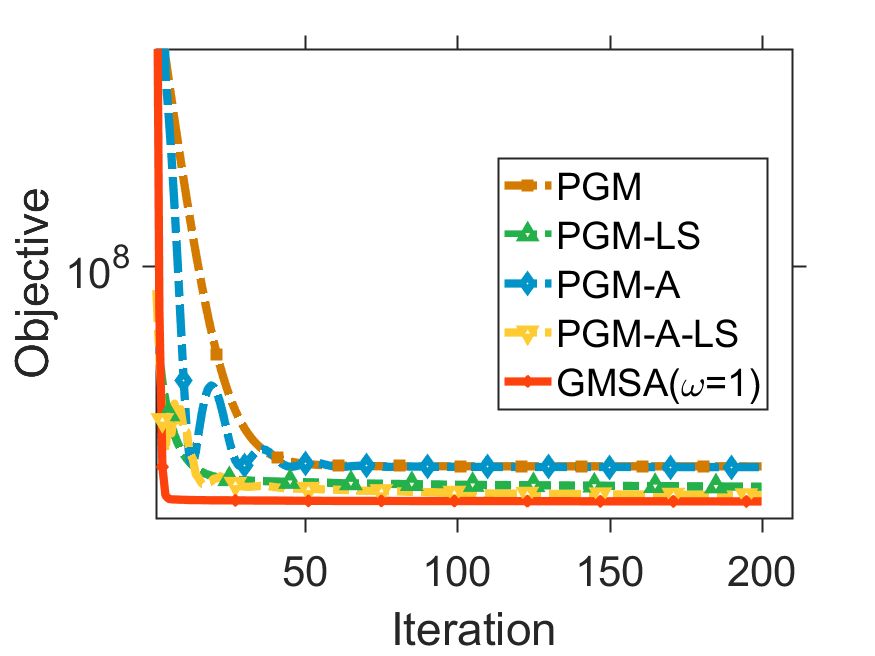}}
\subfloat[ \footnotesize $\lambda=50000$]{\includegraphics[width=0.25\textwidth,height=0.14\textheight]{type2_imagecameraman_lambda5000.png}}


\subfloat[ \footnotesize $\lambda=50$]{\includegraphics[width=0.25\textwidth,height=0.14\textheight]{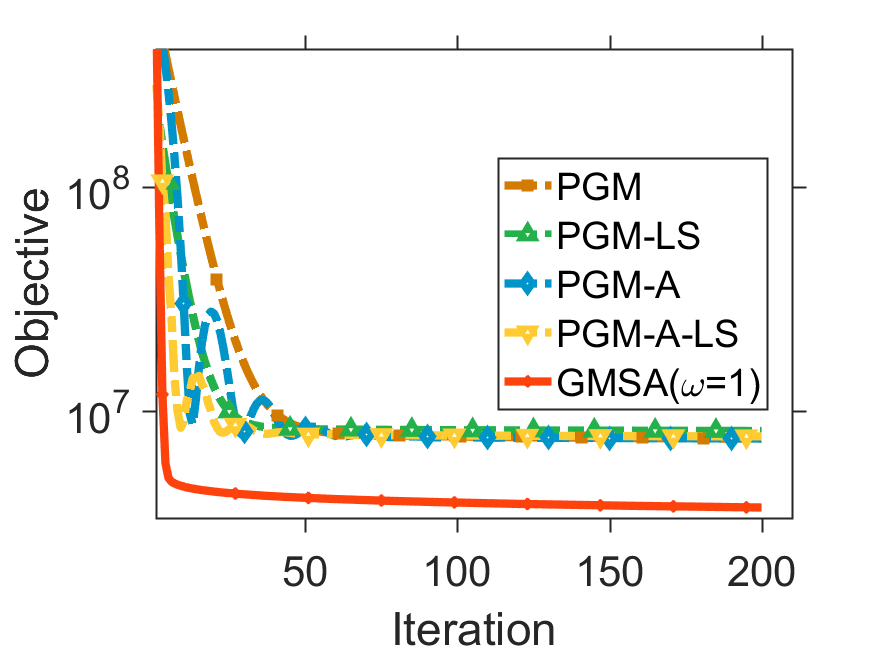}}
\subfloat[ \footnotesize $\lambda=500$]{\includegraphics[width=0.25\textwidth,height=0.14\textheight]{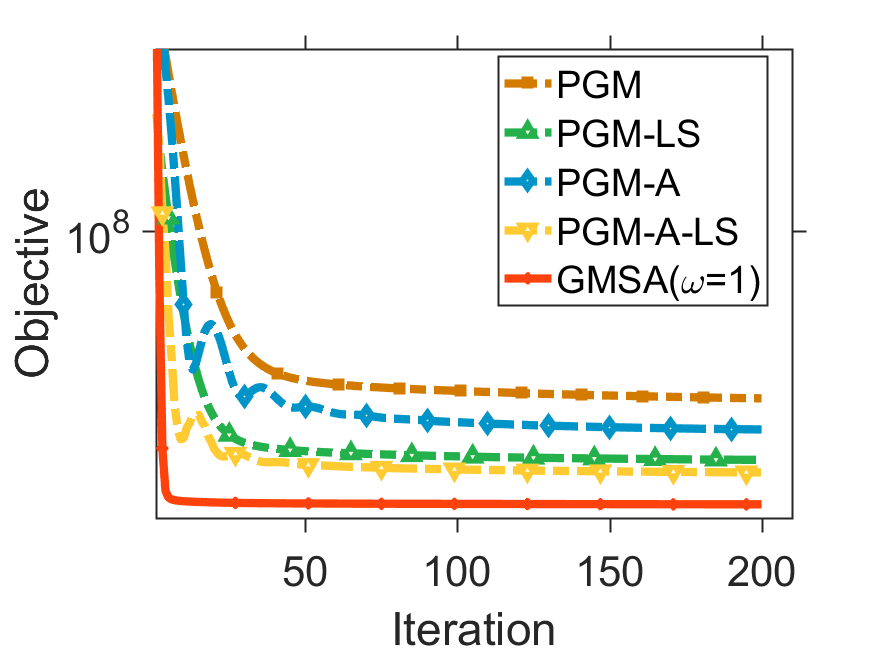}}
\subfloat[ \footnotesize $\lambda=5000$]{\includegraphics[width=0.25\textwidth,height=0.14\textheight]{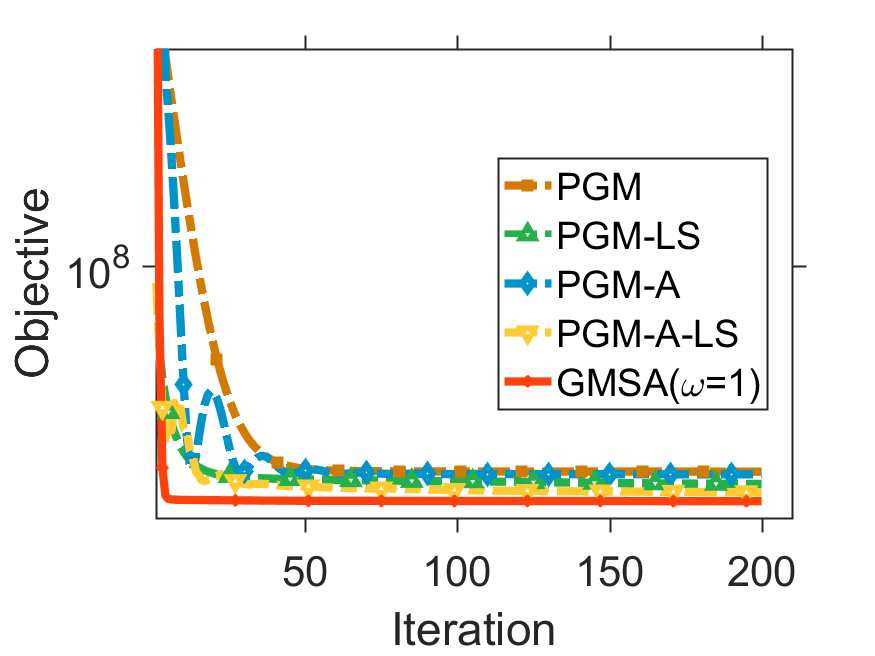}}
\subfloat[ \footnotesize $\lambda=50000$]{\includegraphics[width=0.25\textwidth,height=0.14\textheight]{type4_imagecameraman_lambda5000.png}}

\subfloat[ \footnotesize $\lambda=50$]{\includegraphics[width=0.25\textwidth,height=0.14\textheight]{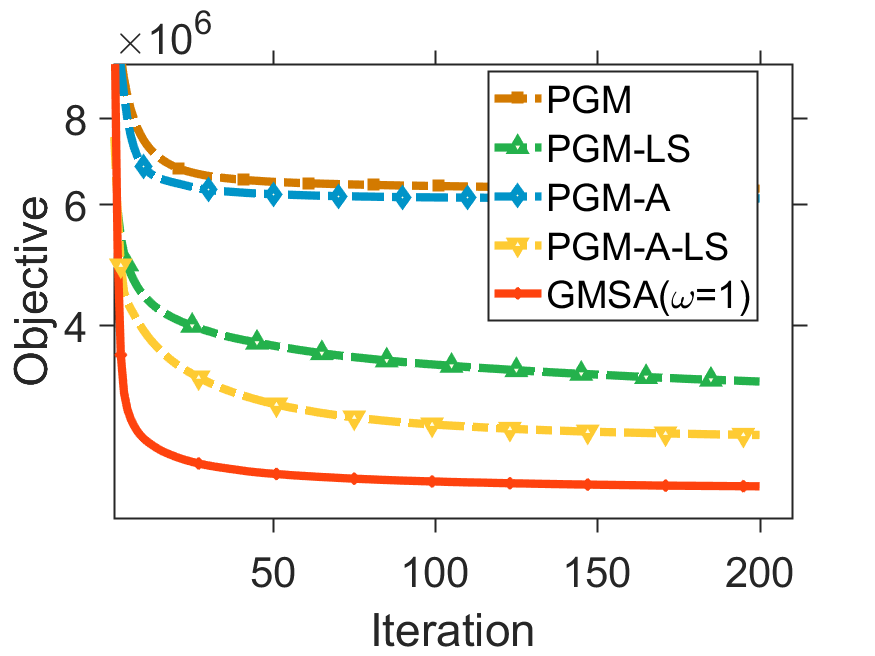}}
\subfloat[ \footnotesize $\lambda=500$]{\includegraphics[width=0.25\textwidth,height=0.14\textheight]{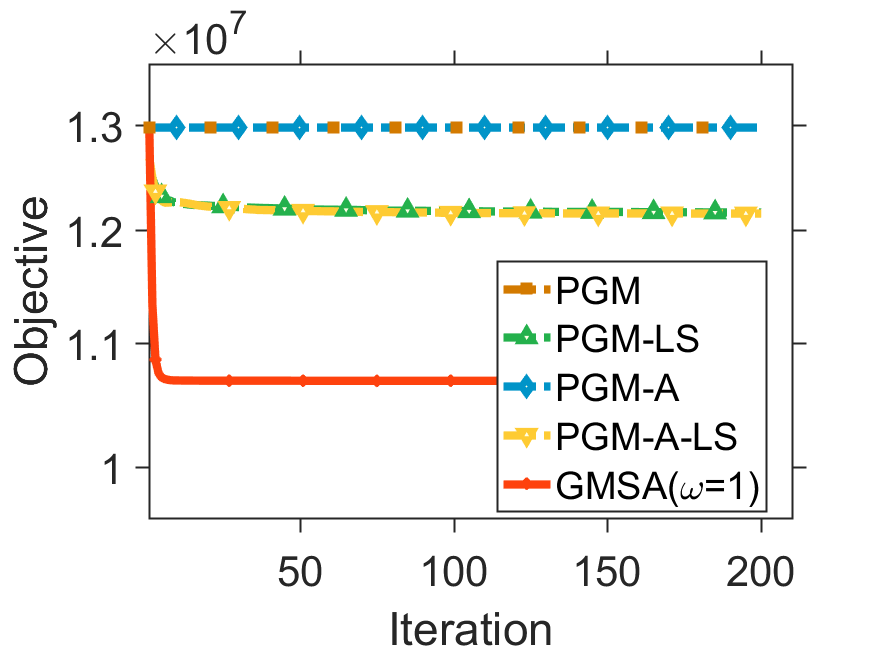}}
\subfloat[ \footnotesize $\lambda=5000$]{\includegraphics[width=0.25\textwidth,height=0.14\textheight]{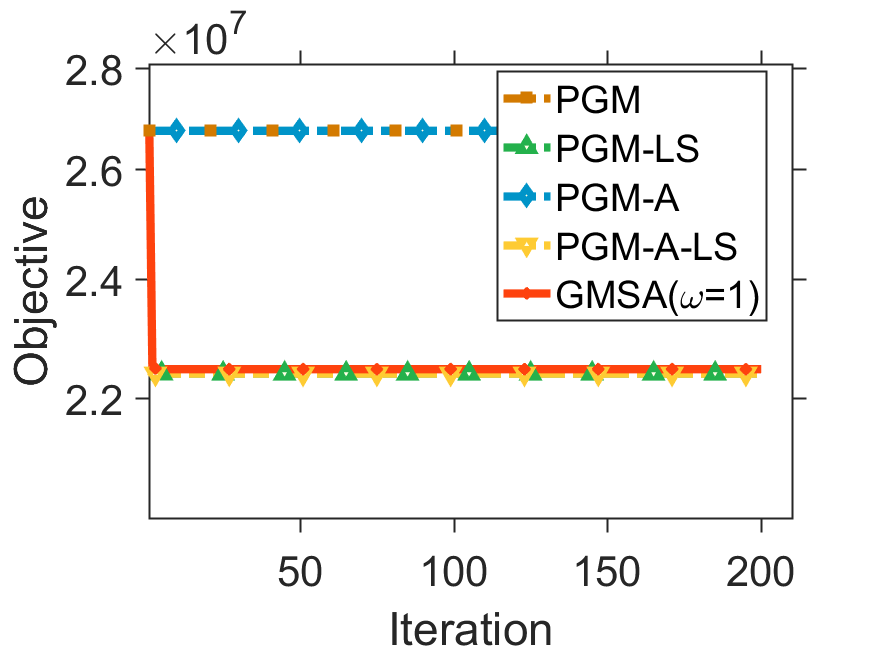}}
\subfloat[ \footnotesize $\lambda=50000$]{\includegraphics[width=0.25\textwidth,height=0.14\textheight]{type1_imagecameraman_lambda5000.png}}

\caption{Convergence behavior for solving (\ref{eq:card:sparse:coding}) with fixing $\bbb{W}$ for different $\lambda$ and initializations. Denoting $\tilde{\bbb{O}}$ as an arbitrary standard Gaussian random matrix of suitable size, we consider the following three initializations for $\bbb{H}$. First row: $\bbb{H}=0.1 \times \tilde{\bbb{O}}$. Second row: $\bbb{H}=10\times \tilde{\bbb{O}}$. Third row: $\bbb{H}$ is set to the output of the orthogonal matching pursuit.}
\label{fig:convergece:obj}
\end{figure*}

We compare the objective values for all methods by \emph{fixing} the variable $\bbb{W}$ to an over-complete DCT dictionary \cite{aharon2006img} and \emph{only} optimizing over $\bbb{H}$. We compare all methods with varying regularization parameter $\lambda$ and different initial points that are either generated by random Gaussian sampling or the Orthogonal Matching Pursuit (OMP) method \cite{tropp2007signal}. In Figure \ref{fig:convergece:obj}, we observe that GMSA converges rapidly in 10 iterations. Moreover, it often generates much better local optimal solutions than the compared methods.

\subsection{$\ell_1$ Norm Regularized Danzig Selectors}

Danzig selectors \cite{candes2007dantzig} can be formulated as the following optimization problem: $\min_{\bbb{x}}~\|\bbb{x}\|_1,~s.t.~\|\bbb{D}^{-1}\bbb{W}^T(\bbb{Wx} - \bbb{y})\|_{\infty}\leq \delta$, where $\bbb{W}\in\mathbb{R}^{m\times n},~\bbb{y}\in\mathbb{R}^{m},~\delta\in \mathbb{R}>0$, and $\bbb{D}\in \mathbb{R}^{n\times n}$ is the diagonal matrix whose diagonal entries are the norm of the columns of $\bbb{W}$. For the ease of discussion, we consider the following equivalent unconstrained optimization problem:
\beq \label{eq:danzig:2}
\min_{\bbb{x}}~\|\bbb{x}\|_1+ \lambda \|\bbb{Qx}-\bbb{s}\|_{\infty}
\eeq
\noi with $\lambda \propto \frac{1}{\delta}$, and $\bbb{Q}=\bbb{D}^{-1}\bbb{W}^T\bbb{W},~\bbb{s}=\bbb{D}^{-1}\bbb{W}^T\bbb{y}$.

We generate the design matrix $\bbb{W}\in \mathbb{R}^{m\times n}$ via sampling from a standard Gaussian distribution. The sparse original signal $\ddot{\bbb{x}}\in \mathbb{R}^{n\times 1}$ is generated via selecting a support set of size 20 uniformly at random and set them to arbitrary number sampled from standard Gaussian distribution. We set $\bbb{y}=\bbb{W}\ddot{\bbb{x}}$. We fix $n=1000$ and consider different choices for $\lambda$ and $m$.

\begin{figure*} [!t]
\centering

\subfloat[ \footnotesize $\lambda=0.05,~m=100$]{\includegraphics[width=0.25\textwidth,height=0.14\textheight]{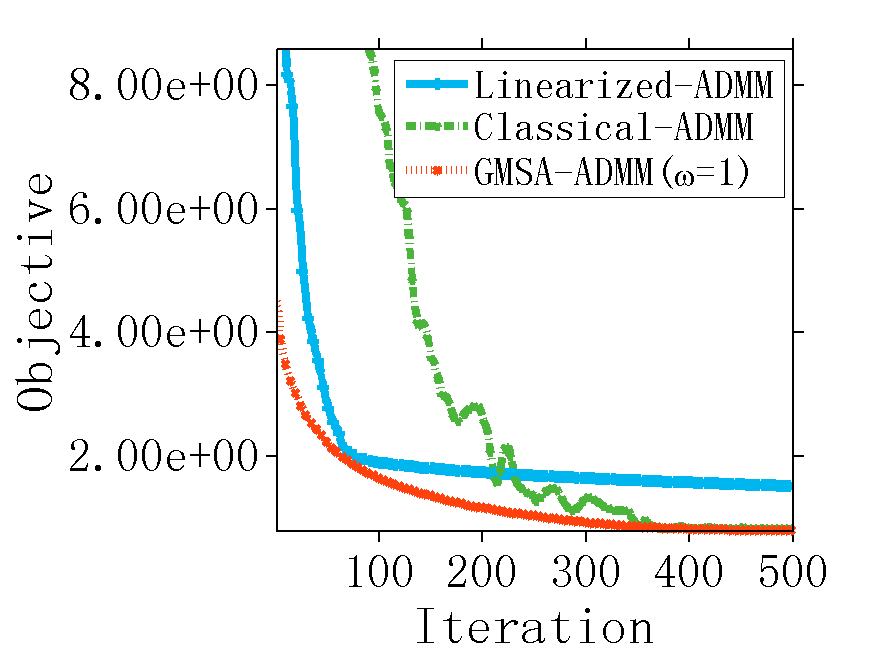}}
\subfloat[ \footnotesize $\lambda=0.05,~m=500$]{\includegraphics[width=0.25\textwidth,height=0.14\textheight]{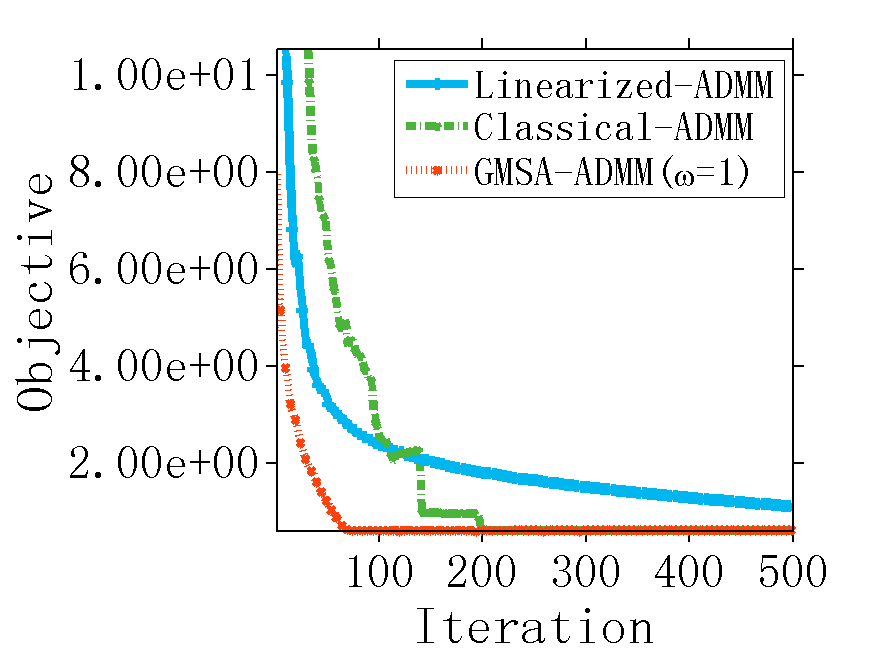}}
\subfloat[ \footnotesize $\lambda=0.05,~m=1000$]{\includegraphics[width=0.25\textwidth,height=0.14\textheight]{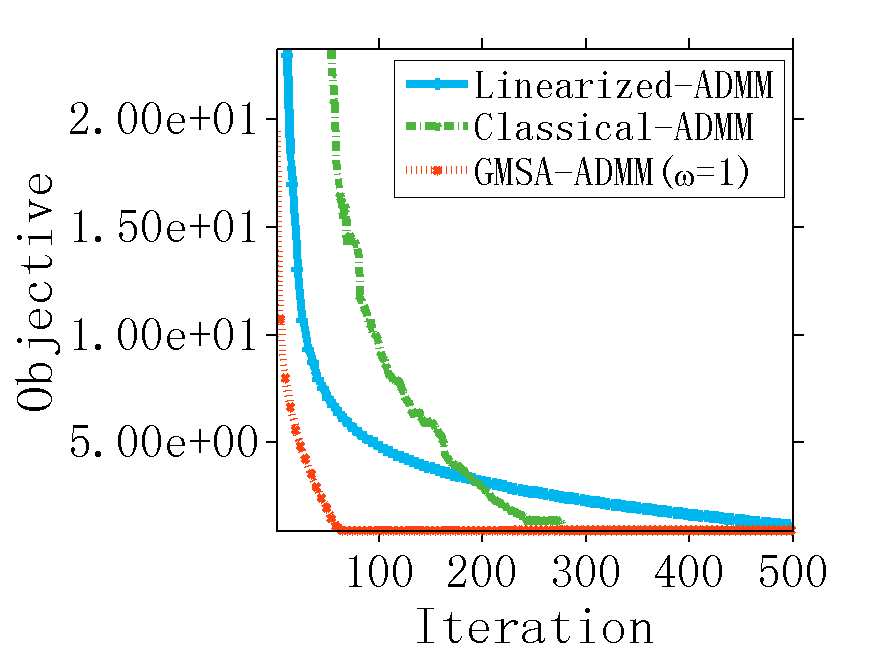}}
\subfloat[ \footnotesize $\lambda=0.05,~m=2000$]{\includegraphics[width=0.25\textwidth,height=0.14\textheight]{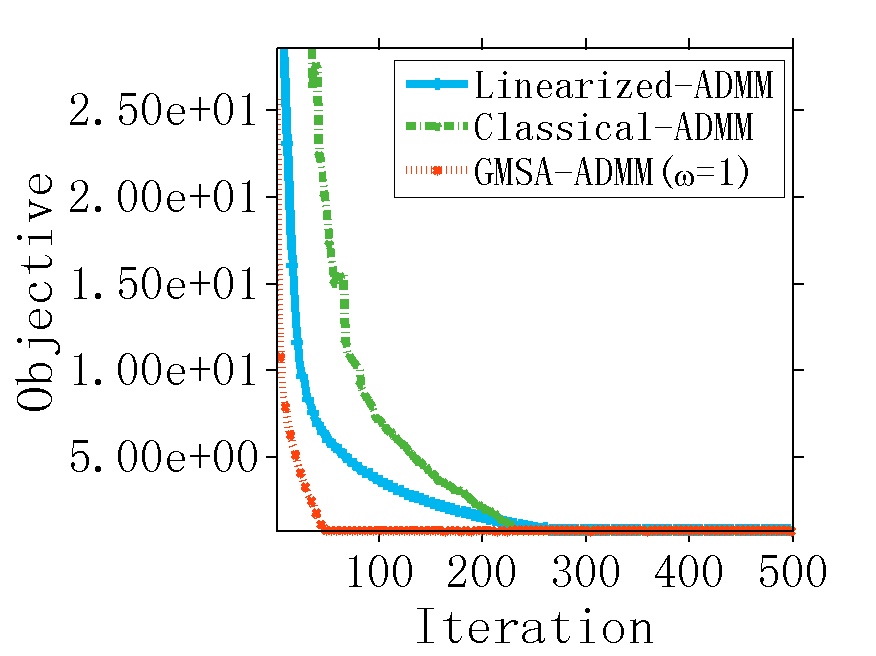}}\\

\subfloat[ \footnotesize $\lambda=0.5,m=100$]{\includegraphics[width=0.25\textwidth,height=0.14\textheight]{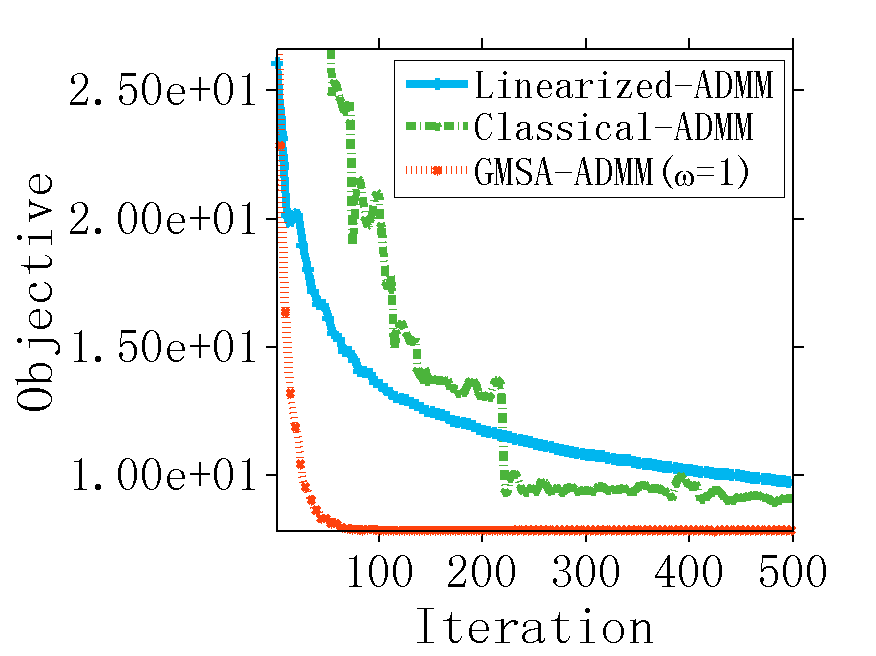}}
\subfloat[ \footnotesize $\lambda=0.5,m=500$]{\includegraphics[width=0.25\textwidth,height=0.14\textheight]{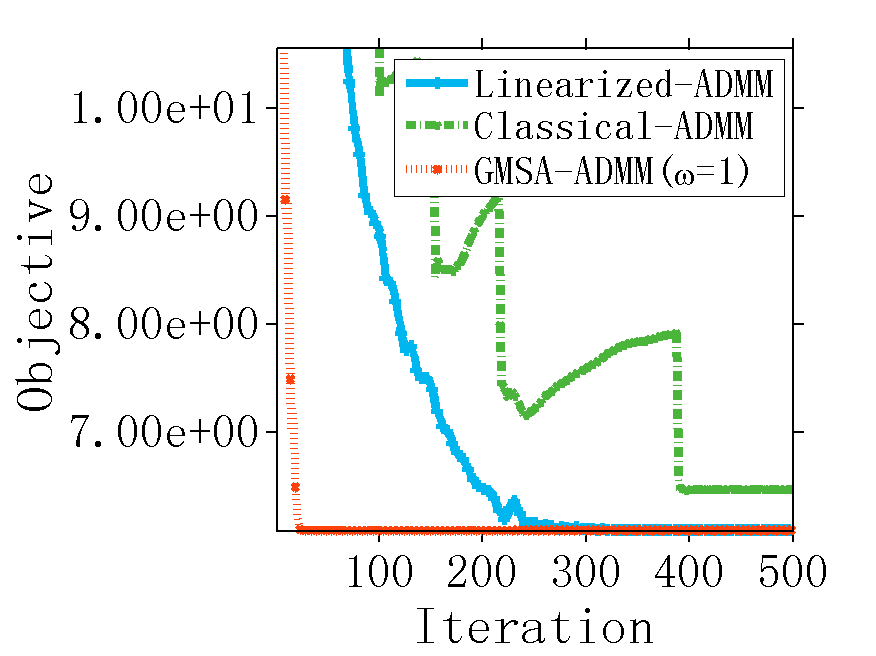}}
\subfloat[ \footnotesize $\lambda=0.5,m=1000$]{\includegraphics[width=0.25\textwidth,height=0.14\textheight]{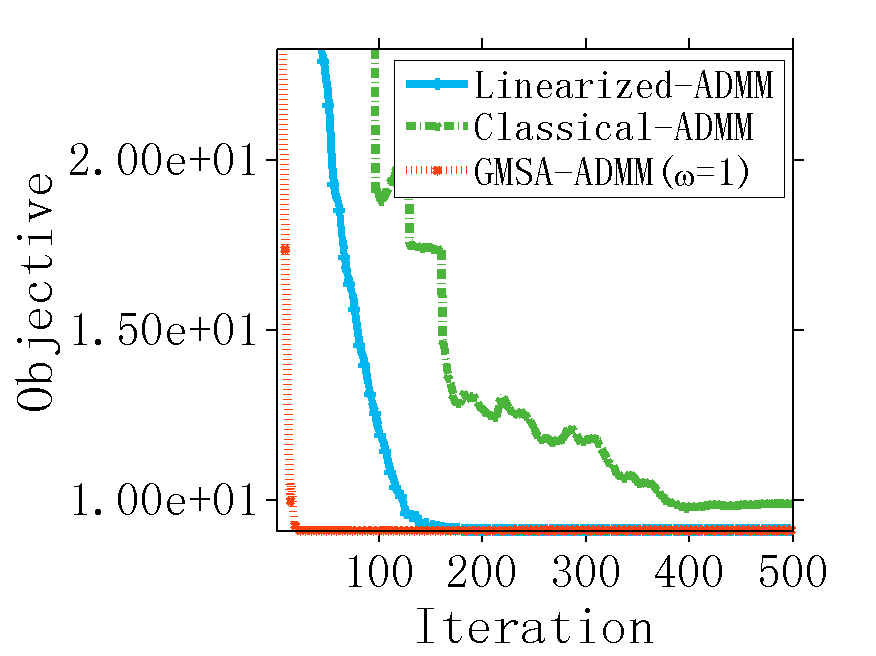}}
\subfloat[ \footnotesize $\lambda=0.5,m=2000$]{\includegraphics[width=0.25\textwidth,height=0.14\textheight]{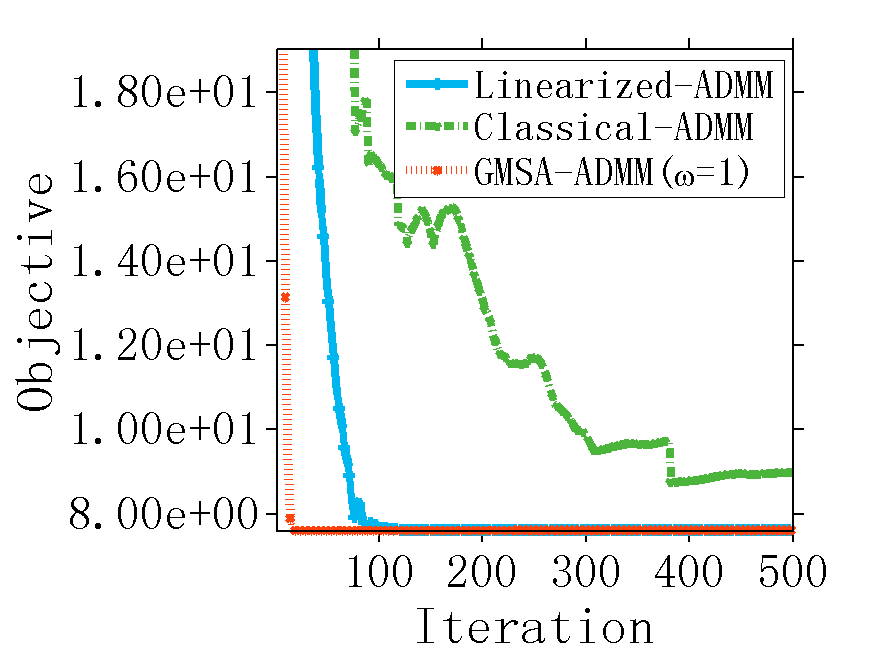}}\\

\subfloat[ \footnotesize $\lambda=5,~m=100$]{\includegraphics[width=0.25\textwidth,height=0.14\textheight]{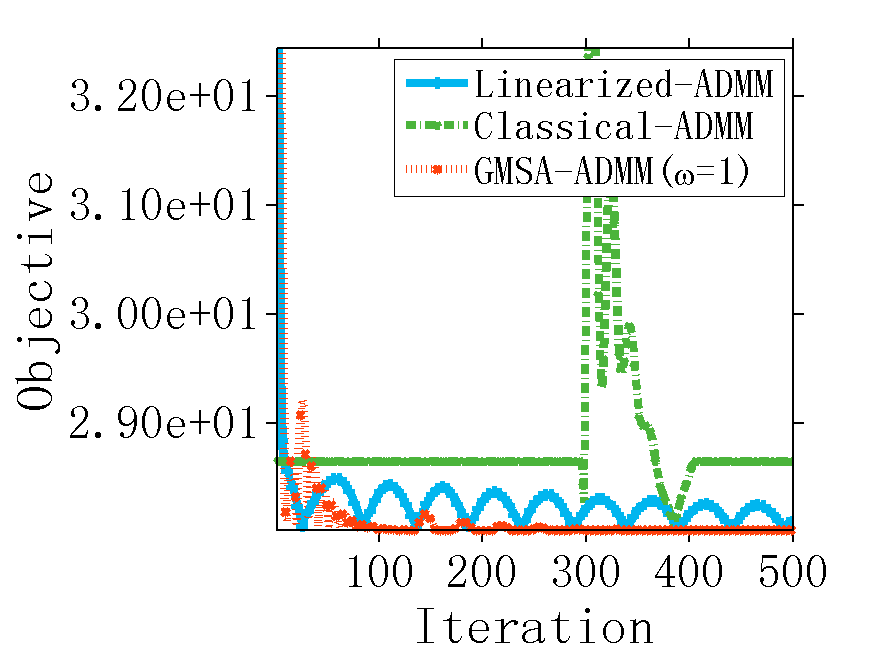}}
\subfloat[ \footnotesize $\lambda=5,~m=500$]{\includegraphics[width=0.25\textwidth,height=0.14\textheight]{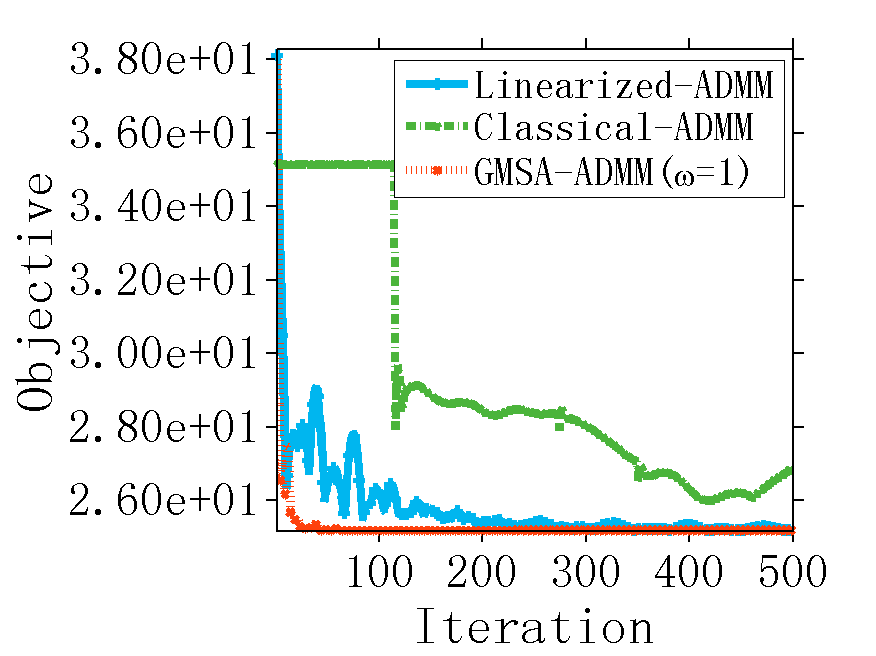}}
\subfloat[ \footnotesize $\lambda=5,~m=1000$]{\includegraphics[width=0.25\textwidth,height=0.14\textheight]{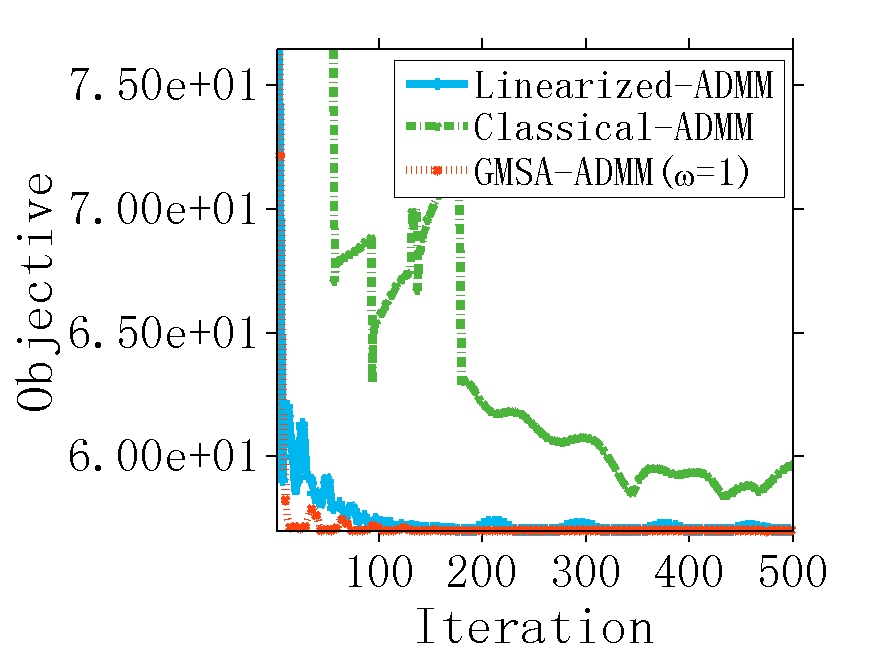}}
\subfloat[ \footnotesize $\lambda=5,~m=2000$]{\includegraphics[width=0.25\textwidth,height=0.14\textheight]{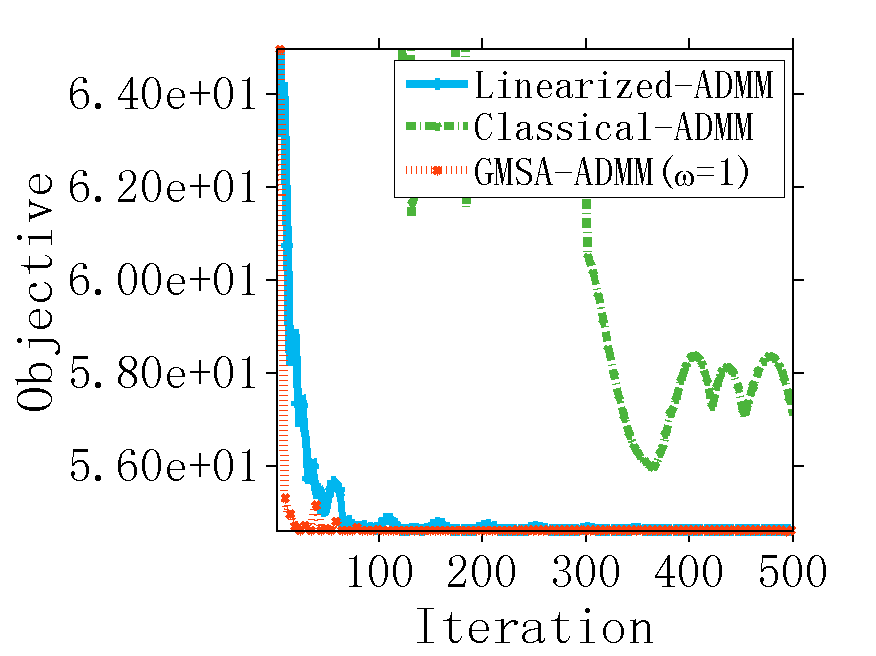}}

\caption{A comparison of linearized ADMM, classical ADMM, and GMSA-ADMM for solving the $\ell_1$ regularized Danzig selectors problem.}
\label{fig:convergece:danzig:obj}
\end{figure*}

We compare the proposed method GMSA-ADMM against linearized ADMM algorithm and classical ADMM algorithm. The penalty parameter $\beta$ is fixed to a constant with $\beta=1$. For linearized ADMM, we use the same splitting strategy as in \cite{WangY12}. For classical ADMM, we introduce addition two variables and rewrite (\ref{eq:danzig:2}) as: $\min_{\bbb{x},\bbb{y},\bbb{z}}~\|\bbb{z}\|_{\infty} + \lambda \|\bbb{x}\|_1,~s.t.~\bbb{x}=\bbb{y}, \bbb{Ay}-\bbb{b} = \bbb{z}$ to make sure that the smooth subproblem of the resulting augmented Lagrangian function is quadratic and can be solved by linear equations. For GMSA-ADMM, we do not solve the $\bbb{x}$-subproblem exactly using GMSA but solve it using one GMSA iteration. We demonstrate the objective values for the comparing methods. It can be been in Figure \ref{fig:convergece:danzig:obj} that our GMSA-ADMM significantly outperforms linearized ADMM and classical ADMM.

\section{Conclusions} 

This paper presents a new generalized matrix splitting algorithm for minimizing composite functions. We rigorously analyze its convergence behavior for convex problems and discuss its several importance extensions. Experimental results on nonnegative matrix factorization, $\ell_0$ norm regularized sparse coding, and $\ell_1$ norm regularized Danzig selector demonstrate that our methods achieve state-of-the-art performance.

%

\section*{Acknowledgments}
\noi This work was supported by the King Abdullah University of Science and Technology (KAUST) Office of Sponsored Research. This work was also supported by the NSF-China (61772570, 61472456, 61522115, 61628212).

\bibliographystyle{plain}
\bibliography{egbib}
\begin{IEEEbiography}[{\includegraphics[width=1in,height=1.25in]{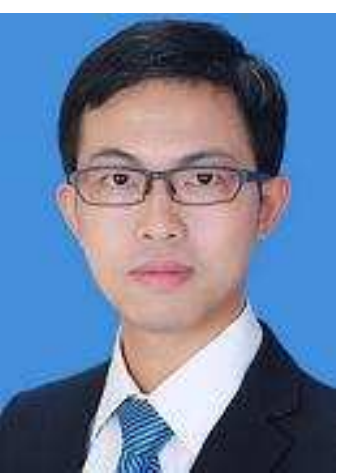}}]{Ganzhao Yuan} was born in Guangdong, China. He received his Ph.D. in School of Computer Science and Engineering, South China University of Technology in 2013. He is currently a research associate professor at School of Data and Computer Science in Sun Yat-sen University. His research interests primarily center around large-scale mathematical optimization and its applications in computer vision and machine learning. He has published technical papers in ICML, SIGKDD, AAAI, CVPR, VLDB, IEEE TPAMI, and ACM TODS.
\end{IEEEbiography}

\begin{IEEEbiography}[{\includegraphics[width=1in,height=1.25in]{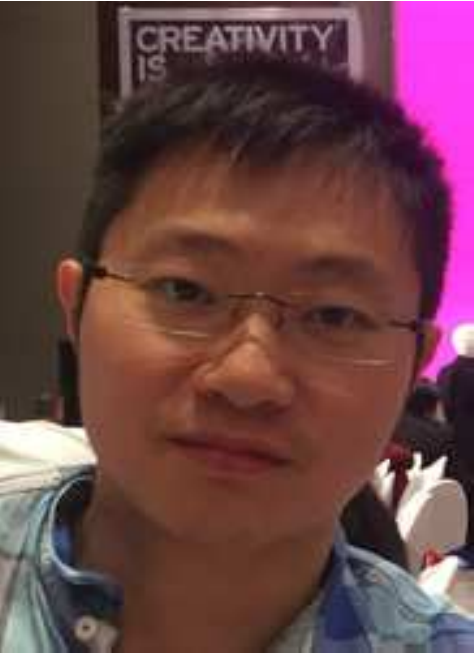}}]{Wei-Shi Zheng} was born in Guangdong, China. He is now a Professor at Sun Yat-sen University. He has now published more than 90 papers, including more than 60 publications in main journals (TPAMI, TIP, TNNLS, PR) and top conferences (ICCV, CVPR, IJCAI). His research interests include person/object association and activity understanding in visual surveillance. He has joined Microsoft Research Asia Young Faculty Visiting Programme. He is a recipient of Excellent Young Scientists Fund of the NSFC, and a recipient of Royal Society-Newton Advanced Fellowship. 
\end{IEEEbiography}

\begin{IEEEbiography}[{\includegraphics[width=1in,height=1.25in]{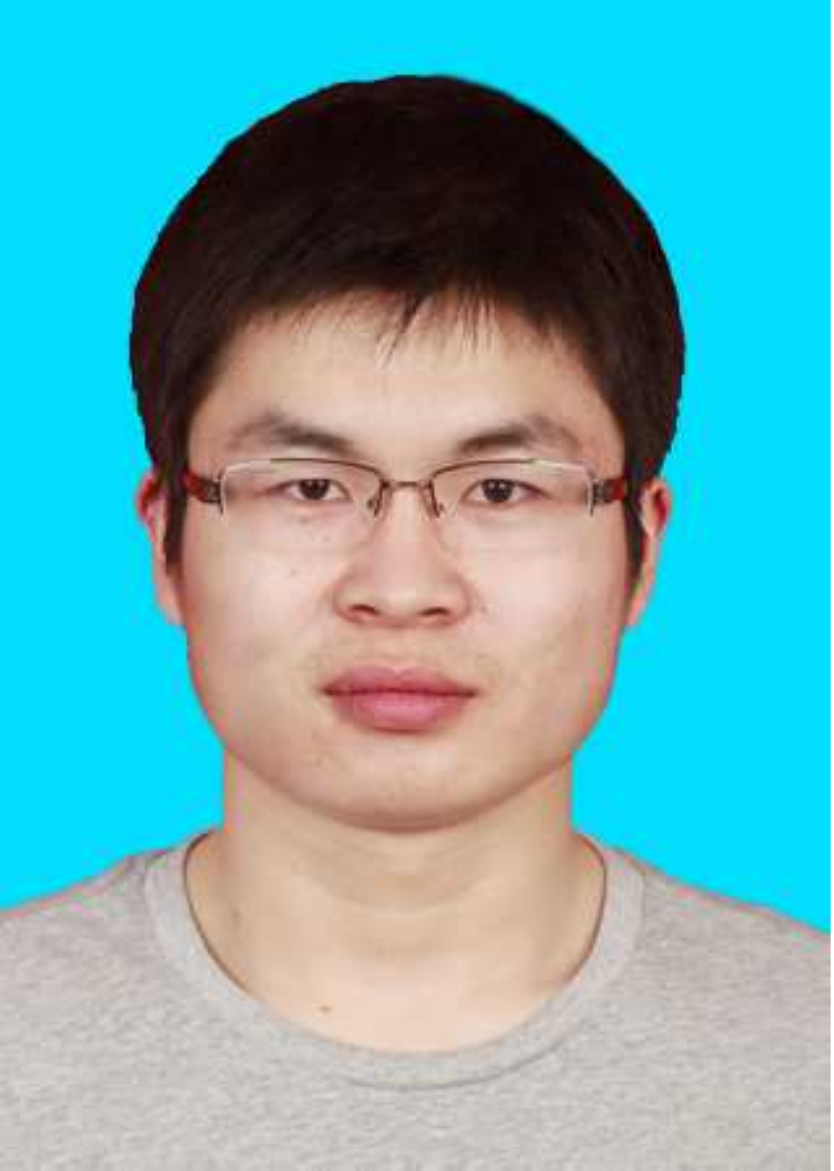}}]{Li Shen} was born in Hubei, China. He received his Ph.D. in School of Mathematics, South China University of Technology in 2017. He is currently a research scientist at Tencent AI Lab, Shenzhen. His research interests include algorithms for nonsmooth optimization, and their applications in statistical machine learning, game theory and reinforcement learning. He has published several papers in ICML and AAAI.
\end{IEEEbiography}

\begin{IEEEbiography}[{\includegraphics[width=1in,height=1.25in]{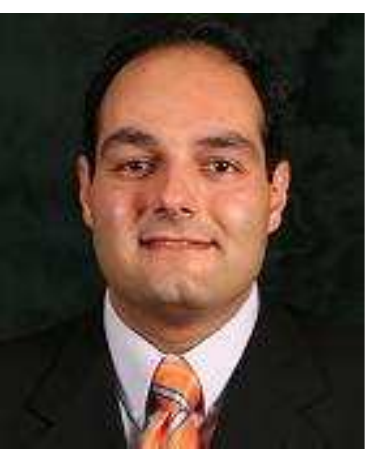}}]{Bernard Ghanem} was born in Betroumine, Lebanon. He received his Ph.D. in Electrical and Computer Engineering from the University of Illinois at Urbana-Champaign (UIUC) in 2010. He is currently an assistant professor at King Abdullah University of Science and Technology (KAUST), where he leads the Image and Video Understanding Lab (IVUL). His research interests focus on designing, implementing, and analyzing approaches to address computer vision problems (e.g. object tracking and action recognition/detection in video), especially at large-scale.
\end{IEEEbiography}

%
%
%
%
%
%
%

\end{document}